\documentclass{amsart}

\usepackage{amssymb}

\usepackage{graphicx}

\usepackage{amsmath, amsfonts, mathtools, tikz, tikz-cd}

\usetikzlibrary{chains}
\newtheorem{theorem}{Theorem}[section]

\newtheorem{prop}[theorem]{Proposition}

\theoremstyle{definition}
\newtheorem{definition}[theorem]{Definition}
\newtheorem{example}[theorem]{Example}

\theoremstyle{remark}

\usepackage{tikz-cd}
\usepackage{ytableau}

\numberwithin{equation}{section}

\newcommand{\dnode}[2][chj]{%
\node[#1] {};}

\newcommand{\dnodenj}[1]{%
\dnode[ch]{#1}}

\newcommand{\dnodebr}[1]{%
\node[chj] {};
}

\newcommand{\dydots}{%
\node[chj,draw=none,inner sep=1pt] {\dots};
}
\tikzset{node distance=2em, ch/.style={circle,draw,on chain,inner sep=2pt},chj/.style={ch,join},every path/.style={shorten >=4pt,shorten <=4pt},line width=1pt,baseline=-1ex}

\newcommand{\veps}{\varepsilon}

\allowdisplaybreaks
\def\ge{\mathfrak{g}}
\def\he{\mathfrak{h}}
\def\al{\alpha}
\def\cd{\cdots}
\def\cV{\mathcal{V}}

\def\C{\mathbb{C}}
\def\Z{\mathbb{Z}}
\def\L{\Lambda}

\begin{document}

\title{On $C_n^{(1)}$-geometric crystal and its Ultra-discretization}


%    author one information
%    author three information
\author{Erica S. Dinkins}
\address{Department of Mathematics,
North Carolina State University, Raleigh, NC 27695-8205, USA}
%\curraddr{}
\email{ecswain@ncsu.edu }
%\thanks{}

%    author two information
\author{Kailash C. Misra}
\address{Department of Mathematics,
North Carolina State University, Raleigh, NC 27695-8205, USA}
%\curraddr{}
\email{misra@ncsu.edu}
\thanks{KCM is partially supported by the Simons Foundation Grant \#636482.}

%    The 2010 edition of the Mathematics Subject Classification is
%    the current definitive version.
%\subjclass[2010] {Primary }

\date{}

\begin{abstract}
Let $\ge$ be an affine Lie algebra with index set $I = \{0, 1, 2, \cdots , n\}$ and $\ge^L$ be its Langlands dual. It is conjectured that  for each Dynkin node $i \in I \setminus \{0\}$ the affine Lie algebra $\ge$ has a positive geometric crystal whose ultra-discretization is isomorphic to the limit of certain coherent family of perfect crystals for the Langland dual $\ge^L$. In this paper we construct positive geometric crystals for 
$\cV(C_n^{(1)})$ in the level zero fundamental spin $C_n^{(1)}$- module $W(\varpi_n)$ for $n = 2, 3,4$ and show that its ultra-discretization is isomorphic to the limit $B^{n, \infty}$ of a coherent family $\{B^{n, l}\}_{l \geq 1}$ of perfect crystals for the Langland dual $D_n^{(2)}$ which proves the conjecture in these cases. 
\end{abstract}

\maketitle

%\newpage
\section{Introduction}

Representation theory of affine Lie algebras $\ge$ with index set $I = \{0, 1, 2, \cdots , n\}$, Cartan matrix $A = (a_{ij})_{i,j \in I}$ \cite{Kac} and their quantized universal enveloping algebra $U_q(\ge)$ \cite{HK} for generic $q$ (called quantum affine algebras) have led to many applications in other areas of mathematics and mathematical physics. Around 1990 Kashiwara and Lusztig independently introduced the concept of crystal bases for integrable representations of a quantum affine algebra (hence affine algebra) \cite{Kas1, Kas2, Lu}. By this theory we associate a combinatorial object $B(\lambda)$ called crystal associated to each highest weight integrable $U_q(\ge)$ - module $V(\lambda)$ for any dominant integral weight $\lambda$ of level $l = \lambda(c)$, where $c$ is the canonical central element in $\ge$. In an effort to give explicit realizations of the affine crystals $B(\lambda)$, Kang et al \cite{KMN1, KMN2} introduced the notion of perfect crystal $B^{k, l}$ for each positive integer $l$ and nonzero Dynkin index $k$ for a quantum affine algebra $U_q(\ge)$ and showed that the affine crystal $B(\lambda)$ can be realized as a suitable sub-crystal of the crystal of the semi-infinite tensor product $B^{k, l} \otimes B^{k, l}  \otimes B^{k, l} \otimes \cdots \cdots$. This realization is called the path realization of the crystal $B(\lambda)$. It turned out that the perfect crystal $B^{k, l}$ is a crystal for certain finite dimensional $U_q(\ge)$-module $W(l\varpi_k)$, called Kirillov-Reshetikhin module (KR-module) \cite{KR} where 
$\varpi_k$ is the $k^{th}$ fundamental weight of level zero \cite{Kas3}. When $\{B^{k,l}\}_{l\geq 1}$ is a coherent family of perfect crystals \cite{KKM}
we denote its limit by $B^\infty (\varpi_k)$ (or just $B^{k,\infty}$).

On the other hand the notion of geometric crystals is defined in \cite{BK} for reductive algebraic groups and extended to Kac-Moody groups in \cite{N}. Let $\ge'$ be the derived Lie algebra of the affine Lie algebra $\ge$ and let  $G$ be the Kac-Moody group associated  with $\ge'$(\cite{KP, PK}).
Let $U_{\alpha}:=\exp\ge_{\alpha}$ $(\alpha\in \Delta^{re})$ be the one-parameter subgroup of $G$ where $\Delta^{re}$ denote the set of real roots. The group $G$ is generated by 
$U_{\alpha}$ $(\alpha\in \Delta^{re})$. Let $U^{\pm}:=\langle U_{\pm\alpha}|\alpha\in\Delta^{re}_+\rangle$ be the subgroup generated by $U_{\pm\alpha}$
($\al\in \Delta^{re}_+$).
For any $i\in I$, there exists a unique homomorphism;
$\phi_i:SL_2(\C)\rightarrow G$ such that
\[
\hspace{-2pt}\phi_i\left(
\left(
\begin{array}{cc}
c&0\\
0&c^{-1}
\end{array}
\right)\right)=c^{\check{\alpha}_i},\,
\phi_i\left(
\left(
\begin{array}{cc}
1&t\\
0&1
\end{array}
\right)\right)=\exp(t e_i),\,
 \phi_i\left(
\left(
\begin{array}{cc}
1&0\\
t&1
\end{array}
\right)\right)=\exp(t f_i).
\]
where $c\in\C^\times$ and $t\in\C$.
Set $\check{\alpha}_i(c):=c^{\check{\alpha}_i}$,
$x_i(t):=\exp{(t e_i)}$, $y_i(t):=\exp{(t f_i)}$, 
$G_i:=\phi_i(SL_2(\C))$,
$H_i:=\phi_i(\{{\rm diag}(c,c^{-1})\mid 
c\in\C \setminus \{0\}\})$. 
%and 
%$N_i:=N_{G_i}(H_i)$
Let $H$  be the subgroup of $G$ generated by $H_i$'s
with the Lie algebra $\he$. 
Then $H$ is called a  maximal torus in $G$, and 
$B^{\pm}=U^{\pm}H$ are the Borel subgroups of $G$.
%We have the isomorphism
%$\phi:W \longrightarrow N/H$ defined by $\phi(s_i)=N_iH/H$.
The element $\bar{s}_i:=x_i(-1)y_i(1)x_i(-1) \in N_G(H)$  is a representative of the simple reflection
$s_i\in W=N_G(H)/H$. 
The geometric crystal for the affine Lie algebra $\ge$ is defined as follows.
\begin{definition}\label{geometric}(\cite{BK},\cite{N}) 
The geometric crystal for the affine Lie algebra $\ge$ is a quadruple $\cV(\ge)=(X, \{e_i\}_{i \in I}, \{\gamma_i\}_{i \in I},$ 
$\{\veps_i\}_{i\in I})$, 
where $X$ is an ind-variety,  $e_i:\C^\times\times
X\longrightarrow X$ $((c,x)\mapsto e^c_i(x))$
are rational $\C^\times$-actions and  
$\gamma_i,\veps_i:X\longrightarrow 
\C$ $(i\in I)$ are rational functions satisfying the following:
\begin{enumerate}
    \item $\{1\} \times X\cap \text{dom}\{e_i\}$ is open dense in $\{1\}\times X$ for any $i\in I$
    \item $\gamma_j(e_i^c(x))=c^{a_{ij}}\gamma_j(x)$
    \item $\varepsilon_i(e_i^c(x))=c^{-1}\varepsilon_i(x)$ and $\varepsilon_i(e_j^c(x))=\varepsilon_i(x)$ if $a_{ij}=a_{ji}=0$
    \item The $e_i$ satisfy:
    \begin{align*}
    \begin{cases}
        &e_i^{c_1}e_j^{c_2}=e_{j}^{c_2}e_i^{c_1}\,\,\,\,\text{if }a_{ij}=a_{ji}=0,\\
        &e_i^{c_1}e_j^{c_1c_2}e_i^{c_2}=e_{j}^{c_2}e_i^{c_1c_2}e_j^{c_1}\,\,\,\,\text{if }a_{ij}=a_{ji}=-1,\\
        &e_i^{c_1}e_j^{c_1^2c_2}e_i^{c_1c_2}e_j^{c_2}=e_{j}^{c_2}e_i^{c_1c_2}e_j^{c_1^2c_{2}}e_i^{c_1}\,\,\,\,\text{if }a_{ij}=-2, \,\,a_{ji}=-1,\\
        &e_i^{c_1}e_j^{c_1^3c_2}e_i^{c_1^2c_2}e_j^{c_1^3c_2^2}e_i^{c_1c_2}e_j^{c_2}=e_{j}^{c_2}e_i^{c_1c_2}e_j^{c_1^3c_{2}^2}e_i^{c_1^2c_2}e_j^{c_1^3c_2}e_i^{c_1}\,\,\,\,\text{if }a_{ij}=-3, \,\,a_{ji}=-1.
    \end{cases}    
    \end{align*}
\end{enumerate}
\end{definition}
The flag variety $X:=G/{B^+}$ is an ind-variety. For $w \in W$, the Schubert cell $X_w$ associated with $w$ has a natural geometric crystal structure \cite{BK, N}. Let $w=s_{i_1}, s_{i_2}, \cdots ,s_{i_l}$ be a reduced expression. For ${\bf i}:=(i_1, i_2, \cdots ,i_l)$, set 
\begin{equation*}
B_{\bf i}^-
:=\{Y_{\bf i}(c_1, c_2, \cdots ,c_l)
:=Y_{i_1}(c_1), Y_{i_2}(c_2),\cdots , Y_{i_l}(c_l)
\,\vert\, c_1, c_2, \cdots ,c_l\in\C^\times\}\subset B^-,
%(Y_i(c):=y_i(\frac{1}{c})\al_i^\vee(c)),
\label{bw1}
\end{equation*}
where $Y_i(c):=y_i(\frac{1}{c})\al^\vee_i(c)$. Then  we have the following result.
\begin{theorem}\label{schubert}\cite{N}
If $I = \{i_1, i_2, \cdots ,i_l\}$, then the set $B_{\bf i}^-$ with the explicit actions of \; $e^c_i$, $\veps_i$, and $\gamma_i$ , for $i \in I, c \in \C^\times$ given by:
\begin{eqnarray}
&& e_i^c(Y_{\bf i}(c_1,\cdots,c_l))
=Y_{\bf i}({\mathcal C}_1,\cdots,{\mathcal C}_l)), \notag\\
&&\text{where} \notag\\
&&{\mathcal C}_j:=
c_j\cdot \frac{\displaystyle \sum_{1\leq m\leq j,i_m=i}
 \frac{c}
{c_1^{a_{i_1,i}}\cdots c_{m-1}^{a_{i_{m-1},i}}c_m}
+\sum_{j< m\leq k,i_m=i} \frac{1}
{c_1^{a_{i_1,i}}\cdots c_{m-1}^{a_{i_{m-1},i}}c_m}}
{\displaystyle\sum_{1\leq m<j,i_m=i} 
 \frac{c}{c_1^{a_{i_1,i}}\cdots c_{m-1}^{a_{i_{m-1},i}}c_m}+ 
\mathop\sum_{j\leq m\leq k,i_m=i}  \frac{1}
{c_1^{a_{i_1,i}}\cdots c_{m-1}^{a_{i_{m-1},i}}c_m}},\\
%\label{eici}\\
&& \veps_i(Y_{\bf i}(c_1,\cd,c_l))=
\sum_{1\leq m\leq l,i_m=i} \frac{1}
{c_1^{a_{i_1,i}}\cdots c_{m-1}^{a_{i_{m-1},i}}c_m},\\
%\label{vep-i}\\
&&\gamma_i(Y_{\bf i}(c_1,\cdots,c_l))
=c_1^{a_{i_1,i}}\cdots c_l^{a_{i_l,i}},
%\label{gamma-i}
\end{eqnarray}
is a geometric crystal isomorphic to $X_w$.
\end{theorem}

Note that if $J:=\{i_1,\cd,i_k\}\subsetneq I$ and  $\ge_J\subsetneq \ge$ 
be the corresponding subalgebra, then
by arguments similar to \cite[4.3]{N} we can define the $\ge_J$-geometric crystal 
structure on $B^-_{\bf i}$.
The geometric crystal $\cV(\ge)=(X, \{e_i\}_{i \in I}, 
\{\gamma_i\}_{i \in I}, \{\veps_i\}_{i\in I})$ is said to be positive if it has a 
positive structure \cite{BK, N}. 
Roughly speaking this means that each of the rational maps 
$e^c_i$, $\veps_i$  and $\gamma_i$ are
given by ratio of polynomial functions with positive coefficients. 
For example, $B_{\bf i}^-$ is a positive geometric crystal.
%For a dominant weight  $\lambda \in P^+$ of level $l$,
%Kashiwara defined the crystal base $(L(\lambda), B(\lambda))$
%\cite{Kas1} (also see \cite{Lu})  for the integrable highest weight $\ge$-module $V(\lambda)$. As shown in \cite{KMN1}, the crystal
%$B(\lambda)$ can be realized as a set of paths in the semi-infinite tensor product $ \cdots \otimes B^l \otimes B^l \otimes B^l$ where
%$B^l$ is a perfect crystal of level $l$. This is called the path realization of the crystal $B(\lambda)$. When a family of perfect crystals 
%$\{B^l\}_{l\geq 1}$ are coherent \cite{KKM}, it has a limit $B^\infty$. 
The positive geometric crystals are related to Kashiwara crystals via the 
ultra-discretization functor $\mathcal{UD}$ \cite{BK, N} which transforms positive rational 
functions to piecewise-linear functions by the simple correspondence:
$$
x \times y \longmapsto x+y, \qquad \frac{x}{y} \longmapsto x - y, 
\qquad x + y \longmapsto {\rm max}\{x, y\}.
$$
It was conjectured in \cite{KNO} that for each affine Lie algebra $\ge$ and 
each Dynkin index $k \in I \setminus \{0\}$, there exists a positive geometric crystal
$\cV(\ge)=(X, \{e_i\}_{i \in I}, \{\gamma_i\}_{i \in I}, \\
\{\veps_i\}_{i\in I})$ whose ultra-discretization $\mathcal{UD}(\cV)$ is isomorphic 
to the limit $B^{k,\infty}$ of a coherent family of perfect crystals for the Langlands dual $\ge^L$. If $\ge$ is simply laced then the Langland dual is 
$\ge$ itself. So far this conjecture has been proved for the Dynkin index $k = 1$ and $\ge = A_n^{(1)}, 
B_n^{(1)}, C_n^{(1)}, D_n^{(1)}, A_{2n-1}^{(2)}, A_{2n}^{(2)},$
$D_{n+1}^{(2)}$ \cite{KNO}, $\ge = D_4^{(3)}$ \cite{IMN}, $\ge = G_2^{(1)}$ \cite{N4}. Also this conjecture has been shown to hold for $\ge = A_n^{(1)}$  at any node $1 \leq k \leq n$ \cite{MN1, MN2}, and for $\ge = D_n^{(1)}$ at the spin node $n$, $n = 5, 6$ \cite{IMP, MP1, MP2}.

In this paper we prove the conjecture in \cite{KNO} for $\ge = C_n^{(1)}$ at the  Dynkin spin node $n$ for $n = 2, 3, 4$. Note that the Lagland dual of $\ge =  C_n^{(1)}$ is the affine algebra $\ge^L = D_{n+1}^{(2)}$. In Section 2, we recall some necessary facts for the affine Lie algebra $C_n^{(1)}$ and show that the level zero fundamental module $W(\varpi_n)$ is not multiplicity free. In the following three sections we explicitly construct the positive geometric crystals and prove the conjecture in \cite{KNO} for the cases $n = 2, 3$ and $4$. As can be observed, since the fundamental representation $W({\varpi}_n)$ is not multiplicity free, the calculations become more and more involved as we go up on the rank.

\section{The affine Lie algebra $C_n^{(1)}$}

From now we denote $\ge$ to be the affine Lie algebra $C_n^{(1)}$ over the field of complex numbers $\mathbb{C}$ with index set $I = \{0,1,...,n\}$, Cartan matrix $A = (a_{ij})_{i,j \in I}$ where $a_{ii} = 2, a_{j,j + 1} = -1 = a_{j+1,j}$, for  $j = 1,...,n-1,\; a_{10} =-2=a_{n-1,n}, \; a_{01} = a_{n,n-1} = -1$, and $a_{ij} = 0$ otherwise and  Dynkin diagram \cite{Kac}:
\begin{center}
\begin{tikzpicture}[start chain]
\dnodenj{1}
\dnodenj{2}
\dydots
\dnode{2}
\dnodenj{1}
\path (chain-1) -- node{\(\Rightarrow\)} (chain-2);
\path (chain-4) -- node{\(\Leftarrow\)} (chain-5);
\end{tikzpicture}
\end{center}
Note that the subalgebra $\ge_i$ associated with the Cartan submatrix $(a_{ij})_{i,j \in {I\setminus\{i\}}}$ is isomorphic to the simple Lie algebra $C_n$ for $i = 0,$ and  $n$.
Let $\{\alpha_0, \alpha_1, ..., \alpha_n\}, \\ \{\check{\alpha_0}, \check{\alpha_1}, ..., \check{\alpha_n}\}$ and $\{\Lambda_0, \Lambda_1, ..., \Lambda_n\}$ denote the set of simple roots, simple coroots and fundamental weights, respectively.
Then $c=\check{\alpha_0}+\check{\alpha_1}+...+\check{\alpha_n}$ and $\delta = \al_0 +2\al_1+...+2\al_{n-1}+\al_n$ are the central element and null root respectively. The sets $P_{cl} = \oplus_{j=0}^n \Z\L_j$ and $P = P_{cl}\oplus\Z\delta$ are called classical weight lattice and weight lattice respectively. The subalgebra $\he = {\rm span}\{\check{\alpha_0}, \check{\alpha_1}, ..., \check{\alpha_n}, d\}$ is the Cartan subalgebra of $\ge$ where $d$ denotes the degree derivation. The Weyl group $W$ for $\ge$ is generated by the simple reflections $\{s_0, s_1, \cdots , s_n\}$ where $s_i(\lambda)=\lambda-\lambda(\check{\alpha}_i)\alpha_i$ for all $\lambda \in P$. 

\subsection{Weight multiplicities of the $C_n$-module $V(\Lambda_n)$}
Let us denote the simple roots and the fundamental weights for the simple Lie algebra $\ge_0 = C_n$ by the same notations $\{\alpha_1, ..., \alpha_n\}$, and $\{\Lambda_1, ..., \Lambda_n\}$ respectively. Let $\{\epsilon_1, \cdots , \epsilon_n\}$ denote the orthonormal base of the dual of the Cartan subalgebra of $C_n$. Then $\alpha_i = \epsilon_i - \epsilon_{i+1}, 1 \leq i \leq n-1$ and $\alpha_n = 2\epsilon_n$. Consider the fundamental $C_n$-module $V(\Lambda_n)$. We use the associated crystal $B(\Lambda_n)$ \cite{KN} to determine the multiplicities of weights in $V(\Lambda_n)$. Let $S = \{i, \overline{i} \mid 1\leq i \leq n\}$ with the order $1\preceq 2\prec...\prec n\prec\overline{n}\prec...\preceq \overline{1}$. Then the elements in $B(\Lambda_n)$ are given by
\begin{equation}\label{crystal-B}
B(\Lambda_n) = 
\small{\left\{(i_1, \cdots , i_n) \mid
 \begin{tabular}{c} 
$(1) 1 \prec i_1 \prec \cdots \prec i_n \prec \overline{1},$ $i_j \in S$ \\
$(2) {\rm if} \ i_k = p \  {\rm and}  \  i_l = \overline{p},  \  {\rm then} \ k+(n-l+1) \leq p $
\end{tabular}\right\}
}
\end{equation}
Here the weight $wt(i) = \epsilon_i, \  wt(\overline{i}) = - \epsilon_i$ for $1\leq i \leq n$ and $wt(i_1, \cdots , i_n) = \sum_{j=1}^n wt(i_j)$.
For $v = (i_1, \cdots , i_n)  \in B(\Lambda_n)$, we say that the entry $i_k = p$ is a {\it pair} if there is another entry $i_l = \overline{p}$, otherwise we say
the entry $i_k = p$ is a {\it single}. We say that a $p$-pair $(i_k, i_l) = (p, \overline{p})$ occurring in $v$ is admissible if $k - l +n+1 \leq p$. Also note that if $(i_k, i_l) = (p, \overline{p})$ then $$wt(v) =  \sum_{\begin{array}{c}j=1\\ j \neq k, l\end{array}}^nwt(i_j)$$
 since $wt(i_k) + wt(i_l) = \epsilon_p -  \epsilon_p = 0$. Of course in $v \in B(\Lambda_n)$ there can be more than one {\it pair} entries. Let $\mu$ be a weight of $V(\Lambda_n)$.
 Then $\mu = \Lambda_n - \sum_{i=1}^n m_i\alpha_i$ for some $m_i \in \mathbb{Z}_{\geq 0}$ and ${\rm dim}(V(\Lambda_n)_\mu = 
 |B(\Lambda_n)_\mu|$. By definition any two elements $v = (i_1, \cdots , i_n) , v' = (i'_1 \cdots , i'_n) \in B(\Lambda_n)_\mu$ contain the same {\it single} entries since each such entry contributes either a distinct $\epsilon_p$ or 
a distinct $- \epsilon_p$ to the weights of $v$ and $v'$. Thus if $v = (i_1, \cdots , i_n) \in B(\Lambda_n)_\mu$ has $n$  single entries then
${\rm dim}(V(\Lambda_n)_\mu) = 1$.

\begin{prop}\label{mult}
Suppose any vector $v = (i_1, \cdots , i_n) \in B(\Lambda_n)_\mu$ has $m < n$ number of single entries, then $n-m = 2k$ even and 
${\rm dim}(V(\Lambda_n)_\mu) = \frac{1}{k+1}{2k \choose k}$.
%Let $V_\lambda$ be a weight space. Then the dimension of $V_\lambda$ is equal to $C_k$, the k-th Catalan number where k is the number of pairs in the vectors in $V_\lambda$.
\end{prop}

\begin{proof}
It is know that the number of $k$ noncrossing arcs on $2k$ nodes is same as the $k$th Catalan number $C(k) =  \frac{1}{k+1}{2k \choose k}$ \cite{Sta}. To prove the proposition we show that there is a bijection between the set of elements $v  \in B(\Lambda_n)$ of weight $\mu$ and the set of  $k$ noncrossing arcs on $2k$ nodes. As stated above any two elements in $B(\Lambda_n)$ of weight $\mu$ have the same number of single entries $m$, say which are completely determined by the weight $\mu$ and they differ by $k = (n - m)/2$ number of pair entries.
In particular, $n - m = 2k$. Suppose the single entries in $v$ are $\{i_1, \cdots i_l, \overline{i_{l+1}} \cdots \overline{i_m}\}$. Denote $J = \{1, 2, \cdots , n\} \setminus \{i_1, \cdots i_l, i_{l+1},  \cdots , i_m\} = \{j_1, j_2, \cdots , j_{2k}\}$ ordered as in $S$. Join these $2k$ nodes by any $k$ noncrossing arcs. Suppose $j_t=p$ is a right end node of any such arc. Then the pair $(j_t, \overline{j_t}) = (p, \overline{p})$ is an admissible pair occurring in $v$ since it satisfies the condition $(2)$ in the definition of $B(\Lambda_n)$. Thus the $k$ end nodes in any set of $k$ noncrossing arcs  determine  the $k$ pairs which completely determine the elements $v$ in $B(\Lambda_n)_\mu$.

Next we show that any set of $k$ admissible pairs $\{p_1, \cdots , p_k, \overline{p_1}, \cdots , \overline{p_k}\}$ in $v \in B(\Lambda_n)$ with $m$ single entries $\{i_1, \cdots i_l, \overline{i_{l+1}} \cdots \overline{i_m}\}$ gives rise to a set of noncrossing arcs joining the $2k$ nodes. Note that the elements $\{q_1, \cdots , q_k\} = \{1, 2, \cdots , n\} \setminus \{i_1, \cdots i_l, i_{l+1},  \cdots , i_m\} \setminus \{p_1, \cdots , p_k\}$ do not appear in $v$. We order the $2k$ elements in the set $\{p_1, \cdots , p_k, q_1, \cdots , q_k\}$ by the order in $S$. We draw an arc from the left most available $p_i$ to the closest available $q_j$  occurring to its left and continue this process for all $p_i$. Suppose for some $p_i$, there is no available $q_j$ to its left. Let $(i_r, i_t) = (p_i, \overline{p_i})$ in the vector $v$. Then $t - r = n - p_i$ and $r-t+n+1 = p_i - n + n +1 = p_i +1 > p_i $ violating the condition $(2)$ in the definition of $B(\Lambda_n)$. Hence for each $p_i$ there is some available $q_j$ to its left which is connected by an arc. This process will create a set of $k$ noncrossing arcs on $2k$ nodes which completes the proof.

\end{proof}

\begin{example}
Consider the elements in the crystal $B(\Lambda_8)$ of the $C_8$- module $V(\Lambda_8)$ of weight $\mu=\varepsilon_3-\varepsilon_7= \Lambda_8-\alpha_1-2\alpha_2-2\alpha_3-3\alpha_4-4\alpha_5-5\alpha_6-7\alpha_7-4\alpha_8$. This implies that each element $v \in B(\Lambda_8)_\mu$ contains the singles $\{3, \overline{7}\}$. Then the following are the noncrossing arc diagrams on the nodes $1,2,4,5,6,8$ and the corresponding 
elements in $B(\Lambda_8)_\mu$:
%\begin{align*}
%&\MatchingB{6}{{1,2,4,5,6,8}}{1/2, 3/4, 5/6} \,\,\, \longrightarrow (2,3,5,8,\overline{8},\overline{7},\overline{5},\overline{2})\\
%&\MatchingB{6}{{1,2,4,5,6,8}}{1/4,2/3,5/6} \,\,\, \longrightarrow (3,4,5,8,\overline{8},\overline{7},\overline{5},\overline{4})\\
%&\MatchingB{6}{{1,2,4,5,6,8}}{1/2, 3/6, 4/5} \,\,\, \longrightarrow (2,3,6,8,\overline{8},\overline{7},\overline{6},\overline{6})\\
%&\MatchingB{6}{{1,2,4,5,6,8}}{1/6, 2/3, 4/5} \,\,\, \longrightarrow (3,4,6,8,\overline{8},\overline{7},\overline{6},\overline{4})\\
%&\MatchingB{6}{{1,2,4,5,6,8}}{1/6, 2/5, 3/4} \,\,\, \longrightarrow (3,5,6,8,\overline{8},\overline{7},\overline{6},\overline{5})  
%\end{align*}
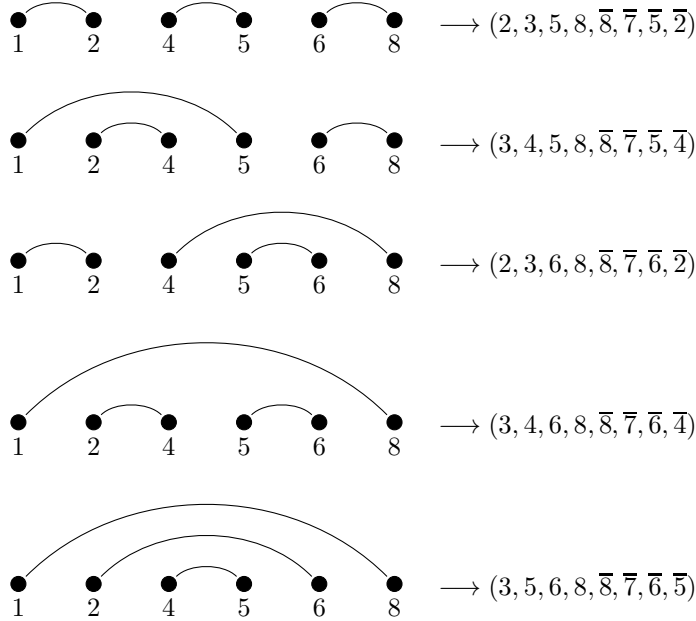
\begin{align*}
&\begin{tikzpicture}
       \draw[circle,fill] (1,0)circle[radius=1mm]node[below]{1};
       \draw[circle,fill] (2,0)circle[radius=1mm]node[below]{2};
       \draw[circle,fill] (3,0)circle[radius=1mm]node[below]{4};
       \draw[circle,fill] (4,0)circle[radius=1mm]node[below]{5};
       \draw[circle,fill] (5,0)circle[radius=1mm]node[below]{6};
       \draw[circle,fill] (6,0)circle[radius=1mm]node[below]{8};
       \draw(1,0) to[bend left=45] (2,0);
   \draw(3,0) to[bend left=45] (4,0);
   \draw(5,0) to[bend left=45] (6,0);
  \end{tikzpicture} \,\,\, \longrightarrow (2,3,5,8,\overline{8},\overline{7},\overline{5},\overline{2})\\
&\begin{tikzpicture}
       \draw[circle,fill] (1,0)circle[radius=1mm]node[below]{1};
       \draw[circle,fill] (2,0)circle[radius=1mm]node[below]{2};
       \draw[circle,fill] (3,0)circle[radius=1mm]node[below]{4};
       \draw[circle,fill] (4,0)circle[radius=1mm]node[below]{5};
       \draw[circle,fill] (5,0)circle[radius=1mm]node[below]{6};
       \draw[circle,fill] (6,0)circle[radius=1mm]node[below]{8};
       \draw(1,0) to[bend left=45] (4,0);
   \draw(2,0) to[bend left=45] (3,0);
   \draw(5,0) to[bend left=45] (6,0);
  \end{tikzpicture} \,\,\, \longrightarrow (3,4,5,8,\overline{8},\overline{7},\overline{5},\overline{4})\\
&\begin{tikzpicture}
       \draw[circle,fill] (1,0)circle[radius=1mm]node[below]{1};
       \draw[circle,fill] (2,0)circle[radius=1mm]node[below]{2};
       \draw[circle,fill] (3,0)circle[radius=1mm]node[below]{4};
       \draw[circle,fill] (4,0)circle[radius=1mm]node[below]{5};
       \draw[circle,fill] (5,0)circle[radius=1mm]node[below]{6};
       \draw[circle,fill] (6,0)circle[radius=1mm]node[below]{8};
       \draw(1,0) to[bend left=45] (2,0);
   \draw(3,0) to[bend left=45] (6,0);
   \draw(4,0) to[bend left=45] (5,0);
  \end{tikzpicture} \,\,\, \longrightarrow (2,3,6,8,\overline{8},\overline{7},\overline{6},\overline{2})\\
&\begin{tikzpicture}
       \draw[circle,fill] (1,0)circle[radius=1mm]node[below]{1};
       \draw[circle,fill] (2,0)circle[radius=1mm]node[below]{2};
       \draw[circle,fill] (3,0)circle[radius=1mm]node[below]{4};
       \draw[circle,fill] (4,0)circle[radius=1mm]node[below]{5};
       \draw[circle,fill] (5,0)circle[radius=1mm]node[below]{6};
       \draw[circle,fill] (6,0)circle[radius=1mm]node[below]{8};
       \draw(1,0) to[bend left=45] (6,0);
   \draw(2,0) to[bend left=45] (3,0);
   \draw(4,0) to[bend left=45] (5,0);
  \end{tikzpicture} \,\,\, \longrightarrow (3,4,6,8,\overline{8},\overline{7},\overline{6},\overline{4})\\
&\begin{tikzpicture}
       \draw[circle,fill] (1,0)circle[radius=1mm]node[below]{1};
       \draw[circle,fill] (2,0)circle[radius=1mm]node[below]{2};
       \draw[circle,fill] (3,0)circle[radius=1mm]node[below]{4};
       \draw[circle,fill] (4,0)circle[radius=1mm]node[below]{5};
       \draw[circle,fill] (5,0)circle[radius=1mm]node[below]{6};
       \draw[circle,fill] (6,0)circle[radius=1mm]node[below]{8};
       \draw(1,0) to[bend left=45] (6,0);
   \draw(2,0) to[bend left=45] (5,0);
   \draw(3,0) to[bend left=45] (4,0);
  \end{tikzpicture} \,\,\, \longrightarrow (3,5,6,8,\overline{8},\overline{7},\overline{6},\overline{5})  
\end{align*}

%These give all the elements in $B(\lambda_8)_\mu$ .
Hence ${\rm dim}(V(\L_8)_\mu) = 5 = C(3)$.
\end{example}

\subsection{Level zero fundamental module}

Let $U_q(\ge)$ denote the quantum affine algebra associated with $\ge$ and $U_q'(\ge)$ be the subalgebra of $U_q(\ge)$ without the degree derivation $d$. For the level $0$ extremal weight $\varpi_n = \L_n - \L_0$ we consider the finite dimensional irreducible $U_q'(\ge)$-module $W(\varpi_n)$ \cite{Kas3} called the $n^{th}$ fundamental representation. It is known that $W(\varpi_n)$ admits a crystal base. Let $B(\varpi_n)$ denote the crystal associated with $W(\varpi_n)$. As a set $B(\varpi_n)$ is same as the set $B(\Lambda_n)$ given in (\ref{crystal-B}). The actions of the Kashiwara operators $\{\Tilde{e}_k , \Tilde{f}_k \mid 0 \leq k \leq n\}$ on any element $(i_1, \cdots , i_n) \in B(\varpi_n)$ are given as follows:

\begin{align*}
    \Tilde{f}_0(i_1,...,i_n)=\begin{cases}
    (1,i_1,...,i_{n-1})& \text{if } i_n=\overline{1}\\
    0&\text{otherwise}
    \end{cases}
\end{align*}

\begin{align*}
    \Tilde{f}_k(i_1,...,i_n)=\begin{cases}
    (i_1,...,i_{j-1},k+1,i_{j+1},...,i_n)& \text{if } i_j=k, i_{j+1}\neq k+1, \  \text{and}\\
   &\text{either } \nexists l>j \text{ s.t } i_l=\overline{k}\\&\text{ or } \exists l>j \text{ s.t } i_l=\overline{k}, i_{l-1}=\overline{k+1}\\
    (i_1,...,i_{j-1},\overline{k},i_{j+1},...,i_n)& i_j=\overline{k+1}, i_{j+1}\neq \overline{k}, \ \text{and}\\&\text{either } \nexists l<j \text{ s.t } i_l=k\\&\text{ or } \exists l<j \text{ s.t } i_l=k, i_{l+1}=k+1\\
    0&\text{otherwise}
    \end{cases}
\end{align*}

\begin{align*}
    \Tilde{f}_n(i_1,...,i_n)=\begin{cases}
    (i_1,...,i_{j-1}, \overline{n},i_{j+1}.,i_{n-1})& \text{if } i_j=n, i_{j+1}\neq \overline{n}\\
    0&\text{otherwise}
    \end{cases}
\end{align*}

\begin{align*}
    \Tilde{e}_0(i_1,...,i_n)=\begin{cases}
    (i_2,...,i_n, \overline{1})& \text{if } i_1=1\\
    0&\text{otherwise}
    \end{cases}
\end{align*}

\begin{align*}
    \Tilde{e}_k(i_1,...,i_n)=\begin{cases}
    (i_1,...,i_{j-1},k,i_{j+1},...,i_n)& \text{if } i_j=k+1, i_{j-1}\neq k, \text{and}\\
    &\text{either } \nexists l>j \text{ s.t }. i_l=\overline{k}\\&\text{ or } \exists l>j \text{ s.t } i_l=\overline{k}, i_{l-1}=\overline{k+1}\\
    (i_1,...,i_{j-1},\overline{k+1},i_{j+1},...,i_n)& i_j=\overline{k}, i_{j-1}\neq \overline{k+1}, \text{and}\\&\text{either } \nexists l<j \text{ s.t } i_l=k+1\\&\text{ or } \exists l<j \text{ s.t } i_l=k+1, i_{l-1}=k\\
    0&\text{otherwise}
    \end{cases}
\end{align*}

\begin{align*}
    \Tilde{e}_n(i_1,...,i_n)=\begin{cases}
    (i_1,...,i_{j-1}, n,i_{j+1}.,i_{n-1})& \text{if } i_j=\overline{n}, i_{j-1}\neq n\\
    0&\text{otherwise}
    \end{cases}
\end{align*}

Note that the element $(1, 2, \cdots , n)$ (resp. $(\overline{n}, \overline{n-1}, \cdots , \overline{1})$) is the unique element with weight $\varpi_n = \L_n - \L_0$ (resp. $\check{\varpi}_n=\L_0-\L_n$) in $B(\varpi_n)$. Also by Proposition \ref{mult} the crystal $B(\varpi_n)$ is not multiplicity free for $n \geq 4$. For $n = 2$ and $3$, the crystal base coincides with the global base and the action of the Chevalley generators $\{e_k, f_k \mid 0 \leq k \leq n\}$ are given by the same formula as the Kashiwara operators $\{\Tilde{e}_k , \Tilde{f}_k \mid 0 \leq k \leq n\}$. For $n=4$, 
% It is known that the translation $t(\varpi_n )$ (resp. $t(\check{\varpi}_n)$) associated with the extremal level $0$ weight $\varpi_n$ (resp. $\check{\varpi}_n=\L_0-\L_n$) is in the extended Weyl group $\widetilde{W}$. We determine these translations explicitly.

\subsection{Translation of some level $0$ extremal weights}

% Let $\mathfrak{g}_0$ be the subalgebra of $\mathfrak{g}$ associated with the index set $I/\{0\}$ and $\mathfrak{g}_n$ be the subalgebra of $\mathfrak{g}$ associated with the index set $I/\{n\}$. Both are isomorphic to $C_n$.

%Let $W(\varpi_n)$ be the level $0$ fundamental $U'_q(\mathfrak{g})$-module associated with the level $0$ weight $\varpi_n = \L_n - \L_0$. $W(\varpi_n)$ is a finite-dimensional irreducible integrable $U'_q(\mathfrak{g})$-module. Similarly, $W(\check{\varpi}_n)$ is a level 0 fundamental $U'_q(\mathfrak{g})$-module where $\check{\varpi}_n=\L_0-\L_n$.

%The Weyl group for $C_n^{(1)}$ is generated by simple reflections $s_i(\lambda)=\lambda-\lambda(\check{\alpha}_i)\alpha_i$. 
It is known that the translation $t(\varpi_n )$ (resp. $t(\check{\varpi}_n)$) associated with the extremal level $0$ weight $\varpi_n = \L_n - \L_0$ (resp. $\check{\varpi}_n=\L_0-\L_n$) is in the extended Weyl group $\widetilde{W}$. Now we determine these translations explicitly.
For any $\lambda \in P_{cl}$ we write $\lambda = (\lambda_0, \lambda_1, \cdots , \lambda_n)$ if $\lambda =\lambda_0\Lambda_0 + \lambda_1\Lambda_1 + \cdots + \lambda_n\Lambda_n$. Then the explicit actions of  the simple reflections on $\lambda$ are given as follows:
  $$s_0(\lambda_0,\lambda_1, ,..., \lambda_n)=(-\lambda_0, 2\lambda_0+\lambda_1,..., \lambda_n)$$
$$s_i(\lambda_0,...,\lambda_i,...,  \lambda_n)=(\lambda_0,...,\lambda_{i-1}+\lambda_i, -\lambda_i,\lambda_{i+1}+\lambda_i,..., \lambda_n)\text{ for }i\in I\setminus\{0,n\}$$
$$s_n(\lambda_0,\lambda_1,..., \lambda_n)=(\lambda_0, \lambda_1,....,\lambda_{n-1}+2\lambda_n, -\lambda_n)$$

Let $\sigma$ denote the  Dynkin diagram automorphism which maps $i \leftrightarrow n - i$, for $i = 0, 1, \cdots \lfloor\frac{n}{2}\rfloor$. 
%Then the translations $t(\varpi_n)$ and $t(\check\varpi_n )$ are as follows.

\begin{prop}
The translations $t(\varpi_n) = \sigma s_ns_{n-1}s_ns_{n-2}s_{n-1}s_n....s_{1}...s_{n}$ and $t(\Check{\varpi}_n) = \sigma s_0s_{1}s_0s_{2}s_{1}s_0....s_{n-1}...s_{0}$.
\end{prop}
\begin{proof}
From the formulas for simple reflections above, we see that $s_i$ only affects the $i-1, i, i+1$ coordinates in $\lambda$ for $i \in I\setminus \{0, n\}$, the $0, 1$ coordinates for $i = 0$ and the $n-1, n$ coordinates for $i = n$. Hence we have
\begin{align*}
    &\sigma s_ns_{n-1}s_ns_{n-2}s_{n-1}s_n....s_{1}...s_{n}(\lambda_0,\lambda_1,...,\lambda_n)\\
    &=\sigma s_ns_{n-1}s_ns_{n-2}s_{n-1}s_n....s_{2}...s_{n}(\lambda_0+...+\lambda_{n-1}+2\lambda_n, -\lambda_1-...-2\lambda_n,\lambda_1,...,\\&\lambda_{n-2},\lambda_n+\lambda_{n-1})\\
    &=\sigma s_ns_{n-1}s_ns_{n-2}s_{n-1}s_n....s_{3}...s_{n}(\lambda_0+...+\lambda_{n-1}+2\lambda_n,\lambda_{n-1},-\lambda_2-...-2\lambda_{n-1}\\&-2\lambda_n,\lambda_1,...,\lambda_{n-2},\lambda_{n-2}+\lambda_{n-1}+\lambda_n)\\
    &\vdots\\
    &=\sigma(\lambda_0+...+\lambda_{n-1}+2\lambda_n,\lambda_{n-1},...,\lambda_1,-\lambda_1-...-\lambda_n)\\
    &=(-\lambda_1-...-\lambda_n,\lambda_1,...,\lambda_{n-1},\lambda_0+...+\lambda_{n-1}+2\lambda_n)\\
    &=(\lambda_0,\lambda_1,...,\lambda_n)+(\lambda_0+...+\lambda_{n-1}+\lambda_n)(\L_n-\L_0)
\end{align*}
The calculations for $t(\check\varpi_n)$ is similar.
\end{proof}

\section{Affine Geometric Crystal associated to \bf{$W(\varpi_n)$}}
In this section we will construct the affine geometric crystal $\mathcal{V}$  explicitly for $n=2,3,4$.
Consider the level $0$ fundamental weight $\varpi_n = \Lambda_n - \Lambda_0$. Note that $\sigma(\varpi_n) = \Check{\varpi}_n = \Lambda_0 - \Lambda_n$. The finite dimensional irreducible $U_q(\ge)$-module $W(\varpi_n)$, called the $n^{th}$ fundamental module admits a global basis with a simple crystal \cite{Kas3}. So we can consider the specialization $q = 1$ and obtain the finite dimensional irreducible $\ge$-module $W(\varpi_n)$, called a fundamental representation of $\ge$ and use the same notations. The explicit action of $\ge$ on $W(\varpi_n)$ is given below.

\subsection{Fundamental representation $W(\varpi_n)$ for $\ge = C_n^{(1)}$}
The $\ge$-module $W(\varpi_n)$ is same as $V(\Lambda_n)$ as a vector space and hence an $({2n \choose n} - {2n \choose {n-2}})$ - dimensional module with the basis:
\begin{equation*}\label{basis}
\small{\left\{(i_1, \cdots , i_n) \mid
 \begin{tabular}{c} 
$(1) 1 \prec i_1 \prec \cdots \prec i_n \prec \overline{1},$ $i_j \in S$ \\
$(2) {\rm if} \ i_k = p \  {\rm and}  \  i_l = \overline{p},  \  {\rm then} \ k+(n-l+1) \leq p $
\end{tabular}\right\}
}
\end{equation*}
where $S = \{i, \overline{i} \mid 1\leq i \leq n\}$ with the order $1\preceq 2\prec...\prec n\prec\overline{n}\prec...\preceq \overline{1}$. By Proposition \ref{mult} for $2 \leq n \leq 4$, each weight space of $W(\varpi_n)$ is one-dimensional except when $n = 4$, the $0$-weight space is two dimensional. When the weight space is one dimensional we can identify the global basis vector with the corresponding crystal basis vector. Therefore, the actions of the Chevalley generators $\{e_i , f_i  \mid 0 \leq i \leq n\}$ on the basis vectors coincide with the actions of the Kashiwara operators $\{\Tilde{e}_i , \Tilde{f}_i  \mid 0 \leq i \leq n\}$ except in the following cases when $n = 4$:
\begin{equation}
e_2(3,4,\bar{4}, \bar{2}) = (2, 4, \bar{4}, \bar{2}) + (3, 4, \bar{4}, \bar{3}) = f_2(2,4,\bar{4}, \bar{3}).
\end{equation}
%\begin{equation}
%f_2(2,4,\bar{4}, \bar{3}) = (2, 4, \bar{4}, \bar{2}) + (3, 4, \bar{4}, \bar{3}).
%\end{equation}
\subsection{$n=2$ Case}
From Proposition 2.4 and 2.5, we know that the translations associated with the geometric crystals $\mathcal{V}_1(x)$ and $\mathcal{V}_2(y)$ will be $\sigma s_2s_1s_2$ and $\sigma s_0s_1s_0$ respectively. Then the varieties associated with these Weyl group elements are defined as follows. 

For $\mathcal{V}_1(x)$ we have:
\begin{align*}
&Y_2(x_{22})Y_1(x_{11})Y_2(x_{21})(1,2)\\
&=Y_2(x_{22})Y_1(x_{11})[x_{21}(1,2)+(1,\overline{2})]\\
&=Y_2(x_{22})[x_{21}(1,2)+x_{11}^2(1,\overline{2})+x_{11}(2,\overline{2})+(2, \overline{1})]\\
&=[x_{21}x_{22}(1,2)+(x_{21}+\frac{x_{11}^2}{x_{22}})(1,\overline{2})+x_{11}(2,\overline{2})+x_{22}(2,\overline{1})+(\overline{2}, \overline{1})]
\end{align*}
and for $V_2(y)$ we have:
\begin{align*}
&Y_0(y_{02})Y_1(y_{11})Y_0(y_{01})(\overline{2},\overline{1})\\
&=Y_0(y_{02})Y_(y_{11})[y_{01}(\overline{2},\overline{1})+(1,\overline{2})]\\
&=Y_0(y_{02})[y_{01}(\overline{2},\overline{1})+y_{11}^2(1,\overline{2})+y_{11}(2,\overline{2})+(2, \overline{1})]\\
&= [y_{01}y_{02}(\overline{2},\overline{1})+(y_{01}+\frac{y_{11}^2}{y_{02}})(1,\overline{2})+y_{11}(2,\overline{2})+y_{02}(2,\overline{1})+(1,2)]
\end{align*}\\
Setting $V_1(x)\cdot a(x)=V_2(y)$, we obtain the following relations:
$$1=a(x)\cdot x_{21}x_{22},\,\,\,y_{01}y_{02}=a(x),\,\,\,(x_{21}+\frac{x_{11}^2}{x_{22}})a(x)=y_{01}+\frac{y_{11}^2}{y_{02}}$$
$$x_{11}a(x)=y_{11},\,\,\,x_{22}a(x)=y_{02} $$

with the following solutions: $a(x)=\frac{1}{x_{21}x_{22}}$, $y_{11}=\frac{x_{11}}{x_{21}x_{22}}$, $y_{02}=\frac{1}{x_{21}}$ and $y_{01}=\frac{1}{x_{22}}$. 
Let $\mathcal{V}_1$ be the geometric crystal corresponding to $\mathfrak{g}_0$ and $\mathcal{V}_2$ be the geometric crystal corresponding to $\mathfrak{g}_n$. This map is clearly birational and bipositive, and  defines the isomorphism $\overline{\sigma}:\mathcal{V}_1 \longrightarrow \mathcal{V}_2$ and $\overline{\sigma^{-1}}:\mathcal{V}_2 \longrightarrow \mathcal{V}_1$.\\
\bigskip
Now we give the actions $e_i^c$, $\varepsilon_i$, and $\gamma_i$ from the general formula in section 5.2 for $i=1,2$.

$$e_1^c(V_1(x_{22},x_{11},x_{21}))=(x_{21},cx_{11}, x_{22})$$

Let $$c_2:=\frac{cx_{21}x_{22}+x_{11}^2}{x_{21}x_{22}+x_{11}^2}$$ 
Then $$e_2^c(V_1(x_{22},x_{11},x_{21})=\left(c_2x_{22},x_{11}, \frac{c}{c_2}x_{21}\right)$$

$$\varepsilon_1(x_{22},x_{11},x_{21})=\frac{x_{22}}{x_{11}}$$
$$\varepsilon_2(x_{22},x_{11},x_{21})=\frac{x_{21}x_{22}+x_{11}^2}{x_{22}^2x_{21}}$$
$$\gamma_1(x_{22},x_{11},x_{21})=\frac{x_{11}^2}{x_{21}x_{22}}$$
$$\gamma_2(x_{22},x_{11},x_{21})=\frac{x_{21}^2x_{22}^2}{x_{11}^2}$$
To obtain the formulas for $e_0^c$, $\varepsilon_0$ and $\gamma_0$ we first obtain the formulas for $\overline{e_0^c}$, $\overline{\varepsilon_0}$ and $\overline{\gamma_0}$ in $\mathcal{V}_2(y)$. These are as follows:
Let
$$\overline{c}_0=\frac{cy_{01}y_{02}+y_{11}^2}{y_{01}y_{02}+y_{11}^2}$$
Then 
$$\overline{e_0^c}(y_{02},y_{11},y_{01})=(\overline{c}_0y_{02}, y_{11}, \frac{c}{\overline{c}_0}y_{01})$$
$$\overline{\varepsilon_0}(y_{02},y_{11},y_{01})=\frac{1}{y_{02}}+\frac{y_{11}^2}{y_{02}^2y_{01}}=\frac{y_{01}y_{02}+y_{11}^2}{y_{02}^2y_{01}}$$
$$\overline{\gamma_0}(y_{02},y_{11},y_{01})=\frac{y_{01}^2y_{02}^2}{y_{11}^2}$$
Now $e_0^c$, $\varepsilon_0$ and $\gamma_0$ are defined as follows:
\begin{align*}
    e_0^c(V_1(x))=\overline{\sigma^{-1}}\circ \overline{e_0^c}\circ \overline{\sigma}(V_1(x))\\
    \gamma_0(V_1(x))=\overline{\gamma_0}(\overline{\sigma}(V_1(x)))\\
    \varepsilon_0(V_1(x))=\overline{\varepsilon_0}(\overline{\sigma}(V_1(x)))
\end{align*}
From these formulas we obtain the following formulas for $e_0^c$, $\varepsilon_0$ and $\gamma_0$.
\begin{align*}
     &e_0^c(x_{22},x_{11},x_{21})=\overline{\sigma^{-1}}\circ \overline{e_0^c}\circ \overline{\sigma}(x_{22},x_{11},x_{21})=
     (\frac{c_2}{c}x_{22}, \frac{x_{11}}{c}, \frac{x_{21}}{c_2})\\
&\varepsilon_0(V_1(x))=\overline{\varepsilon_0}(\overline{\sigma}(V_1(x)))=\overline{\varepsilon_0}\left(\frac{1}{x_{21}},\frac{x_{11}}{x_{21}x_{22}}, \frac{1}{x_{22}}\right)=x_{21}+\frac{x_{11}^2}{x_{22}}\\
&\gamma_0(V_1(x))=\overline{\gamma_0}(\overline{\sigma}(V_1(x)))=\frac{1}{x_{11}^2}
\end{align*}
\begin{theorem}
 $\mathcal{V}_1=(V_1(x),e_i^c,\varepsilon_i,\varphi_i)$ for $i=0,1,2$ is a positive geometric crystal.
\end{theorem}
\begin{proof}
Clearly this is positive, since all coefficients are positive. We only need to check the relations involving $e_0^c,\varphi_0$, and $\varepsilon_0$, since by \cite{BK, nakashima2005geometric} we already know $\mathcal{V}_1=(V_1(x),e_i^c,\varepsilon_i, \varphi_i)$ for $i=1,2$ is a geometric crystal. First we check the relations $\gamma_j(e_i^c(x))=c^{a_{ij}}\gamma_j(x)$:
\begin{align*}
    \gamma_0(e_0^c(x))=\gamma_0(\frac{c_2}{c}x_{22},\frac{x_{11}}{c},\frac{1}{c_2} x_{21})=\frac{c^2}{x_{11}^2}=c^{2}\gamma_0(x)\\
    \gamma_0(e_1^c(x))=\gamma_0(x_{22},cx_{11},x_{21})=\frac{1}{c^2x_{11}^2}=c^{-2}\gamma_0(x)\\
    \gamma_0(e_2^c(x))=\gamma_0(c_2x_{22},x_{11},\frac{c}{c_2}x_{21})=\frac{1}{x_{11}^2}=c^0\gamma_0(x)\\
    \gamma_1(e_0^c(x))=\gamma_1(\frac{c_2}{c}x_{22},\frac{x_{11}}{c},\frac{1}{c_2}x_{21})=\frac{cx_{11}^2}{c^2x_{21}x_{22}}=c^{-1}\gamma_1(x)\\
    \gamma_2(e_0^c(x))=\gamma_2(\frac{c_2}{c}x_{22},\frac{x_{11}}{c},\frac{1}{c_2} x_{21})=\frac{c^2x_{21}^2x_{22}^2}{c^2x_{11}^2}=c^0\gamma_2(x)
\end{align*}
Then we check the relations $\varepsilon_i(e_i^c(x))=c^{-1}\varepsilon_i(x)$,
    $\varepsilon_0(e_2^c(x))=\varepsilon_0(x)$ and
    $\varepsilon_2(e_0^c(x))=\varepsilon_2(x)$.
\begin{align*}
    &\varepsilon_0(e_0^c(x))=\varepsilon_0(\frac{c_2}{c}x_{22}, \frac{x_{11}}{c}, \frac{x_{21}}{c_2})=\frac{x_{21}}{c_2}+\frac{x_{11}^2}{c^2\frac{c_2}{c}x_{22}}=\frac{cx_{22}x_{21}+x_{11}^2}{cc_2x_{22}}\\
    &=\frac{(x_{22}x_{21}+x_{11}^2)}{cx_{22}\frac{cx_{21}x_{22}+x_{11}^2}{x_{21}x_{22}+x_{11}^2}}=\frac{(x_{21}x_{22}+x_{11}^2)}{cx_{22}}=\frac{x_{21}}{c}+\frac{x_{11}^2}{cx_{22}}=c^{-1}\varepsilon_0(x)\\
    &\varepsilon_0(e_2^c(x))=\varepsilon_2(c_2x_{22}, x_{11}, \frac{c}{c_2}x_{21})=\frac{cx_{21}}{c_2}+\frac{x_{11}^2}{c_2x_{22}}=\frac{cx_{21}x_{22}+x_{11}^2}{x_{22}\frac{cx_{21}x_{22}+x_{11}^2}{x_{21}x_{22}+x_{11}^2}}=x_{21}+\frac{x_{11}^2}{x_{22}}\\&=\varepsilon_0(x)\\
    &\varepsilon_2(e_0^c(x))=\varepsilon_2(\frac{c_2}{c}x_{22},\frac{x_{11}}{c}, \frac{x_{21}}{c_2})=\frac{cx_{21}x_{22}+x_{11}^2}{x_{22}^2x_{21}\frac{cx_{21}x_{22}+x_{11}^2}{x_{21}x_{22}+x_{11}^2}}=\frac{x_{21}x_{22}+x_{11}^2}{x_{22}^2x_{21}}=\varepsilon_2(x)
\end{align*}
Now we check the relation $e_0^ce_2^d=e_2^de_0^c$
\begin{align*}
&e_0^ce_2^d(x_{22},x_{11},x_{21})=e_0^c(d_2x_{22},x_{11},\frac{d}{d_2}x_{21})\\&=e_0^c\left(\frac{x_{22}(dx_{21}x_{22}+x_{11}^2)}{x_{21}x_{22}+x_{11}^2},x_{11},\frac{dx_{21}(x_{21}x_{22}+x_{11}^2)}{dx_{21}x_{22}+x_{11}^2}\right)\\
&=\left(\frac{c_2^*}{c}\frac{x_{22}(dx_{21}x_{22}+x_{11}^2)}{x_{21}x_{22}+x_{11}^2},\frac{x_{11}}{c},\frac{1}{c_2^*}\frac{dx_{21}(x_{21}x_{22}+x_{11}^2)}{dx_{21}x_{22}+x_{11}^2}\right)
\end{align*}
We will solve for $c_2^*$ to substitute it in:
\begin{align*}
    c_2^*=\frac{c(\frac{(dx_{21}x_{22}+x_{11}^2)x_{22}}{dx_{21}x_{22}+x_{11}^2}\cdot\frac{(dx_{21}x_{22}+x_{11}^2)x_{22}}{x_{21}x_{22}+x_{11}^2})+x_{11}^2}{\frac{(dx_{21}x_{22}+x_{11}^2)dx_{21}}{x_{21}x_{22}+x_{11}^2}\cdot\frac{dx_{21}(x_{21}x_{22}+x_{11}^2)}{dx_{21}x_{22}+x_{11}^2}+x_{11}^2}=\frac{cdx_{21}x_{22}+x_{11}^2}{dx_{21}x_{22}+x_{11}^2}
\end{align*}
Substituting this in, we obtain the following
\begin{align*}
    &e_0^ce_2^d(x_{22},x_{11},x_{21})=\left(\frac{x_{22}(cdx_{21}x_{22}+x_{11}^2)}{c(x_{21}x_{22}+x_{11}^2)},\frac{x_{11}}{c}, \frac{dx_{21}(x_{21}x_{22}+x_{11}^2)}{cdx_{21}x_{22}+x_{11}^2}\right)
\end{align*}
Now we calculate the other side of the equation:
\begin{align*}
 e_2^de_0^c(x_{22},x_{11},x_{21})=e_2^d\left( \frac{x_{22}(cx_{21}x_{22}+x_{11}^2)}{c(x_{21}x_{22}+x_{11}^2)},\frac{x_{11}}{c}, \frac{x_{21}(x_{21}x_{22}+x_{11}^2)}{cx_{21}x_{22}+x_{11}^2}\right)\\
 =\left(d_2^*\frac{x_{22}(cx_{21}x_{22}+x_{11}^2)}{c(x_{21}x_{22}+x_{11}^2)},\frac{x_{11}}{c},\frac{d}{d_2^*}\frac{x_{21}(x_{21}x_{22}+x_{11}^2)}{cx_{21}x_{22}+x_{11}^2}\right)
\end{align*}
We solve for $c_2^*$ to substitute it in:
\begin{align*}
    d_2^*=\frac{d(\frac{(cx_{21}x_{22}+x_{11}^2)x_{22}}{c(x_{21}x_{22}+x_{11}^2)}\cdot\frac{(x_{21}x_{22}+x_{11}^2)x_{21}}{cx_{21}x_{22}+x_{11}^2})+\frac{x_{11}^2}{c^2}}{\frac{(cx_{21}x_{22}+x_{11}^2)dx_{22}}{c(x_{21}x_{22}+x_{11}^2)}\cdot\frac{x_{21}(x_{21}x_{22}+x_{11}^2)}{cx_{21}x_{22}+x_{11}^2}+\frac{x_{11}^2}{c^2}}=\frac{cdx_{21}x_{22}+x_{11}^2}{cx_{21}x_{22}+x_{11}^2}
\end{align*}
Then we have:
\begin{align*}
&e_2^de_0^c(x_{22},x_{11},x_{21})=\left(\frac{x_{22}(cdx_{21}x_{22}+x_{11}^2)}{c(x_{21}x_{22}+x_{11}^2)},\frac{x_{11}}{c}, \frac{dx_{21}(x_{21}x_{22}+x_{11}^2)}{cdx_{21}x_{22}+x_{11}^2}\right)
\end{align*}
Both sides are equal, so the relation holds. Finally, we check the last relation:
\begin{align*}
    &e_1^ce_0^{c^2d}e_1^{cd}e_0^d(x_{22},x_{11},x_{21})\\&=e_1^ce_0^{c^2d}e_1^{cd}\left(\frac{x_{22}(dx_{21}x_{22}+x_{11}^2)}{d(x_{21}x_{22}+x_{11}^2)},\frac{x_{11}}{d}, \frac{x_{21}(x_{21}x_{22}+x_{11}^2)}{dx_{21}x_{22}+x_{11}^2}\right)\\
    &=e_1^ce_0^{c^2d}\left(\frac{x_{22}(dx_{21}x_{22}+x_{11}^2)}{d(x_{21}x_{22}+x_{11}^2)},cx_{11}, \frac{x_{21}(x_{21}x_{22}+x_{11}^2)}{dx_{21}x_{22}+x_{11}^2}\right)\\
    &=e_1^c\left(\frac{(c^2d)^*}{c^2d}\cdot\frac{x_{22}(dx_{21}x_{22}+x_{11}^2)}{d(x_{21}x_{22}+x_{11}^2)}, \frac{x_{11}}{cd},\frac{1}{(c^2d)^*}\cdot\frac{x_{21}(x_{21}x_{22}+x_{11}^2)}{dx_{21}x_{22}+x_{11}^2}\right)
\end{align*}
Then we solve for $(c^2d)^*$
\begin{align*}
    &(c^2d)^*=\frac{c^2d\left(\frac{x_{22}(dx_{21}x_{22}+x_{11}^2)}{d(x_{21}x_{22}+x_{11}^2)}\cdot\frac{x_{21}(x_{21}x_{22}+x_{11}^2)}{dx_{21}x_{22}+x_{11}^2}\right)+c^2x_{11}^2 }{\frac{x_{22}(dx_{21}x_{22}+x_{11}^2)}{d(x_{21}x_{22}+x_{11}^2)}\cdot\frac{x_{21}(x_{21}x_{22}+x_{11}^2)}{dx_{21}x_{22}+x_{11}^2}+c^2x_{11}^2}\\&=\frac{\frac{c^2dx_{21}x_{22}}{d}+c^2x_{11}^2}{\frac{x_{21}x_{22}}{d}+c^2x_{11}^2}=\frac{c^2d(x_{21}x_{22}+x_{11}^2)}{x_{21}x_{22}+c^2dx_{11}^2}
\end{align*}
Then we have
\begin{align*}
    &e_1^c\left(\frac{x_{22}(dx_{21}x_{22}+x_{11}^2)}{d(x_{21}x_{22}+c^2dx_{11}^2)},\frac{x_{11}}{cd},\frac{x_{21}(x_{21}x_{22}+c^2dx_{11}^2)}{c^2d(dx_{21}x_{22}+x_{11}^2)} \right)\\
    &=\left(\frac{x_{22}(dx_{21}x_{22}+x_{11}^2)}{d(x_{21}x_{22}+c^2dx_{11}^2)},\frac{x_{11}}{d},\frac{x_{21}(x_{21}x_{22}+c^2dx_{11}^2)}{c^2d(dx_{21}x_{22}+x_{11}^2)} \right)
\end{align*}
Now we evaluate the other side of the relation:
\begin{align*}
    &e_0^de_1^{cd}e_0^{c^2d}e_1^c(x_{22},x_{11},x_{21})=e_0^de_1^{cd}e_0^{c^2d}(x_{22},cx_{11},x_{21})\\
    &=e_0^de_1^{cd}\left(\frac{x_{22}(dx_{21}x_{22}+x_{11}^2)}{d(x_{21}x_{22}+c^2x_{11}^2)},\frac{x_{11}}{cd}, \frac{x_{21}(x_{21}x_{22}+c^2x_{11}^2)}{c^2(dx_{21}x_{22}+x_{11}^2)}\right)\\
    &=e_0^d\left(\frac{x_{22}(dx_{21}x_{22}+x_{11}^2)}{d(x_{21}x_{22}+c^2x_{11}^2)},x_{11}, \frac{x_{21}(x_{21}x_{22}+c^2x_{11}^2)}{c^2(dx_{21}x_{22}+x_{11}^2)}\right)\\
    &=\left(\frac{d_2^*}{d}\cdot\frac{x_{22}(dx_{21}x_{22}+x_{11}^2)}{d(x_{21}x_{22}+c^2x_{11}^2)},\frac{x_{11}}{d}, \frac{1}{d_2^*}\cdot\frac{x_{21}(x_{21}x_{22}+c^2x_{11}^2)}{c^2(dx_{21}x_{22}+x_{11}^2)}\right)
\end{align*}
We solve for $d_2^*$:
\begin{align*}
    &d_2^*=\frac{\frac{dx_{21}x_{22}}{c^2d}+x_{11}^2}{\frac{x_{21}x_{22}}{c^2d}+x_{11}^2}=\frac{dx_{21}x_{22}+c^2dx_{11}^2}{x_{21}x_{22}+c^2dx_{11}^2}
\end{align*}
Substituting this expression in, we get:
\begin{align*}
  \left(\frac{x_{22}(dx_{21}x_{22}+x_{11}^2)}{d(x_{21}x_{22}+c^2dx_{11}^2)},\frac{x_{11}}{d},\frac{x_{21}(x_{21}x_{22}+c^2dx_{11}^2)}{c^2d(dx_{21}x_{22}+x_{11}^2)} \right)
\end{align*}

Since all of the relations are satisfied, $\mathcal{V}_1(x)$ along with the actions $e_i^c$, and functions $\gamma_i$ and $\varepsilon_i$ for $i=0,1,2$ is a $C_2^{(1)}$ geometric crystal.
\end{proof}
\subsection{n=3 Case}
The Weyl group elements are $t(\varpi_3)=s_3s_2s_3s_1s_2s_3$ and $t(\check{\varpi}_3)=s_0s_1s_0s_2s_1s_0$. We first compute $V_1(x)$ and $V_2(y)$.
\begin{align*}
& Y_3(x_{33})Y_2(x_{22})Y_3(x_{32})Y_1(x_{11})Y_2(x_{21})Y_3(x_{31})(1,2,3)\\
&=[x_{31}x_{32}x_{33}(1,2,3)+(x_{31}x_{32}+\frac{x_{31}x_{22}^2}{x_{33}}+\frac{x_{21}^2x_{22}^2}{x_{32}x_{33}})(1,2,\overline{3})+(x_{21}x_{11}+\frac{x_{21}^2x_{22}}{x_{32}}\\&+x_{31}x_{22})(1,3,\overline{3})+(\frac{x_{11}^2x_{32}}{x_{22}^2}+x_{31}+\frac{x_{21}^2}{x_{32}}+\frac{2x_{21}x_{11}}{x_{22}})x_{33}(1,3,\overline{2})+(\frac{x_{11}^2x_{32}}{x_{22}^2}+x_{31}\\&+\frac{x_{21}^2}{x_{32}}+\frac{2x_{21}x_{11}}{x_{22}}+\frac{x_{11}^2}{x_{33}})(1,\overline{3},\overline{2})+x_{21}x_{22}(2,3,\overline{3})+ (x_{21}+\frac{x_{11}x_{32}}{x_{22}})x_{33}(2,3,\overline{2})\\&+(x_{21}+\frac{x_{11}x_{32}}{x_{22}}+\frac{x_{11}x_{22}}{x_{33}})(2,\overline{3},\overline{2})+x_{32}x_{33}(2,3,\overline{1})+(x_{32}+\frac{x_{22}^2}{x_{33}})(2,\overline{3},\overline{1})\\&+x_{11}(3,\overline{3},\overline{2})+x_{22}(3,\overline{3},\overline{1})+x_{33}(3,\overline{2},\overline{1})+(\overline{3},\overline{2},\overline{1})]
\end{align*}
\begin{align*}
& V_2(y)=Y_0(y_{03})Y_1(y_{12})Y_0(y_{02})Y_2(y_{21})Y_1(y_{11})Y_0(y_{01})(\overline{3},\overline{2},\overline{1})\\
&=[y_{01}y_{02}y_{03}(\overline{3},\overline{2},\overline{1})+(y_{01}y_{02}+\frac{y_{01}y_{12}^2}{y_{03}}+\frac{y_{11}^2y_{12}^2}{y_{02}y_{03}})(1,\overline{3},\overline{2})+(y_{11}y_{21}+y_{01}y_{12}\\&+\frac{y_{11}^2y_{12}}{y_{02}})(2,\overline{3},\overline{2})+(\frac{2y_{11}y_{21}}{y_{12}}+y_{01}+\frac{y_{11}^2}{y_{02}}+\frac{y_{21}^2y_{02}}{y_{12}^2})y_{03}(2,\overline{3},\overline{1})+(\frac{2y_{11}y_{21}}{y_{12}}+y_{01}\\&+\frac{y_{11}^2}{y_{02}}+\frac{y_{21}^2y_{02}}{y_{12}^2}+\frac{y_{21}^2}{y_{03}})(1,2,\overline{3})+y_{11}y_{12}(3,\overline{3},\overline{2})+(y_{11}+\frac{y_{21}y_{02}}{y_{12}})y_{03}(3,\overline{3},\overline{1})\\&+(\frac{y_{21}y_{12}}{y_{03}}+y_{11}+\frac{y_{21}y_{02}}{y_{12}})(1,3,\overline{3})+y_{21}(2,3,\overline{3})+y_{02}y_{03}(3,\overline{2},\overline{1})\\&+(y_{02}+\frac{y_{12}^2}{y_{03}})(1,3,\overline{2})+y_{12}(2,3,\overline{2})+y_{03}(2,3,\overline{1})+(1,2,3)]
\end{align*}

Now we solve for the coefficients in terms of one another using the equation: $V_1(x)a(x)=V_2(y)$
We get the following solution:
$a(x)=\frac{1}{x_{31}x_{32}x_{33}}$, $y_{21}=\frac{x_{21}x_{22}}{x_{31}x_{32}x_{33}}$, $y_{03}=\frac{1}{x_{31}}$ and $y_{02}=\frac{1}{x_{32}}$, $y_{01}=\frac{1}{x_{33}}$,$y_{12}=\frac{x_{21}x_{22}+x_{11}x_{32}}{x_{31}x_{32}x_{22}}$ and $y_{11}=\frac{x_{11}x_{22}}{x_{33}(x_{21}x_{22}+x_{11}x_{32})}$. We can also see that $x_{31}=\frac{1}{y_{03}}$,  $x_{32}=\frac{1}{y_{02}}$, $x_{33}=\frac{1}{y_{01}}$,  $x_{21}=\frac{y_{21}y_{12}}{y_{03}(y_{21}y_{02}+y_{11}y_{12})}$, $x_{22}=\frac{y_{21}y_{02}+y_{11}y_{12}}{y_{01}y_{02}y_{12}}$ and $x_{11}=\frac{y_{11}y_{12}}{y_{01}y_{02}y_{03}}$.\\
From the general formula, we need to compute the maps for this geometric crystal: $e_i^c$, $\varepsilon_i$ and $\gamma_i$ for $i=0,1,2,3$.

First we compute $e_1^c$, $e_2^c$ and $e_3^c$.
\begin{align*}
&e_1^c(x_{33},x_{22},x_{32},x_{11},x_{21},x_{31})=(x_{33},x_{22},x_{32},cx_{11},x_{21},x_{31})
\\
&e_2^c(x_{33},x_{22},x_{32},x_{11},x_{21},x_{31})
=(x_{33},c_2x_{22},x_{32},x_{11},\frac{c}{c_2}x_{21},x_{31})
\end{align*}

where $$c_2=\frac{cx_{21}x_{22}+x_{32}x_{11}}{x_{21}x_{22}+x_{32}x_{11}}$$
Now let $c_3=x_{32}^2x_{31}x_{33}+x_{22}^2x_{31}x_{32}+x_{22}^2x_{21}^2$, $c_{31}=cx_{32}^2x_{31}x_{33}+x_{22}^2x_{31}x_{32}+x_{22}^2x_{21}^2$ and $c_{32}=cx_{32}^2x_{31}x_{33}+cx_{22}^2x_{31}x_{32}+x_{22}^2x_{21}^2$. Then
$$e_3^c(x_{33},x_{22},x_{32},x_{11},x_{21},x_{31})=\left(\frac{c_{31}}{c_3}x_{33},x_{22},\frac{c_{32}}{c_{31}}x_{32},x_{11},x_{21},\frac{c\cdot c_{3}}{c_{32}}x_{31}\right)$$
Then we compute the $\varepsilon_i$ and $\gamma_i$ actions for $i=1,2,3$
\begin{align*}
    &\varepsilon_1(x_{33},x_{22},x_{32},x_{11},x_{21},x_{31})=\frac{x_{22}}{x_{11}}\\
    &\varepsilon_2(x_{33},x_{22},x_{32},x_{11},x_{21},x_{31})=\frac{x_{33}x_{22}x_{21}+x_{33}x_{32}x_{11}}{x_{22}^2x_{21}}\\
    &\varepsilon_3(x_{33},x_{22},x_{32},x_{11},x_{21},x_{31})=\frac{x_{33}x_{32}^2x_{31}+x_{22}^2x_{32}x_{31}+x_{22}^2x_{21}^2}{x_{33}^2x_{32}^2x_{31}}\\
    &\gamma_1(x_{33},x_{22},x_{32},x_{11},x_{21},x_{31})=\frac{x_{11}^2}{x_{21}x_{22}}\\
    &\gamma_2(x_{33},x_{22},x_{32},x_{11},x_{21},x_{31})=\frac{x_{21}^2x_{22}^2}{x_{31}x_{32}x_{33}x_{11}}\\
    &\gamma_3(x_{33},x_{22},x_{32},x_{11},x_{21},x_{31})=\frac{x_{31}^2x_{32}^2x_{33}^2}{x_{21}^2x_{22}^2}
\end{align*}

To get the formulas for $e_0^c$, $\varepsilon_0$ and $\gamma_0$ we first must get the formulas for $\overline{e_0^c}$, $\overline{\varepsilon_0}$ and $\overline{\gamma_0}$ in $\mathcal{V}_2(y)$. These are as follows:
$$\overline{\varepsilon_0}(y_{03},y_{12},y_{02},y_{21},y_{11},y_{01})==\frac{y_{03}y_{02}^2y_{01}+y_{12}^2y_{02}y_{01}+y_{12}^2y_{11}^2}{y_{03}^2y_{02}^2y_{01}}$$
$$\overline{\gamma_0}(y_{03},y_{12},y_{02},y_{21},y_{11},y_{01})=\frac{y_{01}^2y_{02}^2y_{03}^2}{y_{11}^2y_{12}^2}$$
Let $c_0=y_{03}y_{02}^2y_{01}+y_{12}^2y_{02}y_{01}+y_{12}^2y_{11}^2$, $c_{01}=cy_{03}y_{02}^2y_{01}+y_{12}^2y_{02}y_{01}+y_{12}^2y_{11}^2$, $c_{02}=cy_{03}y_{02}^2y_{01}+cy_{12}^2y_{02}y_{01}+y_{12}^2y_{11}^2$.
Then $$\overline{e_0^c}(y_{03},y_{12},y_{02},y_{21},y_{11},y_{01})=\left(\frac{c_{01}}{c_0}y_{03},y_{12},\frac{c_{02}}{c_{01}}y_{02},y_{21},y_{11},\frac{c\cdot c_0}{c_{02}}y_{01}\right)$$

Using the formulas given previously in 3.1, we obtain $\gamma_0$, $\varepsilon_0$, and $e_0^c$.
$$\gamma_0(V_1(x))=\overline{\gamma_0}(\overline{\sigma}(V_1(x)))=\frac{1}{x_{11}^2}$$
$$\varepsilon_0(V_1(x))=x_{31}+\frac{(x_{21}x_{22}+x_{11}x_{32})^2}{x_{22}^2x_{32}}+\frac{x_{11}^2}{x_{33}}$$

We define $c_{2}'=x_{31}+\frac{x_{21}^2}{x_{32}}+\frac{2x_{21}x_{11}}{x_{22}}+\frac{x_{11}^2x_{32}}{x_{22}^2}+\frac{x_{11}^2}{x_{33}}$, $c_{21}'=cx_{31}+\frac{x_{21}^2}{x_{32}}+\frac{2x_{21}x_{11}}{x_{22}}+\frac{x_{11}^2x_{32}}{x_{22}^2}+\frac{x_{11}^2}{x_{33}}$, $c_{24}'=cx_{31}+c\frac{x_{21}^2}{x_{32}}+c\frac{2x_{21}x_{11}}{x_{22}}+c\frac{x_{11}^2x_{32}}{x_{22}^2}+\frac{x_{11}^2}{x_{33}}$.
\begin{align*}
    &e_0^c(x_{33},x_{22},x_{32},x_{11},x_{21},x_{31})=(\frac{c_{24}'}{c\cdot c_2'}x_{33},x_{22}\frac{x_{21}x_{22}c_{24}'+x_{11}x_{32}c_{21}'}{(x_{21}x_{22}+x_{11}x_{32})c\cdot c_2'},\frac{c_{21}'}{c_{24}'}x_{32},\\&\frac{x_{11}}{c},x_{21}\frac{(x_{21}x_{22}+x_{11}x_{32})\cdot c_2'}{x_{21}x_{22}c_{24}'+x_{11}x_{32}c_{21}'}, \frac{c_2'}{c_{21}'}x_{31})
\end{align*} 
\begin{theorem}
 $\mathcal{V}_1=(V_1(x),e_i^c,\varepsilon_i,\varphi_i)$ for $i=0,1,2,3$ is a positive geometric crystal associated with $C_3^{(1)}$.
\end{theorem}
\begin{proof}
This is clearly positive as all the coefficients are positive. We already know that $\mathcal{V}_1$ is a geometric crystal without the additionally 0-actions, so we only need to check relations with the 0 index.

\begin{align*}
    1. \,\,\,\,&\gamma_0(e_0^c(V_1(x)))=\frac{1}{\frac{x_{11}^2}{c^2}}=\frac{c^2}{x_{11}^2}=c^2\varepsilon_0(V_1(x))\\
    &\gamma_0(e_1^c(V_1(x)))=\gamma_0(x_{33},x_{22},x_{32},cx_{11},x_{21},x_{31})=\frac{1}{c^2x_{11}^2}=\frac{1}{c^2}\gamma_0(V_1(x))\\
    &\gamma_0(e_2^c(V_1(x))=\gamma_0(x_{33},c_2x_{22},x_{32},x_{11},\frac{c}{c_2}x_{21},x_{31})=\frac{1}{x_{11}^2}=\gamma_0(V_1(x))\\
    &\gamma_0(e_3^c(V_1(x)))=\gamma_0\left(\frac{c_{31}}{c_3}x_{33},x_{22},\frac{c_{32}}{c_{31}}x_{32},x_{11},x_{21},\frac{c\cdot c_{3}}{c_{32}}x_{31}\right)=\frac{1}{x_{11}^2}=\gamma_0(V_1(x))\\
2.\,\,\,\,&\gamma_1(e_0^c(V_1(x)))=\frac{\frac{x_{11}^2}{c^2}}{\frac{x_{21}x_{22}}{c}}=c^{-1}\gamma_1(V_1(x))\\
&\gamma_2(e_0^c(V_1(x)))=\frac{\frac{x_{21}^2x_{22}^2}{c^2}}{\frac{x_{11}x_{31}x_{32}x_{33}}{c^2}}=\gamma_2(V_1(x))\\
&\gamma_3(e_0^c(V_1(x)))=\frac{\frac{x_{31}^2x_{32}^2x_{33}^2}{c^2}}{\frac{x_{21}^2x_{22}^2}{c^2}}=\gamma_3(V_1(x))
\end{align*}
Now we prove $\varepsilon_0(e_0^c(V_1(x)))=c^{-1}\varepsilon_0(V_1(x))$:
\begin{align*}
&\varepsilon_0(e_0^c(V_1(x)))=\overline{\varepsilon_0}\overline{\sigma}\overline{\sigma^{-1}}\overline{e_0}^c\overline{\sigma}(V_1(x))=\overline{\varepsilon_0}\left(\frac{c_{01}}{c_0}y_{03},y_{12},\frac{c_{02}}{c_{01}}y_{02},y_{21},y_{11},\frac{c\cdot c_0}{c_{02}}y_{01}\right)\\&=\frac{(\frac{y_{03}y_{02}^2y_{01}c\cdot c_{02}+y_{12}^2c\cdot c_{0}y_{02}y_{01}}{c_{01}}+y_{12}^2y_{11}^2)c_0}{c\cdot c_{02}y_{01}y_{02}^2y_{03}^2}\\
& =\dfrac{c_0K}{cK}=\dfrac{c_0}{c}=\dfrac{y_{03}y_{02}^2y_{01}+y_{12}^2y_{02}y_{01}+y_{11}^2y_{12}^2}{cy_{03}^2y_{02}^2y_{01}}=c^{-1}\overline{\varepsilon_0}(V_2(y))=c^{-1}\varepsilon_0(V_1(x))
\end{align*}

Where $K=c^2y_{03}^2y_{02}^4y_{01}^2+(c^2+c)y_{12}^2y_{03}y_{02}^3y_{01}^2+2cy_{12}^2y_{11}^2y_{03}y_{02}^2y_{01}$ \\$+(c+1)y_{11}^2y_{12}^4y_{02}y_{01}+cy_{12}^4y_{02}^2y_{01}^2+y_{11}^4y_{12}^4$.
\begin{align*}
    4.\,\,\,\,&\varepsilon_2(e_0^c(V_1(x)))=\frac{x_{33}(x_{21}x_{22}+x_{11}x_{32})}{x_{22}^2x_{21}}=\varepsilon_2(V_1(x))\\
    5.\,\,\,\,&\varepsilon_3(e_0^c(V_1(x)))=\frac{\frac{c_{21}'x_{31}x_{32}^2x_{33}}{c_{24}'\cdot c}+\frac{x_{32}x_{31}x_{22}^2(x_{21}x_{22}c_{24}'+x_{11}x_{32}c_{21}')}{(x_{21}x_{22}+x_{11}x_{32})c_{24}'c_2'c^2}+\frac{x_{21}^2x_{22}^2}{c^2}}{\frac{x_{31}x_{32}^2x_{33}^2}{c^2c_2'}}=\varepsilon_3(V_1(x))\\
   6.\,\,\, &\varepsilon_0(e_2^c(V_1(x)))=\varepsilon_0(x_{33},c_2x_{22},x_{32},x_{11},\frac{c}{c_2}x_{21},x_{31})\\&=x_{31}+\frac{x_{11}^2}{x_{33}}+\frac{(cx_{21}x_{22}+x_{11}x_{32})^2}{\frac{x_{32}x_{22}^2(cx_{21}x_{22}+x_{11}x_{32})^2}{(x_{21}x_{22}+x_{11}x_{32})^2}}\\&=x_{31}+\frac{x_{11}^2}{x_{33}}+\frac{(x_{21}x_{22}+x_{11}x_{32})^2}{x_{22}^2x_{32}}=\varepsilon_0(V_1(x))\\
  7.\,\,\, &\varepsilon_0(e_3^c(V_1(x))) =\varepsilon_0\left(\frac{c_{31}}{c_3}x_{33},x_{22},\frac{c_{32}}{c_{31}}x_{32},x_{11},x_{21},\frac{c\cdot c_{3}}{c_{32}}x_{31}\right)\\&=\frac{c\cdot c_{3}}{c_{32}}x_{31}+\frac{x_{11}^2c_{3}}{x_{33}c_{31}}+\frac{x_{21}^2c_{31}}{c_{32}x_{32}}+\frac{2x_{11}x_{21}}{x_{22}}+\frac{x_{11}^2c_{32}x_{32}}{c_{31}x_{22}^2}
\end{align*}
We consider only the terms that have c in them.
\begin{align*}
  &\frac{c\cdot c_{3}}{c_{32}}x_{31}+\frac{x_{11}^2c_{3}}{x_{33}c_{31}}+\frac{x_{21}^2c_{31}}{c_{32}x_{32}}+\frac{x_{11}^2c_{32}x_{32}}{c_{31}x_{22}^2}\\&= \frac{c c_3c_{31} x_{31}x_{32}x_{33}x_{22}^2+x_{21}^2c_{31}^2x_{33}x_{22}^2+x_{11}^2c_3c_{32}x_{32}x_{22}^2+x_{11}^2c_{32}^2x_{32}^2x_{33}}{c_{31}c_{32}x_{32}x_{33}x_{22}^2} 
\end{align*}
If we expand these terms and factor, we get $x_{31}+\frac{x_{21}^2}{x_{32}}+\frac{x_{11}^2x_{32}}{x_{22}}+\frac{x_{11}^2}{x_{33}}$. Along with the term we ignored earlier,  $\frac{2x_{11}x_{21}}{x_{22}}$, we get $\varepsilon_0(V_1(x))$.

Before proving the last 3 identities, we prove a useful statement: $\overline{\sigma}\circ e_i^c=\overline{e_{3-i}}^c\circ\overline{\sigma}$ we prove this for $1\leq i\leq 2$.
For i=1:
\begin{align*}
    &\overline{\sigma}\circ e_1^c(V_1(x))=\overline{\sigma}(x_{33},x_{22},x_{32},cx_{11},x_{21},x_{31})=(y_{03},y_{12},y_{02},cy_{21},y_{11},y_{01})\\
    &\overline{e_{2}}^c\circ\overline{\sigma}(V_1(x))=\overline{e_{2}}^c(y_{03},y_{12},y_{02},y_{21},y_{11},y_{01})=(y_{03},y_{12},y_{02},cy_{21},y_{11},y_{01})
\end{align*}
For i=2:
\begin{align*}
     &\overline{\sigma}\circ e_2^c(V_1(x))=\overline{\sigma}(x_{33},c_2x_{22},x_{32},x_{11},\frac{c}{c_2}x_{21},x_{31})=(y_{03},c_1'y_{12},y_{02},y_{21},\frac{c}{c_1'}y_{11},y_{01})\\
     &\overline{e_{1}}^c\circ\overline{\sigma}(V_1(x))=\overline{e_{1}}^c(y_{03},y_{12},y_{02},y_{21},y_{11},y_{01})=(y_{03},c_1'y_{12},y_{02},y_{21},\frac{c}{c_1'}y_{11},y_{01})
\end{align*}
Where $c_1'=\frac{cy_{12}y_{11}+y_{02}y_{21}}{y_{12}y_{11}+y_{02}y_{21}}$

Using these identities and the fact that $V_2(y)$ is a geometric crystal for $i=0,1,2$, we can prove that
\begin{align*}
    &e_0^ce_2^d=\overline{\sigma^{-1}}\overline{e_0}^c\overline{\sigma}e_2^d=\overline{\sigma^{-1}}\overline{e_0}^c\overline{e_1}^d\overline{\sigma}=\overline{\sigma^{-1}}\overline{e_1}^d\overline{e_0}^c\overline{\sigma}=e_2^d\overline{\sigma^{-1}}\overline{e_0}^c\overline{\sigma}=e_2^de_0^c\\
   &e_1^c\overline{\sigma^{-1}}\overline{e_0}^{c^2d}\overline{\sigma}e_1^{cd}\overline{\sigma^{-1}}\overline{e_0}^d\overline{\sigma}=\overline{\sigma^{-1}}\overline{e_2}^c\overline{e_0}^{c^2d}\overline{\sigma}\overline{\sigma^{-1}}\overline{e_2}^{cd}\overline{e_0}^d\overline{\sigma}=\overline{\sigma^{-1}}\overline{e_2}^c\overline{e_0}^{c^2d}\overline{e_2}^{cd}\overline{e_0}^d\overline{\sigma}\\&=\overline{\sigma^{-1}}\overline{e_0}^d\overline{e_2}^{cd}\overline{e_0}^{c^2d}\overline{e_2}^c\overline{\sigma}=\overline{\sigma^{-1}}\overline{e_0}^d\overline{e_2}^{cd}\overline{\sigma}\overline{\sigma^{-1}}\overline{e_0}^{c^2d}\overline{\sigma}e_1^c\\&=\overline{\sigma^{-1}}\overline{e_0}^d\overline{\sigma}e_1^{cd}\overline{\sigma^{-1}}\overline{e_0}^{c^2d}\overline{\sigma}e_1^c=e_0^de_1^{cd}e_0^{c^2d}e_1^c
\end{align*}
We used a computer algebra system to prove the final relation, $e_0^ce_4^d=e_4^de_0^c$. Therefore, all relations hold, so it is a geometric crystal.
\end{proof}
\subsection{n=4 Case}
This case is the first case where the representation is not multiplicity free. The two vectors $(3,4,\overline{4},\overline{3})$ and $(2,4,\overline{4},\overline{2})$ both have weight 0. As a result, we used the global base to consider the e and f actions. By direct computation, we see that $f_2(2,4,\overline{4},\overline{3})=(2,4,\overline{4},\overline{2})+(3,4,\overline{4},\overline{3})$. All other actions are the same for the global and crystal bases. Additionally, we note that if $f_i^2\neq 0$, then we have $f_i^2(b)=2b'$. 

 $V_1(x)$ is equal to
\begin{small}
$$Y_{4}(x_{44})Y_{3}(x_{33})Y_4(x_{43})Y_2(x_{22})Y_3(x_{32})Y_4(x_{42})Y_1(x_{11})Y_2(x_{21})Y_3(x_{31})Y_4(x_{41})(1,2,3,4)$$
\end{small}
And $V_2(y)$ is equal to
\begin{small}
$$Y_0(y_{04})Y_1(y_{13})Y_0(y_{03})Y_2(y_{22})Y_1(y_{12})Y_0(y_{02})Y_3(y_{31})Y_2(y_{21})Y_1(y_{11})Y_0(y_{01})(\overline{4},\overline{3},\overline{2},\overline{1})$$
\end{small}
  Setting up the system of equations $V_1(x)\cdot a(x)=V_2(y)$. The solutions to this system are $x_{41}=\dfrac{1}{y_{04}}$, $x_{42}=\dfrac{1}{y_{03}}$, $x_{43}=\dfrac{1}{y_{02}}$, $x_{44}=\dfrac{1}{y_{01}}$,\\ $x_{33}=\dfrac{y_{31}y_{02}y_{12}+y_{21}y_{22}y_{02}+y_{11}y_{12}y_{22}}{y_{12}y_{22}y_{01}y_{02}}$, \\$x_{32}=\dfrac{(y_{31}y_{03}y_{22}+y_{31}y_{12}y_{13}+y_{21}y_{22}y_{13})y_{12}}{y_{03}y_{13}(y_{31}y_{02}y_{12}+y_{21}y_{22}y_{02}+y_{11}y_{12}y_{22})}$,\\ $x_{31}=\dfrac{y_{31}y_{22}y_{13}}{y_{04}(y_{31}y_{03}y_{22}+y_{31}y_{12}y_{13}+y_{21}y_{22}y_{13})}$, 
  \\$x_{22}=\dfrac{y_{21}y_{22}y_{02}y_{03}+y_{11}y_{12}y_{22}y_{03}+y_{11}y_{12}^2y_{13})}{y_{12}y_{13}y_{01}y_{02}y_{03}}$, \\$x_{21}=\dfrac{y_{21}y_{22}y_{12}y_{13}}{y_{04}(y_{21}y_{22}y_{02}y_{03}+y_{11}y_{12}y_{22}y_{03}+y_{11}y_{12}^2y_{13})}$, and $x_{11}=\dfrac{y_{11}y_{12}y_{13}}{y_{01}y_{02}y_{03}y_{04}}$.
  We can also see that $x_{41}=\frac{1}{y_{04}}$, $x_{42}=\frac{1}{y_{03}}$, $x_{43}=\frac{1}{y_{02}}$, $x_{44}=\frac{1}{y_{01}}$,\\ $x_{33}=\frac{y_{31}y_{02}y_{12}+y_{21}y_{22}y_{02}+y_{11}y_{12}y_{22}}{y_{12}y_{22}y_{01}y_{02}}$, $x_{32}=\frac{(y_{31}y_{03}y_{22}+y_{31}y_{12}y_{13}+y_{21}y_{22}y_{13})y_{12}}{y_{03}y_{13}(y_{31}y_{02}y_{12}+y_{21}y_{22}y_{02}+y_{11}y_{12}y_{22})}$, $x_{31}=\frac{y_{31}y_{22}y_{13}}{y_{04}(y_{31}y_{03}y_{22}+y_{31}y_{12}y_{13}+y_{21}y_{22}y_{13})}$, $x_{22}=\frac{y_{21}y_{22}y_{02}y_{03}+y_{11}y_{12}y_{22}y_{03}+y_{11}y_{12}^2y_{13})}{y_{12}y_{13}y_{01}y_{02}y_{03}}$,\\ $x_{21}=\frac{y_{21}y_{22}y_{12}y_{13}}{y_{04}(y_{21}y_{22}y_{02}y_{03}+y_{11}y_{12}y_{22}y_{03}+y_{11}y_{12}^2y_{13})}$, and $x_{11}=\frac{y_{11}y_{12}y_{13}}{y_{01}y_{02}y_{03}y_{04}}$.
  
Now we need to compute the actions of the geometric crystal. First we compute $e_1^c$, $e_2^c$, $e_3^c$, and $e_4^c$. \begin{align*}
&e_1^c(V_1(x))=(x_{44},x_{33},x_{43},x_{22},x_{32},x_{42},cx_{11},x_{21},x_{31},x_{41})\\
    &e_2^c(V_1(x))=(x_{44},x_{33},x_{43},\frac{c_{21}}{c_2}x_{22},x_{32},x_{42},x_{11},\frac{cc_2}{c_{21}}x_{21},x_{31},x_{41})
\end{align*}
where $\dfrac{c_{21}}{c_2}=\dfrac{\dfrac{cx_{33}x_{22}x_{21}+x_{32}x_{33}x_{11}}{x_{22}^2x_{21}}}{\dfrac{x_{33}x_{22}x_{21}+x_{32}x_{33}x_{11}}{x_{22}^2x_{21}}}=\dfrac{cx_{22}x_{21}+x_{32}x_{11}}{x_{22}x_{21}+x_{32}x_{11}}$

Let $c_3=x_{44}x_{33}x_{32}^2x_{31}+x_{44}x_{43}x_{22}x_{32}x_{31}+x_{44}x_{43}x_{42}x_{22}x_{21}$,\\ $c_{31}=cx_{44}x_{33}x_{32}^2x_{31}+x_{44}x_{43}x_{22}x_{32}x_{31}+x_{44}x_{43}x_{42}x_{22}x_{21}$ \\ and $c_{32}=cx_{44}x_{33}x_{32}^2x_{31}+cx_{44}x_{43}x_{22}x_{32}x_{31}+x_{44}x_{43}x_{42}x_{22}x_{21}$. Then \begin{align*}&e_3^c(V_1(x))=(x_{44},\frac{c_{31}}{c_3}x_{33},x_{43},x_{22},\frac{c_{32}}{c_{31}}x_{32},x_{42},x_{11},x_{21},\frac{cc_3}{c_{32}}x_{31},x_{41})\end{align*}
\\ Let $c_4=x_{41}x_{42}^2x_{43}^2x_{44}+x_{33}^2x_{43}x_{42}^2x_{41}+x_{33}^2x_{32}^2x_{42}x_{41}+x_{33}^2x_{32}^2x_{31}^2$,\\ $c_{41}=cx_{41}x_{42}^2x_{43}^2x_{44}+x_{33}^2x_{43}x_{42}^2x_{41}+x_{33}^2x_{32}^2x_{42}x_{41}+x_{33}^2x_{32}^2x_{31}^2$, \\$c_{42}=cx_{41}x_{42}^2x_{43}^2x_{44}+cx_{33}^2x_{43}x_{42}^2x_{41}+x_{33}^2x_{32}^2x_{42}x_{41}+x_{33}^2x_{32}^2x_{31}^2$,\\ and $c_{43}=cx_{41}x_{42}^2x_{43}^2x_{44}+cx_{33}^2x_{43}x_{42}^2x_{41}+cx_{33}^2x_{32}^2x_{42}x_{41}+x_{33}^2x_{32}^2x_{31}^2$. Then 

\begin{align*}&e_4^c(V_1(x))=\left(\frac{c_{41}}{c_{4}}x_{44},x_{33},\frac{c_{42}}{c_{41}}x_{43},x_{22},x_{32},\frac{c_{43}}{c_{42}}x_{42},x_{11},x_{21},x_{31},\frac{cc_4}{c_{43}}x_{41}\right)\end{align*}
Now we define $\varepsilon_i$ and $\gamma_i$ for $i=1..4$.
\begin{align*}
    &\varepsilon_1(V_1(x))=\frac{x_{22}}{x_{11}}\\
    &\varepsilon_2(V_1(x))=\frac{x_{33}}{x_{22}}+\frac{x_{33}x_{32}x_{11}}{x_{22}^2x_{21}}\\&
    \varepsilon_3(V_1(x))=\frac{x_{44}}{x_{33}}+\frac{x_{44}x_{43}x_{22}}{x_{33}^2x_{32}}+\frac{x_{44}x_{43}x_{42}x_{22}x_{21}}{x_{33}^2x_{32}^2x_{31}}\\
    &\varepsilon_4(V_1(x))=\frac{1}{x_{44}}+\frac{x_{33}^2}{x_{43}x_{44}^2}+\frac{x_{33}^2x_{32}^2}{x_{44}^2x_{43}^2x_{42}}+\frac{x_{33}^2x_{32}^2x_{31}^2}{x_{44}^2x_{43}^2x_{42}^2x_{41}}\\
    &\gamma_1(V_1(x))=\frac{x_{11}^2}{x_{21}x_{22}}\\&
    \gamma_2(V_1(x))=\frac{x_{21}^2x_{22}^2}{x_{11}x_{31}x_{32}x_{33}}\\&
    \gamma_3(V_1(x))=\frac{x_{31}^2x_{32}^2x_{33}^2}{x_{21}x_{22}x_{41}x_{42}x_{43}x_{44}}\\&
    \gamma_4(V_1(x))=\frac{x_{41}^2x_{42}^2x_{43}^2x_{44}^2}{x_{31}^2x_{32}^2x_{33}^2}
\end{align*}
To get the formulas for $e_0^c$, $\varepsilon_0$ and $\gamma_0$ we first must get the formulas for $\overline{e_0^c}$, $\overline{\varepsilon_0}$ and $\overline{\gamma_0}$ in $\mathcal{V}_2(y)$. These are as follows:
\begin{align*}&\overline{e_0^c}(V_2(y))=(\frac{c_{01}}{c_0}y_{04},y_{13},\frac{c_{02}}{c_{01}}y_{03},y_{22},y_{12},\frac{c_{03}}{c_{02}}y_{02},y_{31},y_{21},y_{11},\frac{cc_0}{c_{03}}y_{01})\end{align*}
where $c_0=\frac{1}{y_{04}}+\frac{y_{13}^2}{y_{04}^2y_{03}}+\frac{y_{13}^2y_{12}^2}{y_{04}^2y_{03}^2y_{02}}+\frac{y_{13}^2y_{12}^2y_{11}^2}{y_{04}^2y_{03}^2y_{02}^2y_{01}}$, $c_{01}=\frac{c}{y_{04}}+\frac{y_{13}^2}{y_{04}^2y_{03}}+\frac{y_{13}^2y_{12}^2}{y_{04}^2y_{03}^2y_{02}}+\frac{y_{13}^2y_{12}^2y_{11}^2}{y_{04}^2y_{03}^2y_{02}^2y_{01}}$, $c_{02}=\frac{c}{y_{04}}+\frac{cy_{13}^2}{y_{04}^2y_{03}}+\frac{y_{13}^2y_{12}^2}{y_{04}^2y_{03}^2y_{02}}+\frac{y_{13}^2y_{12}^2y_{11}^2}{y_{04}^2y_{03}^2y_{02}^2y_{01}}$ and $c_{03}=\frac{c}{y_{04}}+\frac{cy_{13}^2}{y_{04}^2y_{03}}+\frac{cy_{13}^2y_{12}^2}{y_{04}^2y_{03}^2y_{02}}+\frac{y_{13}^2y_{12}^2y_{11}^2}{y_{04}^2y_{03}^2y_{02}^2y_{01}}$.
$$\overline{\varepsilon_0}(V_2(y))=\frac{1}{y_{04}}+\frac{y_{13}^2}{y_{04}^2y_{03}}+\frac{y_{13}^2y_{12}^2}{y_{04}^2y_{03}^2y_{02}}+\frac{y_{13}^2y_{12}^2y_{11}^2}{y_{04}^2y_{03}^2y_{02}^2y_{01}}$$
$$\overline{\gamma_0}(V_2(y))=\frac{y_{01}^2y_{02}^2y_{03}^2y_{04}^2}{y_{11}^2y_{12}^2y_{13}^2}$$

Using the formulas given previously, we obtain $\gamma_0$, $\varepsilon_0$, and $e_0^c$:

\begin{align*}
&\gamma_0(V_1(x))=\frac{y_{01}^2y_{02}^2y_{03}^2y_{04}^2}{y_{11}^2y_{12}^2y_{13}^2}=\frac{1}{x_{11}^2}\\
&\varepsilon_0(V_1(x))=\frac{1}{y_{04}}+\frac{y_{13}^2}{y_{04}^2y_{03}}+\frac{y_{13}^2y_{12}^2}{y_{04}^2y_{03}^2y_{02}}+\frac{y_{13}^2y_{12}^2y_{11}^2}{y_{04}^2y_{03}^2y_{02}^2y_{01}}\\&=x_{41}+\frac{(\frac{x_{11}x_{42}}{x_{22}}+\frac{x_{21}x_{42}}{x_{32}}+x_{31})^2}{x_{42}}+\frac{(x_{21}+\frac{x_{11}x_{43}}{x_{33}}+\frac{x_{11}x_{32}}{x_{22}})^2}{x_{43}}+\frac{x_{11}^2}{x_{44}}
\end{align*}
Let \begin{align*}
    &A=x_{11}x_{22}x_{43}c_{42}'+c_{43}'x_{21}x_{22}x_{33}+c_{43}'x_{11}x_{32}x_{33}\\
    &B=x_{11}x_{22}x_{43}+x_{21}x_{22}x_{33}+x_{11}x_{32}x_{33}\\
    &C=c_{41}'x_{11}x_{32}x_{42}+c_{41}'x_{21}x_{22}x_{42}+c_{42}'x_{22}x_{31}x_{32}\\
    &D=x_{11}x_{32}x_{42}+x_{21}x_{22}x_{42}+x_{22}x_{31}x_{32}\\
    &G=c_{41}'(x_{11}^2x_{22}x_{32}x_{42}x_{43}+x_{11}^2x_{32}^2x_{33}x_{42}+x_{11}x_{21}x_{22}x_{32}x_{33}x_{42})\\&+c_{42}'(x_{11}x_{21}x_{22}^2x_{42}x_{43}+x_{11}x_{22}^2x_{31}x_{32}x_{43}+x_{11}x_{22}x_{31}x_{32}^2x_{33})\\&+c_{43}'(x_{11}x_{21}x_{22}x_{32}x_{33}x_{42}+x_{21}^2x_{22}^2x_{33}x_{42}+x_{21}x_{22}^2x_{31}x_{32}x_{33})
\end{align*}
\begin{align*}
    &e_0^c(V_1(x))\\&= (x_{44}\frac{c_{43}'}{cc_4'},\frac{x_{33}A}{cc_4'B},\frac{x_{43}c_{42}'}{c_{43}'},\frac{x_{22}G}{BDcc_4'},\frac{x_{32}CB}{AD}, \frac{c_{41}'x_{42}}{c_{42}'},\frac{x_{11}}{c},\frac{x_{21}c_4'DB}{G}, \frac{x_{31}c_4'D}{C}, \frac{x_{41}c_4'}{c_{41}'})
\end{align*}
where $c_4'=x_{41}+\frac{(\frac{x_{11}x_{42}}{x_{22}}+\frac{x_{21}x_{42}}{x_{32}}+x_{31})^2}{x_{42}}+\frac{(x_{21}+\frac{x_{11}x_{43}}{x_{33}}+\frac{x_{11}x_{32}}{x_{22})})^2}{x_{43}}+\frac{x_{11}^2}{x_{44}}$, $c_{41}'=cx_{41}+\frac{(\frac{x_{11}x_{42}}{x_{22}}+\frac{x_{21}x_{42}}{x_{32}}+x_{31})^2}{x_{42}}+\frac{(x_{21}+\frac{x_{11}x_{43}}{x_{33}}+\frac{x_{11}x_{32}}{x_{22}})^2}{x_{43}}+\frac{x_{11}^2}{x_{44}}$, $c_{42}'=cx_{41}+\frac{c(\frac{x_{11}x_{42}}{x_{22}}+\frac{x_{21}x_{42}}{x_{32}}+x_{31})^2}{x_{42}}+\frac{(x_{21}+\frac{x_{11}x_{43}}{x_{33}}+\frac{x_{11}x_{32}}{x_{22}})^2}{x_{43}}+\frac{x_{11}^2}{x_{44}}$,\\ and $c_{43}'=cx_{41}+\frac{c(\frac{x_{11}x_{42}}{x_{22}}+\frac{x_{21}x_{42}}{x_{32}}+x_{31})^2}{x_{42}}+\frac{c(x_{21}+\frac{x_{11}x_{43}}{x_{33}}+\frac{x_{11}x_{32}}{x_{22}})^2}{x_{43}}+\frac{x_{11}^2}{x_{44}}$.

\begin{theorem}
 $\mathcal{V}_1=(V_1(x),e_i^c,\varepsilon_i,\varphi_i)$ for $i=0,1,2,3,4$ is a positive geometric crystal associated with $C_4^{(1)}$.
\end{theorem}
\begin{proof}
We know this is positive because all of the coefficients are positive. We already know that $\mathcal{V}_1$ is a geometric crystal without the additionally 0-actions, so we only need to check relations with the 0 index.

We first check the relation $\gamma_0(e_i^c(V_1(x)))=c^{a_{i0}}\gamma_0(V_1(x))$
\begin{itemize}
    \item $\gamma_0(e_0^c(V_1(x)))=\frac{c^2}{x_{11}^2}=c^2\gamma_0(V_1(x))$
    \item $\gamma_0(e_1^c(V_1(x)))=\frac{1}{c^2x_{11}^2}=c^{-2}\gamma_0(V_1(x))$
    \item $\gamma_0(e_2^c(V_1(x)))=\frac{1}{x_{11}^2}=\gamma_0(V_1(x))$
    \item $\gamma_0(e_3^c(V_1(x)))=\frac{1}{x_{11}^2}=\gamma_0(V_1(x))$
    \item $\gamma_0(e_4^c(V_1(x)))=\frac{1}{x_{11}^2}=\gamma_0(V_1(x))$
\end{itemize}
 Next we check that $\gamma_i(e_0^c(V_1(x)))=c^{a_{0i}}\gamma_i(V_1(x))$.
\begin{itemize}
    \item $\gamma_1(e_0^c(V_1(x)))=\dfrac{\frac{x_{11}^2}{c^2}}{\frac{x_{21}x_{22}}{c}}=c^{-1}\gamma_0(V_1(x))$
    \item $\gamma_2(e_0^c(V_1(x)))=\dfrac{\frac{x_{21}^2x_{22}^2}{c^2}}{\frac{x_{11}x_{31}x_{32}x_{33}}{c^2}}=\gamma_2(V_1(x))$
    \item $\gamma_3(e_0^c(V_1(x)))=\dfrac{\frac{x_{31}^2x_{32}^2x_{33}^2}{c^2}}{\frac{x_{11}x_{41}x_{42}x_{43}x_{44}}{c^2}}=\gamma_3(V_1(x))$
    \item $\gamma_4(e_0^c(V_1(x)))=\dfrac{\frac{x_{41}^2x_{42}^2x_{43}^2x_{44}^2}{c^2}}{\frac{x_{31}^2x_{32}^2x_{33}^2}{c^2}}=\gamma_4(V_1(x))$
\end{itemize}
Now we check $\varepsilon_0(e_0^c(V_1(x)))=c^{-1}\varepsilon_0(V_1(x))$. We can see that this relation is equivalent to another relation as follows: $\varepsilon_0(e_0^c(V_1(x)))=\overline{\varepsilon}_0\circ\overline{\sigma}\circ \overline{\sigma^{-1}}\circ\overline{e}_0^c\circ\overline{\sigma}=\overline{\varepsilon}_0\circ\overline{e}_0^c\circ\overline{\sigma}(V_1(x))=\overline{\varepsilon}_0\circ\overline{e}_0^c(V_2(y))$
So we can just consider $\overline{\varepsilon}_0\circ\overline{e}_0^c(V_2(y))$. 
\begin{align*}
 &\overline{\varepsilon_0}(\frac{c_{01}y_{04}}{c_0},y_{13},\frac{c_{02}y_{03}}{c_{01}},y_{22},y_{12},\frac{c_{03}y_{02}}{c_{02}},y_{31},y_{21},y_{11},\frac{cc_0y_{01}}{c_{03}})\\&=\dfrac{\frac{c_{03}}{c_{01}}y_{04}y_{03}^2y_{02}^2y_{01}+\frac{cc_0c_{03}}{c_{01}c_{02}}y_{13}^2y_{03}y_{02}^2y_{01}+\frac{cc_0}{c_{02}}y_{13}^2y_{12}^2y_{02}y_{01}+y_{11}^2y_{12}^2y_{13}^2}{\frac{cc_{03}}{c_0}y_{04}^2y_{03}^2y_{02}^2y_{01}}
\end{align*}
If we expand the $c_{01},c_{02}, c_{03}$ terms and factor, we obtain $$=\dfrac{y_{04}y_{03}^2y_{02}^2y_{01}+y_{13}^2y_{03}y_{02}^2y_{01}+y_{13}^2y_{12}^2y_{02}y_{01}+y_{11}^2y_{12}^2y_{13}^2}{y_{04}^2y_{03}^2y_{02}^2y_{01}}=\overline{\varepsilon}_0(V_2(y))$$
Next we check $\varepsilon_2(e_0^c(V_1(x)))=\varepsilon_2(V_1(x))$. 

\begin{align*}
 &\varepsilon_2(e_0^c(V_1(x)))=\dfrac{\frac{x_{33}x_{22}x_{21}A}{Bc^2c_4'}+\frac{x_{11}x_{32}x_{33}C}{c^2c_4'D}}{\frac{x_{21}x_{22}^2G}{c^2c_4'BD}}=   \dfrac{x_{33}x_{22}x_{21}AD+x_{11}x_{32}x_{33}BC}{x_{21}x_{22}^2G}\\
 &=\dfrac{x_{33}(x_{21}x_{22}+x_{11}x_{32})G}{Gx_{22}^2x_{21}}=\varepsilon_2(V_1(x))
\end{align*}

Now we check the relation $\varepsilon_3(e_0^c(V_1(x)))=\varepsilon_3(V_1(x))$. 
\begin{align*}
    \varepsilon_3(e_0^c(V_1(x)))=&\dfrac{\frac{x_{44}c_{43}'x_{33}x_{32}^2x_{31}CB}{c^2c_4'AD}+\frac{x_{44}x_{43}c_{42}'x_{32}x_{31}c_4'Bx_{22}G}{cc_4'ABDcc_4'}+\frac{x_{44}x_{43}x_{42}c_{41}'x_{22}x_{21}}{c^2c_4'}}{\frac{x_{31}x_{32}^2x_{33}^2C}{c^2c_4'D}}\\
    &=\dfrac{\frac{x_{44}c_{43}'x_{33}x_{32}^2x_{31}CB}{AD}+\frac{x_{44}x_{43}c_{42}'x_{32}x_{31}x_{22}G}{AD}+x_{44}x_{43}x_{42}c_{41}'x_{22}x_{21}}{\frac{x_{31}x_{32}^2x_{33}^2C}{D}}
\end{align*}
Expanding the first two terms in the numerator and factoring, we get
\begin{align*}
    &\dfrac{x_{44}c_{43}'x_{33}x_{32}^2x_{31}CB+x_{44}x_{43}c_{42}'x_{32}x_{31}x_{22}G}{AD}\\
    &=\dfrac{c_{41}'(x_{11}x_{22}x_{31}x_{32}^2x_{42}x_{43}x_{44}+x_{11}x_{31}x_{32}^3x_{33}x_{42}x_{44}+x_{21}x_{22}x_{31}x_{32}^2x_{33}x_{42}x_{44})}{D}\\
    &+\dfrac{c_{42}'(x_{21}x_{22}^2x_{31}x_{32}x_{42}x_{43}x_{44}+x_{22}^2x_{31}^2x_{32}^2x_{43}x_{44}+x_{22}x_{31}^2x_{32}^3x_{33}x_{44})}{D}
\end{align*}
Putting this expression back in the earlier expression we get
\begin{align*}
    &\dfrac{c_{41}'x_{11}x_{22}x_{31}x_{32}^2x_{42}x_{43}x_{44}+c_{41}'x_{11}x_{31}x_{32}^3x_{33}x_{42}x_{44}+c_{41}'x_{21}x_{22}x_{31}x_{32}^2x_{33}x_{42}x_{44}}{Cx_{31}x_{32}^2x_{33}^2}\\
    &+\dfrac{c_{42}'x_{21}x_{22}^2x_{31}x_{32}x_{42}x_{43}x_{44}+c_{42}'x_{22}^2x_{31}^2x_{32}^2x_{43}x_{44}+c_{42}'x_{22}x_{31}^2x_{32}^3x_{33}x_{44}}{Cx_{31}x_{32}^2x_{33}^2}\\
    &+\dfrac{x_{44}x_{43}x_{42}x_{22}x_{21}c_{41}'(x_{11}x_{32}x_{42}+x_{21}x_{22}x_{42}+x_{22}x_{31}x_{32})}{Cx_{31}x_{32}^2x_{33}^2}\\
    &=\dfrac{C(x_{22}x_{31}x_{32}x_{43}x_{44}+x_{31}x_{32}^2x_{33}x_{44}+x_{21}x_{22}x_{42}x_{43}x_{44})}{Cx_{31}x_{32}^2x_{33}^2}=\varepsilon_3(V_1(x))
\end{align*}
Next we check $\varepsilon_4(e_0^c(V_1(x))=\varepsilon_4(V_1(x))$.
\begin{align*}
    &\varepsilon_4(e_0^c(V_1(x)))=\dfrac{\frac{x_{44}x_{43}^2x_{42}^2x_{41}c_{41}'}{cc_{43}'}+\frac{x_{33}^2A^2x_{43}x_{42}^2x_{41}c_{41}'c_4'}{c^2c_4'^2B^2c_{43}'c_{42}'}+\frac{x_{33}^2x_{32}^2C^2x_{41}x_{42}c_4'}{c^2c_4'^2D^2c_{42}'}+\frac{x_{31}^2x_{32}^2x_{33}^2}{c^2}}{\frac{x_{44}^2x_{43}^2x_{42}^2x_{41}c_{41}'}{c^2c_4'}}\\
    &=\dfrac{\frac{x_{44}x_{43}^2x_{42}^2x_{41}c_{41}'}{c_{43}'}+\frac{x_{33}^2A^2x_{43}x_{42}^2x_{41}c_{41}'}{B^2c_{43}'c_{42}'c}+\frac{x_{33}^2x_{32}^2C^2x_{41}x_{42}}{cc_4'D^2c_{42}'}+\frac{x_{31}^2x_{32}^2x_{33}^2}{c}}{\frac{x_{44}^2x_{43}^2x_{42}^2x_{41}c_{41}'}{cc_4'}}=\varepsilon_4(V_1(x))
\end{align*}
We obtain $\varepsilon_4(V_1(x))$ by expanding the terms and factoring the expanded rational function.

Now we check $\varepsilon_0(e_2^c(V_1(x)))=\varepsilon_0(V_1(x))$.
\begin{align*}
   & \varepsilon_0(e_2^c(V_1(x)))=x_{41}+\dfrac{(\frac{c_2x_{11}x_{42}}{c_{21}x_{22}}+\frac{x_{21}x_{42}cc_2}{c_{21}x_{32}}+x_{31})^2}{x_{42}}\\
    &+\dfrac{(\frac{cc_2x_{21}}{c_{21}}+\frac{x_{11}x_{43}}{x_{33}}+\frac{c_2x_{11}x_{32}}{c_{21}x_{22}})^2}{x_{43}}+\dfrac{x_{11}^2}{x_{44}}\\
\end{align*}
We focus on the terms $\frac{c_2x_{11}x_{42}}{c_{21}x_{22}}+\frac{x_{21}x_{42}cc_2}{c_{21}x_{32}}+x_{31}$ and $\frac{cc_2x_{21}}{c_{21}}+\frac{x_{11}x_{43}}{x_{33}}+\frac{c_2x_{11}x_{32}}{c_{21}x_{22}}$. 

Then we can simplify the first expression as follows:
\begin{align*}
  &\frac{c_2x_{11}x_{42}}{c_{21}x_{22}}+\frac{x_{21}x_{42}cc_2}{c_{21}x_{32}}+x_{31}= \dfrac{c_2x_{11}x_{32}x_{42}+cc_2x_{21}x_{22}x_{42}+c_{21}x_{22}x_{31}x_{32}}{c_{21}x_{22}x_{32}}\\
  &=\dfrac{(cx_{21}x_{22}+x_{11}x_{32})(x_{21}x_{22}x_{42}+x_{11}x_{32}x_{42}+x_{22}x_{31}x_{32})}{(cx_{21}x_{22}+x_{11}x_{32})x_{22}x_{32}}\\
  &=\frac{x_{11}x_{42}}{x_{22}}+\frac{x_{21}x_{42}}{x_{32}}+x_{31}
\end{align*}
And the second expression as follows:
\begin{align*}
    &\frac{cc_2x_{21}}{c_{21}}+\frac{x_{11}x_{43}}{x_{33}}+\frac{c_2x_{11}x_{32}}{c_{21}x_{22}}=\dfrac{cc_2x_{21}x_{22}x_{33}+x_{11}x_{22}x_{43}c_{21}+c_2x_{11}x_{32}x_{33}}{c_{21}x_{22}x_{33}}\\
    &=\dfrac{(cx_{21}x_{22}+x_{11}x_{32})(x_{21}x_{22}x_{33}+x_{11}x_{32}x_{33}+x_{11}x_{22}x_{43})}{(cx_{21}x_{22}+x_{11}x_{32})x_{22}x_{33}}\\
    &=x_{21}+\frac{x_{11}x_{43}}{x_{22}}+\frac{x_{11}x_{32}}{x_{22}}
\end{align*}
Substituting these expressions into our earlier computation we obtain
\begin{align*}
    =x_{41}+\frac{(\frac{x_{11}x_{42}}{x_{22}}+\frac{x_{21}x_{42}}{x_{32}}+x_{31})^2}{x_{42}}+\frac{(x_{21}+\frac{x_{11}x_{43}}{x_{33}}+\frac{x_{11}x_{32}}{x_{22}})^2}{x_{43}}+\frac{x_{11}^2}{x_{44}}=\varepsilon_0(V_1(x))
\end{align*}
Next we check the relation $\varepsilon_0(e_3^c(V_1(x)))=\varepsilon_0(V_1(x))$.
\begin{align*}
    &\varepsilon_0(e_3^c(V_1(x)))=x_{41}+\dfrac{(\frac{x_{11}x_{42}}{x_{22}}+\frac{c_{31}x_{21}x_{42}}{x_{32}c_{32}}+\frac{cc_3x_{31}}{c_{32}})^2}{x_{42}}\\&+\dfrac{(x_{21}+\frac{x_{11}x_{43}c_3}{c_{31}x_{33}}+\frac{x_{11}x_{32}c_{32}}{x_{22}c_{31}})^2}{x_{43}}+\dfrac{x_{11}^2}{x_{44}}
\end{align*}
 The first and last terms are already equal to the terms in $\varepsilon_0(V_1(x))$, so we focus on the middle two terms. We will expand and simplify $\frac{x_{11}x_{42}}{x_{22}}+\frac{c_{31}x_{21}x_{42}}{x_{32}c_{32}}+\frac{cc_3x_{31}}{c_{32}}$ and $x_{21}+\frac{x_{11}x_{43}c_3}{c_{31}x_{33}}+\frac{x_{11}x_{32}c_{32}}{x_{22}c_{31}}$. We can look at only the components of these expressions that contain $c_{3}, c_{31}$, or $c_{32}$. Thus we will rewrite and simplify $\frac{c_{31}x_{21}x_{42}}{x_{32}c_{32}}+\frac{cc_3x_{31}}{c_{32}}$ and $\frac{x_{11}x_{43}c_3}{c_{31}x_{33}}+\frac{x_{11}x_{32}c_{32}}{x_{22}c_{31}}$. The first expression simplifies as follows:
 
 \begin{align*}
    &\frac{c_{31}x_{21}x_{42}}{x_{32}c_{32}}+\frac{cc_3x_{31}}{c_{32}}=\dfrac{c_{31}x_{21}x_{42}+cc_3x_{31}x_{32}}{x_{32}c_{32}}\\
    &=\dfrac{(x_{21}x_{42}+x_{31}x_{32})(cx_{44}x_{33}x_{32}^2x_{31}+cx_{44}x_{43}x_{22}x_{32}x_{31}+x_{44}x_{43}x_{42}x_{22}x_{21})}{(cx_{44}x_{33}x_{32}^2x_{31}+cx_{44}x_{43}x_{22}x_{32}x_{31}+x_{44}x_{43}x_{42}x_{22}x_{21})x_{32}}\\&=\frac{x_{21}x_{42}}{x_{32}}+x_{31}
 \end{align*}
And the second simplifies as follows:
\begin{align*}
    &\frac{x_{11}x_{43}c_3}{c_{31}x_{33}}+\frac{x_{11}x_{32}c_{32}}{x_{22}c_{31}}=\dfrac{x_{11}}{x_{33}x_{22}}+\dfrac{x_{43}x_{22}c_3+x_{32}x_{33}c_{32}}{c_{31}}\\
    &=\dfrac{x_{11}(x_{22}x_{43}+x_{32}x_{33})(cx_{44}x_{33}x_{32}^2x_{31}+x_{44}x_{43}x_{22}x_{32}x_{31}+x_{44}x_{43}x_{42}x_{22}x_{21})}{x_{32}x_{33}(cx_{44}x_{33}x_{32}^2x_{31}+x_{44}x_{43}x_{22}x_{32}x_{31}+x_{44}x_{43}x_{42}x_{22}x_{21})}\\
    &=\frac{x_{11}x_{43}}{x_{33}}+\frac{x_{32}x_{11}}{x_{22}}
\end{align*}

Substituting these reduced expressions into $\varepsilon_0(e_3^c)$ we get:
\begin{align*}
    &x_{41}+\dfrac{(\frac{x_{11}x_{42}}{x_{22}}+\frac{c_{31}x_{21}x_{42}}{x_{32}c_{32}}+\frac{cc_3x_{31}}{c_{32}})^2}{x_{42}}+\dfrac{(x_{21}+\frac{x_{11}x_{43}c_3}{c_{31}x_{33}}+\frac{x_{11}x_{32}c_{32}}{x_{22}c_{31}})^2}{x_{43}}+\dfrac{x_{11}^2}{x_{44}}\\
    &= x_{41}+\dfrac{(\frac{x_{11}x_{42}}{x_{22}}+\frac{x_{21}x_{42}}{x_{32}}+x_{31})^2}{x_{42}}+\dfrac{(x_{21}+\frac{x_{11}x_{43}}{x_{33}}+\frac{x_{11}x_{32}}{x_{22}})^2}{x_{43}}+\dfrac{x_{11}^2}{x_{44}}=\varepsilon_0(V_1(x))
\end{align*}
 Next we prove $\varepsilon_0(e_4^c(V_1(x)))=\varepsilon_0(V_1(x))$.
\begin{align*}
    &\varepsilon_0(e_4^c(V_1(x)))=\frac{c_4'x_{41}}{c_{41}'}+\dfrac{(\frac{x_{11}x_{42}c_{41}'}{c_{42}'x_{22}}+\frac{x_{21}x_{42}c_{41}'}{x_{32}c_{42}'}+x_{31})^2c_{42}'}{c_{41}'x_{42}}\\&+\dfrac{(x_{21}+\frac{x_{11}x_{43}c_{42}'}{c_{43}'x_{33}}+\frac{x_{11}x_{32}}{x_{22}})^2c_{43}'}{x_{43}c_{42}'}+\dfrac{x_{11}^2cc_4'}{x_{44}c_{43}'}
\end{align*}
If we expand and factor these terms, the result is $\varepsilon_0(V_1(x))$.

To prove $e_0^ce_2^d=e_2^de_0^c$,
$e_0^ce_3^d=e_3^de_0^c$, and
   $e_1^ce_0^{c^2d}e_1^{cd}e_0^d=e_0^de_1^{cd}e_0^{c^2d}e_1^c$, we need the following identity: $\overline{\sigma}\circ e_i^c=\overline{e_{4-i}}^c\circ \overline{\sigma}$. It can be shown by direct calculation that this relation is true. As a result, $e_i^c\circ \overline{\sigma^{-1}}=\overline{\sigma^{-1}}\circ \overline{e_{4-i}}^c$.
Using this identity, we can prove the relations above.
\begin{itemize}
    \item  $e_0^ce_2^d=\overline{\sigma^{-1}}\overline{e_0}^c\overline{\sigma}e_2^d=\overline{\sigma^{-1}}\overline{e_0}^c\overline{e_2}^d\overline{\sigma}=\overline{\sigma^{-1}}\overline{e_2}^d\overline{e_0}^c\overline{\sigma}=e_2^d\overline{\sigma^{-1}}\overline{e_0}^c\overline{\sigma}=e_2^de_0^c$
    \item $e_0^ce_3^d=\overline{\sigma^{-1}}\overline{e_0}^c\overline{\sigma}e_3^d=\overline{\sigma^{-1}}\overline{e_0}^c\overline{e_1}^d\overline{\sigma}=\overline{\sigma^{-1}}\overline{e_1}^d\overline{e_0}^c\overline{\sigma}=e_3^d\overline{\sigma^{-1}}\overline{e_0}^c\overline{\sigma}=e_3^de_0^c$
    \item \begin{align*}
        &e_1^ce_0^{c^2d}e_1^{cd}e_0^d=e_1^c\overline{\sigma^{-1}}\overline{e_0}^{c^2d}\overline{\sigma}e_1^{cd}\overline{\sigma^{-1}}e\overline{e_0}^d\overline{\sigma}=\overline{\sigma^{-1}}\overline{e_3}^c\overline{e_0}^{c^2d}\overline{\sigma}\overline{\sigma^{-1}}\overline{e_3}^{cd}\overline{e_0}^d\overline{\sigma}\\
        &\overline{\sigma^{-1}}\overline{e_3}^c\overline{e_0}^{c^2d}\overline{e_3}^{cd}\overline{e_0}^d\overline{\sigma}=\overline{\sigma^{-1}}\overline{e_0}^d\overline{e_3}^{cd}\overline{e_0}^{c^2d}\overline{e_3}^c\overline{\sigma}=\overline{\sigma^{-1}}\overline{e_0}^d\overline{e_3}^{cd}\overline{\sigma}\overline{\sigma^{-1}}\overline{e_0}^{c^2d}\overline{e_3}^c\overline{\sigma}\\
        & \overline{\sigma^{-1}}\overline{e_0}^d\overline{\sigma}e_1^{cd}\overline{\sigma^{-1}}\overline{e_0}^{c^2d}\overline{\sigma}e_1^c=e_0^de_1^{cd}e_0^{c^2d}e_1^c
    \end{align*} 
    
\end{itemize}
The final relation, $e_0^ce_4^d=e_4^de_0^c$ can be shown via direct computation.
Thus, since all of the relations hold, $\mathcal{V}_1(x)$ is an affine geometric crystal.
\end{proof}

%%%%%%%%%%%%%%%%%
\section{Perfect Crystals and Ultradiscretization}
\subsection{Definitions and Background}
Let $\mathcal{B}$ be a classical crystal, a finite crystal associated with the Cartan datum $(A, \Pi, \Pi^\vee,\Bar{P}, \Bar{P}^\vee)$. It is also called a $U_q'(\mathfrak{g})$ crystal. Then $$\varepsilon(b)=\sum\limits_i \varepsilon_i(b)\Lambda_i\,\,\,\,\,\,\,\varphi(b)=\sum\limits \varphi_i(b)\Lambda_i$$
Therefore $\text{wt}(b)=\varphi(b)-\varepsilon(b)$. Also let $\Bar{P}^+_l=\{\lambda\in P|\langle c,\lambda\rangle=l\}$. With this notation, we can present the definition of a perfect crystal, as presented in \cite{HK}
\begin{definition}
For a positive integer $l>0$, $\mathcal{B}$ is a perfect crystal of level $l$ if it satisfies the following conditions:
\begin{enumerate}
    \item there exists a finite dimensional $U_q'(\mathfrak{g})$ module with a crystal basis whose crystal graph is isomorphic to $\mathcal{B}$
    \item $\mathcal{B}\otimes\mathcal{B}$ is connected.
    \item there exists a classical weight $\lambda_0\in\Bar{P}$ such that $$\text{wt}(\mathcal{B})\subset \lambda_0+\sum\limits_{i\neq 0} \mathbb{Z}_{\geq 0}\alpha_i$$ where $\#(\mathcal{B}_{\lambda_0})=1$
    \item for any $b\in\mathcal{B}$, we have $\langle c, \varepsilon(b)\rangle\geq l$
    \item for each $\lambda\in\Bar{P}_l^+$, there exist unique vectors $b^\lambda\in \mathcal{B}$, $b_\lambda\in \mathcal{B}$ such that $\varepsilon(b^\lambda)=\lambda$ and $\varphi(b_\lambda)=\lambda$
\end{enumerate}
\end{definition}

We also define $\mathcal{B}^\text{min} =\{b\in\mathcal{B}|\langle c, \varepsilon(b)\rangle=l\}$. Additionally, we can observe from the definitions that $\varepsilon,\varphi:\mathcal{B}^\text{min} \to \Bar{P}^+_l$ are both bijective. Now we define the limit of a family of perfect crystals:
\begin{definition}[\cite{KKM}]
A crystal $\mathcal{B}_\infty$ with element $b_\infty$ is called a limit of a family of perfect crystals $\{\mathcal{B}_l\}_{l\geq 1}$ if 
\begin{enumerate}
    \item $\text{wt }(b_\infty)=\varepsilon(b_\infty)=\varphi(b_\infty)=0$
    \item For any $(l,b)\in J$, there exists an embedding of crystals:
    $$f_{(l,b)}:T_{\varepsilon(b)}\otimes \mathcal{B}_l\otimes T_{-\varphi(b)}\hookrightarrow \mathcal{B}_\infty$$
    $$t_{\varepsilon(b)}\otimes b \otimes t_{-\varphi(b)}\mapsto b_\infty$$
    where $T_\lambda:=\{t_\lambda\}$ for $\lambda\in P$ and $\Tilde{e}_i(t_\lambda)=\Tilde{f}_i(t_\lambda)=0$ and $\varepsilon_i(t_\lambda)=\varphi_i(t_\lambda)=\infty$ and $\text{wt }(t_\lambda)=\lambda$
    \item $\mathcal{B}=\bigcup_{(l,b)\in J} \text{Im}_{f_{(l,b)}}$
\end{enumerate}
\end{definition}
If a limit exists, then $\{\mathcal{B}_l\}_{l\geq 1}$ is called a coherent family of perfect crystals. In the following sections we define the limit of perfect crystals corresponding to the $n=2,3,4$ cases. These limits are proven in []. Then we ultradiscretize the geometric crystals constructed in the previous sections and then we prove the conjecture for each of the three cases.
\subsection{n=2 Case}
The Langlands dual of $C_2^{(1)}$ is $D_3^{(2)}$. We can reparametrize the level $l$ $U_q(D_3^{(2)})$ perfect crystal $B(l\Lambda_2)$ that was presented in \cite{KMN2} and its limit $B^{2,\infty}$ as follows: 
$$B^{2,l}=\left\{(b_{ij})_{1\leq i\leq 2,i\leq j\leq i+2} \, \, \, \, \vline \begin{array}{l@{}}
(b_{ij})\in\mathbb{Z}^{\geq 0}, \sum_{j=i}^{i+2} b_{ij}=l, 1\leq i\leq 2\\
\sum_{j=1}^t b_{1j}\geq\sum_{j=2}^{t+1}b_{2j}, 1\leq t\leq 3\\
b_{11}=b_{22}+b_{23}, b_{24}=b_{12}+b_{13}
\end{array}\right\}$$

 $$B^{2,\infty}=\left\{(b_{ij})_{1\leq i\leq 2,i\leq j\leq i+2}\, \, \, \, \vline \begin{array}{l@{}}
(b_{ij})\in\mathbb{Z}, \sum_{j=i}^{i+2} b_{ij}=0, 1\leq i\leq 2\\
b_{11}=b_{22}+b_{23}, b_{24}=b_{12}+b_{13}
\end{array}\right\}$$
Associated with these crystals are the following actions and maps, $\Tilde{e}_i, \Tilde{f}_i, \varepsilon_i, \varphi_i$ for $i=0,1,2$. The conditions in parentheses correspond only to the conditions required for $B^{2,l}$, not $B^{2,\infty}$.
For $0\leq k\leq 2$, $\Tilde{e}_k(b)=(b_{ij}')$, where
\[ \begin{cases}
k=0 & b_{11}'=b_{11}-1, b_{24}'=b_{24}+1\, \begin{cases}
b_{23}'=b_{23}-1, b_{13}'=b_{13}+1& b_{23}>b_{12}, \\
b_{22}'=b_{22}-1, b_{12}'=b_{12}+1 & b_{23}\leq b_{12}, 
\end{cases}\\
k=1 & b_{11}'=b_{11}+1, b_{12}'=b_{12}-1, b_{23}'=b_{23}+1, b_{24}'=b_{24}-1 \,\,\,  \\
k=2 & \begin{cases}
b_{13}'=b_{13}-1, b_{12}'=b_{12}+1& b_{12}\leq b_{23}, \\
b_{23}'=b_{23}-1, b_{22}'=b_{22}+1 & b_{12}> b_{23}, 
\end{cases}
\end{cases}
\]

and $b_{ij}'=b_{ij}$ otherwise.

For $0\leq k \leq 2$, $\Tilde{f}_k(b)=(b_{ij}')$ where

\[ \begin{cases}
k=0 & b_{11}'=b_{11}+1, b_{24}'=b_{24}-1\, \text{and} \begin{cases}
b_{23}'=b_{23}+1, b_{13}'=b_{13}-1& b_{23}\geq b_{12}, \\
b_{22}'=b_{22}+1, b_{12}'=b_{12}-1 & b_{23}< b_{12}, 
\end{cases}\\
k=1 & b_{11}'=b_{11}-1, b_{12}'=b_{12}+1, b_{23}'=b_{23}-1, b_{24}'=b_{24}+1 \,\,\,  \\
k=2 & \begin{cases}
b_{13}'=b_{13}+1, b_{12}'=b_{12}-1& b_{12}> b_{23}, \\
b_{23}'=b_{23}+1, b_{22}'=b_{22}-1 & b_{12}\leq b_{23}, 
\end{cases}
\end{cases}\] 
and $b_{ij}=b_{ij}'$ otherwise.
We now give the maps $\varepsilon_i,\varphi_i$ for $i=0,1,2$.
\begin{align*}
  & \varepsilon_1(b)=b_{12},\,\,\varphi_1(b)=b_{23}\\
  & \varepsilon_2(b)=b_{13}+\text{max}(b_{23}-b_{12},0)\\
   &\varphi_2(b)=b_{22}+\text{max}(b_{12}-b_{23},0))\\
   &\varepsilon_0(b)=\begin{cases}
l-b_{24}-\text{min}(b_{12},b_{23}) & b\in B^{2,l}\\
-b_{24}-\text{min}(b_{12},b_{23}) & b\in B^{2,\infty}
\end{cases}\\
&\varphi_0(b)=\begin{cases}
l-b_{11}-\text{min}(b_{12},b_{23}) & b\in B^{2,l}\\
-b_{11}-\text{min}(b_{12},b_{23}) & b\in B^{2,\infty}
\end{cases}
\end{align*}

In [], we showed this is a limit of a coherent family of perfect crystals. Next we ultradiscretize the geometric crystal from the previous section.
Explicitly, for $C_2^{(1)}$, if $\mathcal{X}=\mathcal{UD(V)}$ then $\mathcal{X}=\mathbb{Z}^3$ as sets, and $\mathcal{UD(V)}$ is equipped with the following functions:
  \[ \text{wt}_i(x)=\begin{cases} 
      -2x_{11} & i=0 \\
      2x_{11}-x_{21}-x_{22} & i=1 \\
      2x_{21}+2x_{22}-2x_{11} & i=2
   \end{cases}\]
 \[\varepsilon_i(x)=\begin{cases}
\text{max}\{x_{21}, 2x_{11}-x_{22}\} & i=0\\
x_{22}-x_{11} & i=1\\
\text{max}\{-x_{22}, 2x_{11}-2x_{22}-x_{21}\} & i=2
\end{cases}\]
\[ e_i^c(x)=\begin{cases} 
    (C_2+x_{22}-c,x_{11}-c,x_{21}-C_2) & i=0\\
      (x_{22},c+x_{11},x_{21}) & i=1 \\
     (C_2+x_{22},x_{11},C+x_{21}-C_2) & i=2 
   \end{cases}
\]

where $C_2=\text{max}\{c+x_{21}+x_{22},2x_{11}\}-\text{max}\{x_{21}+x_{22},2x_{11}\}$.  

When we restrict $c$ to compute $\Tilde{e_k}$ and $\Tilde{f_k}$ for $k=0,1,2$, we get the following actions:
\begin{align*}
    &\Tilde{f_0}(x)=\begin{cases} 
      (x_{22}+1,x_{11}+1,x_{21}) & \text{if }2x_{11}-x_{22}\leq x_{21} \\
      (x_{22},x_{11}+1,x_{21}+1) & \text{if }2x_{11}-x_{22}> x_{21}\\
   \end{cases}\\
   &\Tilde{f_1}(x)=(x_{22},x_{11}-1,x_{21})\\
   & f_2(x)=\begin{cases} 
      (x_{22},x_{11},x_{21}-1) & \text{if }2x_{11}-x_{22}\geq x_{21} \\
      (x_{22}-1,x_{11},x_{21}) & \text{if }2x_{11}-x_{22}< x_{21}\\
   \end{cases}
\end{align*}
 \begin{align*}
     &\Tilde{e_0}(x)=\begin{cases} 
      (x_{22}-1,x_{11}-1,x_{21}) &\text{if } 2x_{11}-x_{22}> x_{21} \\
      (x_{22},x_{11}-1,x_{21}-1) &\text{if } 2x_{11}-x_{22}\leq x_{21}\\
   \end{cases}\\
   &\Tilde{e_1}(x)=(x_{22},x_{11}+1,x_{21})\\
   & \Tilde{e_2}(x)=\begin{cases} 
      (x_{22},x_{11},x_{21}+1) & \text{if }2x_{11}-x_{22}> x_{21} \\
      (x_{22}+1,x_{11},x_{21}) &\text{if } 2x_{11}-x_{22}\leq x_{21}\\
   \end{cases}
 \end{align*}
\begin{theorem}
 Let $\Omega:\mathcal{UD}(\mathcal{V})\to B^{2,\infty}$ be given by 
    \begin{align*}
      b_{11}=x_{11},\,\,\,  b_{12}=x_{22}-x_{11}\\
      b_{13}=-x_{22},\,\,\,b_{22}=x_{21}\\
      b_{23}=x_{11}-x_{21},\,\,\,b_{24}=-x_{11}
    \end{align*} and $\Omega^{-1}:B^{2,\infty}\to\mathcal{UD}(\mathcal{V})$ be given by $$x_{11}=b_{11},\,\, x_{21}=b_{22},\,\, x_{22}=b_{11}+b_{12}$$ 
    $\Omega$ is an isomorphism of crystals. 
\end{theorem}
\begin{proof}
Clearly we see that the map is bijective. We need to prove the following conditions to show that $\Omega$ is an isomorphism:
\begin{enumerate}
\item $\Omega(\Tilde{f_k}(b)=\Tilde{f_k}\Omega(b)$
\item $\Omega(\Tilde{e_k}(b)=\Tilde{e_k}\Omega(b)$
\item $\text{wt}_k(\Omega(b))=\text{wt}_k(b)$
\item $\varepsilon_k(\Omega(b)=\varepsilon_k(b)$
\end{enumerate}
Clearly if these hold, then $\varphi_k(\Omega(b))=\varepsilon_k(\Omega(b))+\text{wt}_k(\Omega(b))=\varepsilon_k(b)+\text{wt}_k(b)=\varphi_k(b)$, which will prove the isomorphism. Now we check the above conditions for $k=0,1,2$.
\begin{enumerate}
    \item For $k=0$ \[b_{11}'=b_{11}+1, b_{24}'=b_{24}-1\, \text{and} \begin{cases}
b_{23}'=b_{23}+1, b_{13}'=b_{13}-1& b_{23}\geq b_{12}, b_{23}>0\\
b_{22}'=b_{22}+1, b_{12}'=b_{12}-1 & b_{23}< b_{12}, b_{22}>0
\end{cases}
    \]
Applying the isomorphism we get:
    \[b_{11}'=x_{11}+1, b_{24}'=-x_{11}-1\, \text{and}\] \[\begin{cases}
b_{23}'=x_{11}-x_{21}+1, b_{13}'=-x_{22}-1& 2x_{11}-x_{22}\leq x_{21}\\
b_{22}'=x_{21}+1, b_{12}'=x_{22}-x_{11}-1 & 2x_{11}-x_{22}>x_{21}
\end{cases}
    \]
    This is equivalent to
    \[ \begin{cases} 
      (x_{22}+1,x_{11}+1,x_{21}) & 2x_{11}-x_{22}\leq x_{21} \\
      (x_{22},x_{11}+1,x_{21}+1) & 2x_{11}-x_{22}> x_{21}\\
   \end{cases}
\]
For $k=1$, $\Tilde{f_1}(b)=b'$ such that $b_{11}'=b_{11}-1, b_{12}'=b_{12}+1, b_{23}'=b_{23}-1, b_{24}'=b_{24}+1 \,\,\,  b_{23},b_{11}>0$.
Applying the isomorphism, we get $b_{11}'=x_{11}-1, b_{12}'=x_{22}-x_{11}+1, b_{23}'=x_{11}-x_{21}-1, b_{24}'=-x_{11}+1$. This is equivalent to $\Tilde{f_1}(x)=(x_{22},x_{11}-1,x_{21})$.

For $k=2$,  $\Tilde{f_2}(b)=b'$ such that
\[ \begin{cases}
b_{13}'=b_{23}+1, b_{12}'=b_{12}-1& b_{12}> b_{23}, b_{12}>0\\
b_{23}'=b_{23}+1, b_{22}'=b_{22}-1 & b_{12}\leq b_{23}, b_{22}>0
\end{cases}\] 

Applying the isomorphism, we get
\[ \begin{cases}
b_{13}'=-x_{22}+1, b_{12}'=x_{22}-x_{11}-1& 2x_{11}-x_{22}<x_{21}\\
b_{23}'=x_{11}-x_{21}+1, b_{22}'=x_{21}-1 & 2x_{11}-x_{22}\geq x_{21}
\end{cases}\] This is equivalent to the action we get from the ultra-discretized geometric crystal.
\item For $k=0$, applying isomorphism we get:
    \[b_{11}'=x_{11}-1, b_{24}'=-x_{11}+1\, \text{and}\] \[\begin{cases}
b_{23}'=x_{11}-x_{21}-1, b_{13}'=-x_{22}+1& 2x_{11}-x_{22}> x_{21}\\
b_{22}'=x_{21}-1, b_{12}'=x_{22}-x_{11}+1 & 2x_{11}-x_{22}\leq x_{21}
\end{cases}
    \]
    This is equivalent to
    \[ \begin{cases} 
      (x_{22}-1,x_{11}-1,x_{21}) & 2x_{11}-x_{22}> x_{21} \\
      (x_{22},x_{11}-1,x_{21}-1) & 2x_{11}-x_{22}\leq x_{21}\\
   \end{cases}
\]
For $k=1$, $\Tilde{e_1}(b)=b'$ such that $b_{11}'=b_{11}+1, b_{12}'=b_{12}-1, b_{23}'=b_{23}+1, and b_{24}'=b_{24}-1$ with $b_{23},b_{11}>0$.
Applying the isomorphism, we get $b_{11}'=x_{11}+1, b_{12}'=x_{22}-x_{11}-1, b_{23}'=x_{11}-x_{21}+1, b_{24}'=-x_{11}-1$. This is equivalent to $\Tilde{f_1}(x)=(x_{22},x_{11}+1,x_{21})$.

For $k=2$,  $\Tilde{e_2}(b)=b'$ such that
\[ \begin{cases}
b_{13}'=b_{23}-1, b_{12}'=b_{12}+1& b_{12}\geq b_{23}, b_{12}>0\\
b_{23}'=b_{23}-1, b_{22}'=b_{22}+1 & b_{12}< b_{23}, b_{22}>0
\end{cases}\] 

Applying the isomorphism, we get
\[ \begin{cases}
b_{13}'=-x_{22}-1, b_{12}'=x_{22}-x_{11}+1& 2x_{11}-x_{22}\geq x_{21}\\
b_{23}'=x_{11}-x_{21}-1, b_{22}'=x_{21}+1 & 2x_{11}-x_{22}< x_{21}
\end{cases}\] This is equivalent to the action we get from the ultra-discretized geometric crystal.
\item \begin{align*}
   &\Omega(\varepsilon_1(b))=\Omega(b_{12})=x_{22}-x_{11}\\
   &\varepsilon_2(b)=\begin{cases}
   b_{13} & b_{23}<b_{12}\\
   b_{13}+b_{23}-b_{12} & b_{12}\leq b_{23}
   \end{cases}=\Omega(\varepsilon_2(b))\\
   &\varepsilon_0(b)=\begin{cases} -b_{24}-b_{12} & b_{23}>b_{12}\\
   -b_{24}-b_{23} & b_{23}\leq b_{12}\end{cases}
   =\begin{cases}
2x_{11}-x_{22} &  2x_{11}-x_{22}> x_{21}\\
x_{21} & 2x_{11}-x_{22}\leq x_{21}
\end{cases}=\Omega(\varepsilon_0(b))
\end{align*}
\item Finally we check $\text{wt}_k(b)$.
\begin{align*}
    &\Omega(\text{wt}_0(b))=\Omega(b_{24}-b_{11})=-2x_{11}\\
    &\Omega(\text{wt}_1(b))=\Omega(b_{23}-b_{12})=2x_{11}-x_{22}-x_{21}\\
    &\Omega(\text{wt}_2(b))=\Omega(b_{12}+b_{22}-b_{13}-b_{23})=2x_{22}+2x_{21}-2x_{11}\\
\end{align*}
\end{enumerate}
Since all 4 conditions hold, $\Omega$ is an isomorphism.
\end{proof}
\subsection{n=3 Case}
The Langlands dual of $C_3^{(1)}$ is $D_4^{(2)}$. We give a parametrization for the perfect crystal associated with $V(l\Lambda_3)$ along with its associated limit below: 
$$B^{3,l}=\left\{(b_{ij})_{1\leq i\leq 3,i\leq j\leq i+3} \, \, \, \, \vline \begin{array}{l@{}}
(b_{ij})\in\mathbb{Z}^{\geq 0}, \sum_{j=i}^{i+3} b_{ij}=l, 1\leq i\leq 3\\
\sum_{j=i}^t b_{ij}\geq\sum_{j=i+1}^{t+1}b_{i+1,j}, 1\leq i\leq 2, 1\leq t\leq 4\\
b_{11}=b_{33}+b_{34}+b_{35}, b_{36}=b_{12}+b_{13}+b_{14}\\
b_{22}=b_{33}+b_{34},\,\,\,b_{25}=b_{13}+b_{14}
\end{array}\right\}$$
$$B^{3,\infty}=\left\{(b_{ij})_{1\leq i\leq 3,i\leq j\leq i+3} \, \, \, \, \vline \begin{array}{l@{}}
(b_{ij})\in\mathbb{Z}, \sum_{j=i}^{i+3} b_{ij}=0, 1\leq i\leq 3\\
b_{11}=b_{33}+b_{34}+b_{35}, b_{36}=b_{12}+b_{13}+b_{14}\\
b_{22}=b_{33}+b_{34},\,\,\,b_{25}=b_{13}+b_{14}
\end{array}\right\}$$

For $0\leq k\leq 3$, $b\in B^{3,\infty}$ $(B^{3,l})$, $\Tilde{e}_k(b)=(b_{ij}')$, where
\begin{itemize}
    \item $k=0$ for each case $b_{11}'=b_{11}-1, b_{36}'=b_{36}+1$ 
    \[ \begin{cases}
b_{12}'=b_{12}+1, b_{22}'=b_{22}-1, b_{23}'=b_{23}+1, b_{33}'=b_{33}-1 & b_{23}\geq b_{34}, 2b_{12}\geq b_{23}+b_{34},\\& b_{12}+b_{13}\geq b_{34}+b_{35}, \\&(b_{11},b_{22},b_{33}>0)\\
b_{22}'=b_{22}-1, b_{12}'=b_{12}+1, b_{24}'=b_{24}+1, b_{34}'=b_{34}-1 & b_{12}+b_{13}\geq b_{23}+ b_{35}, b_{34}>b_{23},\\& b_{12}+b_{24}\geq b_{23}+b_{35}, \\&(b_{11},b_{22},b_{34}>0)\\
b_{13}'=b_{13}+1, b_{23}'=b_{23}-1, b_{25}'=b_{25}+1, b_{35}'=b_{35}-1 & b_{13}\geq b_{24}, b_{35}+b_{34}> b_{12}+b_{24},\\& b_{12}<b_{23}+b_{35},\\& (b_{11},b_{23},b_{35}>0)\\
b_{14}'=b_{14}+1, b_{24}'=b_{24}-1, b_{25}'=b_{25}+1, b_{35}'=b_{35}-1 & b_{24}>b_{13}, 2b_{35}> b_{13}+b_{24},\\& b_{35}+b_{34}>b_{12}+b_{13},\\& (b_{11},b_{24},b_{35}> 0)
\end{cases}\]
\item $k=1$ 
$$ b_{12}'=b_{12}-1, b_{11}'=b_{11}+1, b_{36}'=b_{36}-1, b_{35}'=b_{35}+1,\,\,\, (b_{12},b_{36}>0)$$
\item $k=2$
\[\begin{cases} b_{12}'=b_{12}+1, b_{13}'=b_{13}-1, b_{24}'=b_{24}+1, b_{25}'=b_{25}-1 & b_{23}+b_{35}\leq b_{12}+b_{24}, \\&(b_{25},b_{13}>0)\\
b_{22}'=b_{22}+1, b_{23}'=b_{23}-1, b_{34}'=b_{34}+1, b_{35}'=b_{35}-1 & b_{23}+b_{35}>b_{12}+b_{24}, \\&(b_{23},b_{35}>0)
\end{cases}
\]
\item $k=3$
\[\begin{cases}b_{13}'=b_{13}+1, b_{14}'=b_{14}-1 & b_{24}\leq b_{13}, b_{13}+b_{23}\geq b_{24}+b_{34}, (b_{14}>0)\\
b_{23}'=b_{23}+1, b_{24}'=b_{24}-1 & b_{24}>b_{13}, b_{34}\leq b_{23}, (b_{24}>0)\\
b_{33}'=b_{33}+1, b_{34}'=b_{34}-1 & b_{34}>b_{23}, b_{24}+b_{34}>b_{13}+b_{23}, (b_{34}>0)
\end{cases}\]
\end{itemize}
and $b_{ij}'=b_{ij}$ otherwise.\\
For $0\leq k\leq 3$, $b\in B^{3,\infty}$ $(B^{3,l})$, $\Tilde{f}_k(b)=(b_{ij}')$, where
\begin{itemize}
    \item $k=0$ for each case $b_{11}'=b_{11}+1, b_{36}'=b_{36}-1$ 
    \[ \begin{cases}
b_{12}'=b_{12}-1, b_{22}'=b_{22}+1, b_{23}'=b_{23}-1, b_{33}'=b_{33}+1 & b_{23}> b_{34}, b_{12}+b_{13}> b_{34}+b_{35}, \\&2b_{12}>b_{23}+b_{34},\\&b_{12}+b_{24}\geq b_{23}+b_{35}, \\&(b_{36},b_{12},b_{23}>0)\\
b_{22}'=b_{22}+1, b_{12}'=b_{12}-1, b_{24}'=b_{24}-1, b_{34}'=b_{34}+1 & b_{12}+b_{13}>b_{23}+ b_{35}, b_{34}\geq b_{23}, \\ & b_{35}+b_{23}< b_{24}+b_{12}, \\&(b_{36},b_{12},b_{24}>0)\\
b_{13}'=b_{13}-1, b_{23}'=b_{23}+1, b_{25}'=b_{25}-1, b_{35}'=b_{35}+1 & b_{13}> b_{24}, b_{35}+b_{34}\geq b_{12}+b_{24}\\ & b_{12}\leq b_{23}+b_{35},\\& (b_{36},b_{13},b_{25}> 0)\\
b_{14}'=b_{14}-1, b_{24}'=b_{24}+1, b_{25}'=b_{25}-1, b_{35}'=b_{35}+1 & b_{24}\geq b_{13}, b_{35}+b_{34}\geq b_{12}+b_{13},\\& 2b_{35}\geq b_{24}+b_{13},\\& (b_{36},b_{14},b_{25}> 0)
\end{cases}\]
\item $k=1$ 
$$ b_{12}'=b_{12}+1, b_{11}'=b_{11}-1, b_{36}'=b_{36}+1, b_{35}'=b_{35}-1,\,\,\,(b_{11},b_{35}>0)$$
\item $k=2$
\[\begin{cases} b_{12}'=b_{12}-1, b_{13}'=b_{13}+1, b_{24}'=b_{24}-1, b_{25}'=b_{25}+1 & b_{23}+b_{35}< b_{12}+b_{24}, \\&(b_{24},b_{12}>0)\\
b_{22}'=b_{22}-1, b_{23}'=b_{23}+1, b_{34}'=b_{34}-1, b_{35}'=b_{35}+1 & b_{23}+b_{35}\geq b_{12}+b_{24}, \\&(b_{22},b_{34}>0)
\end{cases}
\]
\item $k=3$
\[\begin{cases}b_{13}'=b_{13}-1, b_{14}'=b_{14}+1 & b_{13}> b_{24}, b_{13}+b_{23}>b_{24}+b_{34}, (b_{14}>0)\\
b_{23}'=b_{23}-1, b_{24}'=b_{24}+1 & b_{23}>b_{34}, b_{13}\leq b_{24}, (b_{23}>0)\\
b_{33}'=b_{33}-1, b_{34}'=b_{34}+1 & b_{23}\leq b_{34}, b_{24}+b_{34}\geq b_{13}+b_{23}, (b_{34}>0)
\end{cases}\]
\end{itemize}
and $b_{ij}'=b_{ij}$ otherwise.

The actions $\varepsilon_i$, $\varphi_i$ and $\text{wt}_i$ are defined as follows:
$ \varepsilon_1(b)=b_{12}$\\
  $ \varepsilon_2(b)=b_{13}+(b_{23}-b_{12})_+$\\
  $\varepsilon_3(b)=b_{14}+\max\{b_{24}-b_{13},b_{34}+b_{24}-b_{13}-b_{23},0\}$\\
  $\varepsilon_0(b)=\begin{cases}
l-b_{36}-\text{min}\{b_{12}+b_{13},b_{23}+b_{24}, b_{34}+b_{35}, b_{12}+b_{24}, b_{23}+b_{35}\} & b\in B^{2,l}\\
-b_{36}- \text{min}\{b_{12}+b_{13},b_{23}+b_{24}, b_{34}+b_{35}, b_{12}+b_{24}, b_{23}+b_{35}\}& b\in B^{2,\infty}
\end{cases}$\\
$\varphi_0(b)=\begin{cases}
l-b_{11}-\text{min}\{b_{12}+b_{13},b_{23}+b_{24}, b_{34}+b_{35}, b_{12}+b_{24}, b_{23}+b_{35}\} & b\in B^{2,l}\\
-b_{11}-\text{min}\{b_{12}+b_{13},b_{23}+b_{24}, b_{34}+b_{35}, b_{12}+b_{24}, b_{23}+b_{35}\} & b\in B^{2,\infty}
\end{cases}$\\
$\varphi_1(b)=b_{35}$\\
$\varphi_2(b)=b_{22}-b_{33}+(b_{12}-b_{23})_+$\\
$\varphi_3(b)=b_{33}+\max\{b_{23}-b_{34}, b_{13}+b_{23}-b_{24}-b_{34},0\}$\\
$\text{wt}_0(b)=b_{36}-b_{11}$\\
$\text{wt}_1(b)=b_{35}-b_{12}$\\
$\text{wt}_2(b)=b_{22}-b_{33}-b_{25}+b_{14}+b_{12}-b_{23}$\\
$\text{wt}_3(b)=b_{33}-b_{14}+b_{13}-b_{24}+b_{23}-b_{34}$\\

In [], we showed this is a limit of a coherent family of perfect crystals. Now we apply the ultradiscretization functor on $\mathcal{V}_1$. For $C_3^{(1)}$, if $\mathcal{X}=\mathcal{UD(V)}$ then $\mathcal{X}=\mathbb{Z}^6$ as sets, and $\mathcal{UD(V)}$ is equipped with the following functions:
  \[ \text{wt}_i(x)=\begin{cases} 
      -2x_{11} & i=0 \\
      2x_{11}-x_{21}-x_{22} & i=1 \\
      2x_{21}+2x_{22}-x_{11}-x_{31}-x_{32}-x_{33} & i=2\\
      2x_{31}+2x_{32}+2x_{33}-2x_{21}-2x_{22} & i=3
   \end{cases}
\]
\[ \varepsilon_i(x)=\begin{cases}
\text{max}\{x_{31}, 2x_{21}-x_{32}, x_{21}+x_{11}-x_{22}, 2x_{11}+x_{32}-2x_{22}, 2x_{11}-x_{33}\} & i=0\\
x_{22}-x_{11} & i=1\\
\text{max}\{x_{33}-x_{22}, x_{33}+x_{32}+x_{11}-2x_{22}-x_{21}\} & i=2\\
\max\{-x_{33}, 2x_{22}-2x_{33}-x_{32}, 2x_{22}+2x_{21}-2x_{33}-2x_{32}-x_{31}\} & i=3
\end{cases}
\]
\begin{align*}
    &e_0^c(x)=(C_{24}'+x_{33}-c-C_2',x_{22}+\max\{x_{21}+x_{22}+C_{24}',x_{11}+x_{32}+C_{21}'\}-c-C_2'\\&-\max\{x_{21}+x_{22},x_{11}+x_{32}\}, C_{21}'+x_{32}-C_{24}',x_{11}-c, x_{21}-\max\{x_{21}+x_{22}+C_{24}',x_{11}+x_{32}\\&+C_{21}'\}+C_2'+\max\{x_{21}+x_{22},x_{11}+x_{32}\}, C_2'+x_{31}-C_{21}')
\end{align*} 
where $C_2'=\max\{x_{31},2\max \{x_{21}+x_{22},x_{11}+x_{32}\}-2x_{22}-x_{32},2x_{11}-x_{33}\}$,\\ $C_{21}'=\max\{c+x_{31},2\max \{x_{21}+x_{22},x_{11}+x_{32}\}-2x_{22}-x_{32},2x_{11}-x_{33}\}$ and \\$C_{24}'=\max\{c+x_{31},c+2\max\{x_{21}+x_{22},x_{11}+x_{32}\}-2x_{22}-x_{32},2x_{11}-x_{33}\}$.
$$e_1^c(x)=(x_{33},x_{22},x_{32},x_{11}+c, x_{21},x_{31})$$
$$e_2^c(x)=(x_{33},C_2+x_{22},x_{32},x_{11},c+x_{21}-C_2,x_{31})$$
where $C_2=\max\{c+x_{21}+x_{22},x_{11}+x_{32}\}-\max\{x_{21}+x_{22},x_{11}+x_{32}\}$.
$$e_3^c(x)=(C_{31}+x_{33}-C_3, x_{22}, C_{32}+x_{32}-C_{31}, x_{11},x_{21}, c+C_3-C_32+x_{31})$$
where $C_3=\max\{2x_{32}+x_{31}+x_{33}, 2x_{22}+x_{31}+x_{32}, 2x_{21}+2x_{22}\}$, $C_{31}=\max\{c+2x_{32}+x_{31}+x_{33}, 2x_{22}+x_{31}+x_{32}, 2x_{21}+2x_{22}\}$, and $C_{32}=\max\{c+2x_{32}+x_{31}+x_{33}, c+2x_{22}+x_{31}+x_{32}, 2x_{21}+2x_{22}\}$.

We show the simplification of $e_0^c(x)$ when we set $c=1$. Plugging in, we get 
\begin{align*}
    &(C_{24}'+x_{33}-1-C_2', x_{22}+\max\{x_{11}+x_{32}+C_{21}',x_{21}+x_{22}+C_{24}'\}-1-C_2'\\&-\max\{x_{21}+x_{22},x_{11}+x_{32}\},C_{21}'+x_{32}-C_{24}',x_{11}-1, x_{21}-\max\{x_{21}+x_{22}+C_{24}',x_{11}+x_{32}+C_{21}'\}\\&+C_2'+\max\{x_{21}+x_{22},x_{11}+x_{32}\}, C_2'+x_{31}-C_{21}')
\end{align*}
along with $C_2'=\max\{x_{31},2x_{21}-x_{32},2x_{11}-2x_{22}+x_{32},2x_{11}-x_{33}\}$, $C_{21}'=\max\{x_{31}+1,2x_{21}-x_{32},2x_{11}-2x_{22}+x_{32},2x_{11}-x_{33}\}$ and $C_{24}'=\max\{x_{31}+1,2x_{21}-x_{32}+1,2x_{11}-2x_{22}+x_{32}+1, 2x_{11}-x_{33}\}$. Then based on which of the 4 components is maximal, we get 4 cases for $\Tilde{e_0}$.
 \[ \Tilde{e_0}(x)=\begin{cases} 
      (x_{33},x_{22},x_{32}, x_{11}-1,x_{21}-1, x_{31}-1) & x_{31}>2x_{21}-x_{32},x_{31}>2x_{11}-2x_{22}+x_{32},\\&x_{31}>2x_{11}-x_{33} \\
    (x_{33},x_{22},x_{32}-1, x_{11}-1,x_{21}-1, x_{31}) & 2x_{21}-x_{32}\geq x_{31}, \\&2x_{21}-x_{32}>2x_{11}-2x_{22}+x_{32},\\&2x_{21}-x_{32}>2x_{11}-x_{33}\\
    (x_{33},x_{22}-1,x_{32}-1, x_{11}-1,x_{21}, x_{31}) & 2x_{11}-2x_{22}+x_{32}\geq x_{31},\\&2x_{11}-2x_{22}+x_{32}\geq 2x_{21}-x_{32},\\&2x_{11}-2x_{22}+x_{32}>2x_{11}-x_{33} \\
    (x_{33}-1,x_{22}-1,x_{32}, x_{11}-1,x_{21}, x_{31}) & 2x_{11}-x_{33}\geq x_{31},2x_{11}-x_{33}\geq 2x_{21}-x_{32},\\&2x_{11}-x_{33}\geq 2x_{11}-2x_{22}+x_{32} \\
   \end{cases}
\]
When we restrict $c$ to compute $\Tilde{e_k}$ and $\Tilde{f_k}$ for $k=0,1,2,3$, we get the following actions (besides $\Tilde{e_0}$ which is given above):

 \[ \Tilde{e_1}(x)=(x_{33},x_{22},x_{32},x_{11}+1,x_{21},x_{31})
\]
 \[ \Tilde{e_2}(x)=\begin{cases} 
      (x_{33},x_{22}+1,x_{32},x_{11},x_{21},x_{31}) & x_{21}+x_{22}\geq x_{11}+x_{32} \\
     (x_{33},x_{22},x_{32},x_{11},x_{21}+1,x_{31}) & x_{21}+x_{22}< x_{11}+x_{32}
   \end{cases}
\]
\[ \Tilde{e_3}(x)=\begin{cases}
(x_{33}+1,x_{22},x_{32},x_{11},x_{21},x_{31}) & 2x_{32}+x_{31}+x_{33}\geq 2x_{22}+x_{31}+x_{32},\\& 2x_{32}+x_{31}+x_{33}\geq 2x_{22}+2x_{21} \\
(x_{33},x_{22},x_{32}+1,x_{11},x_{21},x_{31}) & 2x_{22}+x_{31}+x_{32}>2x_{32}+x_{31}+x_{33},\\& 2x_{22}+x_{31}+x_{32}\geq 2x_{22}+2x_{21} \\
(x_{33},x_{22},x_{32},x_{11},x_{21},x_{31}+1) & 2x_{22}+2x_{21}>2x_{32}+x_{31}+x_{33},\\&2x_{22}+2x_{21}> 2x_{22}+x_{31}+x_{32}\\
\end{cases}\]
 \[ \Tilde{f_0}(x)=\begin{cases} 
      (x_{33},x_{22},x_{32}, x_{11}+1,x_{21}+1, x_{31}+1) & x_{31}\geq 2x_{21}-x_{32},x_{31}\geq 2x_{11}-2x_{22}+x_{32},\\& x_{31}\geq 2x_{11}-x_{33} \\
    (x_{33},x_{22},x_{32}+1, x_{11}+1,x_{21}+1, x_{31}) & 2x_{21}-x_{32}>x_{31},2x_{21}-x_{32}\geq 2x_{11}-2x_{22}+x_{32},\\&2x_{21}-x_{32}\geq 2x_{11}-x_{33}\\
    (x_{33},x_{22}+1,x_{32}+1, x_{11}+1,x_{21}, x_{31}) & 2x_{11}-2x_{22}+x_{32}> x_{31},\\&2x_{11}-2x_{22}+x_{32}>2x_{21}-x_{32},\\&2x_{11}-2x_{22}+x_{32}\geq 2x_{11}-x_{33} \\
    (x_{33}+1,x_{22}+1,x_{32}, x_{11}+1,x_{21}, x_{31}) & 2x_{11}-x_{33}>x_{31},2x_{11}-x_{33}>2x_{21}-x_{32},\\& 2x_{11}-x_{33}>2x_{11}-2x_{22}+x_{32},\\
   \end{cases}
\]
 \[ \Tilde{f_1}(x)=(x_{33},x_{22},x_{32},x_{11}-1,x_{21},x_{31})
\]
 \[ \Tilde{f_2}(x)=\begin{cases} 
      (x_{33},x_{22},x_{32},x_{11},x_{21}-1,x_{31}) & x_{21}+x_{22}\geq x_{11}+x_{32} \\
     (x_{33},x_{22}-1,x_{32},x_{11},x_{21},x_{31}) & x_{21}+x_{22}< x_{11}+x_{32}
   \end{cases}
\]
\[ \Tilde{f_3}(x)=\begin{cases}
(x_{33}-1,x_{22},x_{32},x_{11},x_{21},x_{31}) & 2x_{32}+x_{31}+x_{33}>2x_{22}+x_{31}+x_{32},\\&2x_{32}+x_{31}+x_{33}>2x_{22}+2x_{21} \\
(x_{33},x_{22},x_{32}-1,x_{11},x_{21},x_{31}) & 2x_{22}+x_{31}+x_{32}\geq 2x_{32}+x_{31}+x_{33},\\& 2x_{22}+x_{31}+x_{32}>2x_{22}+2x_{21} \\
(x_{33},x_{22},x_{32},x_{11},x_{21},x_{31}-1) & 2x_{22}+x_{21}\geq 2x_{32}+x_{31}+x_{33},\\& 2x_{22}+x_{21}\geq 2x_{22}+x_{31}+x_{32} \\
\end{cases}\]
\begin{theorem}
 Let $\Omega:\mathcal{UD}(\mathcal{V})\to B^{3,\infty}$ be given by 
    \begin{align*}
      &b_{11}=x_{11},\,\,\,  b_{12}=x_{22}-x_{11},\,\,\,
      b_{13}=-x_{22}+x_{33},\,\,\,b_{14}=-x_{33}\\
      &b_{22}=x_{21},\,\,\,b_{23}=-x_{21}+x_{32},\,\,\,
      b_{24}=x_{22}-x_{32},\,\,\,b_{25}=-x_{22}\\
      &b_{33}=x_{31},\,\,\,b_{34}=x_{21}-x_{31},\,\,\,
      b_{35}=x_{11}-x_{21},\,\,\,b_{36}=-x_{11}
    \end{align*} and $\Omega^{-1}:B^{3,\infty}\to\mathcal{UD}(\mathcal{V})$ be given by $$x_{11}=b_{11},\,\, x_{21}=b_{22},\,\, x_{22}=b_{11}+b_{12},\,\,\, x_{31}=b_{33},\,\,\,x_{32}=b_{22}+b_{23},\,\,\,x_{33}=-b_{14}$$ 
    $\Omega$ is an isomorphism of crystals. 
\end{theorem}
\begin{proof}
Clearly we see that the map is bijective. We need to prove the following conditions to show that $\Omega$ is an isomorphism:
\begin{enumerate}
\item $\Omega(\Tilde{f_k}(b)=\Tilde{f_k}\Omega(b)$
\item $\Omega(\Tilde{e_k}(b)=\Tilde{e_k}\Omega(b)$
\item $\text{wt}_k(\Omega(b))=\text{wt}_k(b)$
\item $\varepsilon_k(\Omega(b)=\varepsilon_k(b)$
\end{enumerate}
Clearly if these hold, then $\varphi_k(\Omega(b))=\varepsilon_k(\Omega(b))+\text{wt}_k(\Omega(b))=\varepsilon_k(b)+\text{wt}_k(b)=\varphi_k(b)$. Now we check the above conditions for each $k=0,1,2,3$.
First we check 1. There are 4 cases for $\Tilde{f_0}(b)$.
\begin{enumerate}
    \item $b_{11}'=b_{11}+1,\,\,\, b_{36}'=b_{36}-1,\,\,\,b_{12}'=b_{12}-1,\,\,\,b_{22}'=b_{22}+1,\,\,\,b_{23}'=b_{23}-1,\,\,\,b_{33}'=b_{33}+1$ if $b_{23}>b_{34},\,\,\,b_{12}+b_{13}>b_{34}+b_{35},\,\,\,2b_{12}>b_{23}+b_{24}$. Applying the isomorphism, we get that all $x_{ij}$ stay the same except $x_{11}'=x_{11}+1$, $x_{21}'=x_{21}+1$ and $x_{31}'=x_{31}+1$. The inequalities give us the following:
    $b_{23}>b_{34}\Longleftrightarrow x_{31}>2x_{21}-x_{32}$, $b_{12}+b_{13}>b_{34}+b_{35}\Longleftrightarrow x_{31}>2x_{11}-x_{33}$ and $2b_{12}> b_{23}+b_{34}\Longleftrightarrow x_{31}> 2x_{11}-2x_{22}+x_{32}$. 
    \item $b_{11}'=b_{11}+1,\,\,\, b_{36}'=b_{36}-1,\,\,\,b_{12}'=b_{12}-1,\,\,\,b_{22}'=b_{22}+1,\,\,\,b_{24}'=b_{24}-1,\,\,\,b_{34}'=b_{34}+1$ if $b_{23}\leq b_{34},\,\,\,b_{12}+b_{13}> b_{23}+b_{35},\,\,\,b_{12}+b_{24}>b_{23}+b_{35}$.
     Applying the isomorphism, we get that all $x_{ij}$ stay the same except $x_{11}'=x_{11}+1$, $x_{21}'=x_{21}+1$ and $x_{32}'=x_{32}+1$. The inequalities give us the following:
    $b_{23}\leq b_{34}\Longleftrightarrow x_{31}
    \leq 2x_{21}-x_{32}$, $b_{12}+b_{13}> b_{23}+b_{35}\Longleftrightarrow 2x_{21}-x_{32}> 2x_{11}-x_{33}$ and $b_{12}+b_{24}> b_{23}+b_{35}\Longleftrightarrow 2x_{21}-x_{32}> 2x_{11}-2x_{22}+x_{32}$.
    \item $b_{11}'=b_{11}+1,\,\,\, b_{36}'=b_{36}-1,\,\,\,b_{13}'=b_{13}-1,\,\,\,b_{23}'=b_{23}+1,\,\,\,b_{25}'=b_{25}-1,\,\,\,b_{35}'=b_{35}+1$ if $b_{13}> b_{24},\,\,\,b_{12}+b_{24}\leq b_{34}+b_{35},\,\,\,b_{12}+b_{24}\leq b_{23}+b_{35}$.
     Applying the isomorphism, we get that all $x_{ij}$ stay the same except $x_{11}'=x_{11}+1$, $x_{22}'=x_{22}+1$ and $x_{32}'=x_{32}+1$. The inequalities give us the following:
    $b_{13}> b_{24}\Longleftrightarrow 2x_{11}-2x_{22}+x_{32}>2x_{11}-x_{33} $, $b_{12}+b_{24}\leq b_{34}+b_{35}\Longleftrightarrow 2x_{11}-2x_{22}+x_{32}\geq x_{31}$ and $b_{12}+b_{24}\leq b_{23}+b_{35}\Longleftrightarrow 2x_{21}-x_{32}\leq 2x_{11}-2x_{22}+x_{32}$.
    \item $b_{11}'=b_{11}+1,\,\,\, b_{36}'=b_{36}-1,\,\,\,b_{14}'=b_{14}-1,\,\,\,b_{24}'=b_{24}+1,\,\,\,b_{25}'=b_{25}-1,\,\,\,b_{35}'=b_{35}+1$ if $b_{13}\leq b_{24},\,\,\,b_{12}+b_{13}\leq b_{34}+b_{35},\,\,\,2b_{35}\geq b_{13}+b_{24}$.
     Applying the isomorphism, we get that all $x_{ij}$ stay the same except $x_{11}'=x_{11}+1$, $x_{22}'=x_{22}+1$ and $x_{33}'=x_{33}+1$. The inequalities give us the following:
    $b_{13}\leq b_{24}\Longleftrightarrow 2x_{11}-2x_{22}+x_{32}\leq 2x_{11}-x_{33} $, $b_{12}+b_{13}\leq b_{34}+b_{35}\Longleftrightarrow 2x_{11}-x_{33}\geq x_{31}$ and $2b_{35}\geq b_{13}+b_{24}\Longleftrightarrow 2x_{11}-x_{33}\geq 2x_{21}-x_{32}$.
\end{enumerate}
Next we look at $\Tilde{f_1}(b)$. We have $b_{11}'=b_{11}-1,\,\,\,b_{12}'=b_{12}+1,\,\,\,b_{35}'=b_{35}-1,\,\,\,b_{36}'=b_{36}+1$. This gives $x_{11}'=x_{11}-1$. Every other $x_{ij}$ stays the same.
Next we consider the two cases of $\Tilde{f_2}(b)$. 
\begin{enumerate}
    \item $b_{12}'=b_{12}-1,\,\,\, b_{13}'=b_{13}+1,\,\,\,b_{24}'=b_{24}-1,\,\,\,b_{25}'=b_{25}+1$ if $b_{12}+b_{24}\geq b_{23}+b_{35}$. Applying the isomorphism, we have $x_{22}-1$ all others remaining the same if $x_{21}+x_{22}\geq x_{11}+x_{32}$.
    \item $b_{22}'=b_{22}-1,\,\,\, b_{23}'=b_{23}+1,\,\,\,b_{34}'=b_{34}-1,\,\,\,b_{35}'=b_{35}+1$ if $b_{12}+b_{24}< b_{23}+b_{35}$. Applying the isomorphism, we have $x_{22}-1$ all others remaining the same if $x_{21}+x_{22}<x_{11}+x_{32}$.
\end{enumerate}
Finally we consider the three cases of $\Tilde{f_3}(b)$.
\begin{enumerate}
    \item $b_{13}'=b_{13}-1,\,\,\,b_{14}'=b_{14}+1$ if $b_{13}+b_{23}>b_{24}+b_{34}$ and $b_{13}>b_{24}$. Applying the isomorphism, we get $x_{33}'=x_{33}-1$ and everything else stays the same. The inequalities give us:$b_{13}>b_{24}\Longleftrightarrow x_{32}+x_{33}>2x_{22}$ and $b_{13}+b_{23}>b_{24}+b_{34}\Longleftrightarrow x_{31}+2x_{32}+x_{33}>2x_{21}+2x_{22}$.
    \item $b_{23}'=b_{23}-1,\,\,\,b_{24}'=b_{24}+1$ if $b_{23}>b_{34}$ and $b_{13}\leq b_{24}$. Applying the isomorphism, we get $x_{32}'=x_{32}-1$ and everything else stays the same. The inequalities give us:$b_{13}\leq b_{24}\Longleftrightarrow x_{32}+x_{33} \leq 2x_{22}$ and $b_{23}>b_{34}\Longleftrightarrow 2x_{21}<x_{31}+x_{32}$.
   \item $b_{33}'=b_{33}-1,\,\,\,b_{34}'=b_{34}+1$ if $b_{13}+b_{23}\leq b_{24}+b_{34}$ and $b_{13}\leq b_{24}$. Applying the isomorphism, we get $x_{33}'=x_{33}-1$ and everything else stays the same. The inequalities give us:$b_{23}\leq b_{34}\Longleftrightarrow x_{32}+x_{33}\leq 2x_{22}$ and $b_{13}+b_{23}\leq b_{24}+b_{34}\Longleftrightarrow x_{31}+2x_{32}+x_{33}\leq 2x_{21}+2x_{22}$.
\end{enumerate}
This proves that $\Omega(\Tilde{f_k}(b)=\Tilde{f_k}\Omega(b)$.

Now we consider 2.
There are 4 cases for $\Tilde{e_0}(b)$.
\begin{enumerate}
    \item $b_{11}'=b_{11}-1,\,\,\, b_{36}'=b_{36}+1,\,\,\,b_{12}'=b_{12}+1,\,\,\,b_{22}'=b_{22}-1,\,\,\,b_{23}'=b_{23}+1,\,\,\,b_{33}'=b_{33}-1$ if $b_{23}\geq b_{34},\,\,\,b_{12}+b_{13}\geq b_{34}+b_{35},\,\,\,2b_{12}\geq b_{23}+b_{34}$. Applying the isomorphism, we get that all $x_{ij}$ stay the same except $x_{11}'=x_{11}-1$, $x_{21}'=x_{21}-1$ and $x_{31}'=x_{31}-1$. The inequalities give us the following:
    $b_{23}>b_{34}\Longleftrightarrow x_{31}>2x_{21}-x_{32}$, $b_{12}+b_{13}>b_{34}+b_{35}\Longleftrightarrow x_{31}>2x_{11}-x_{33}$ and 2$b_{12}\geq b_{23}+b_{34}\Longleftrightarrow x_{31}\geq 2x_{11}-2x_{22}+x_{32}$. The final inequality repeats these conditions.
    \item $b_{11}'=b_{11}-1,\,\,\, b_{36}'=b_{36}+1,\,\,\,b_{12}'=b_{12}+1,\,\,\,b_{22}'=b_{22}-1,\,\,\,b_{24}'=b_{24}+1,\,\,\,b_{34}'=b_{34}-1$ if $b_{34}>b_{23},\,\,\,b_{12}+b_{24}\geq b_{23}+b_{35},\,\,\, b_{12}+b_{13}\geq b_{23}+b_{35}$.
     Applying the isomorphism, we get that all $x_{ij}$ stay the same except $x_{11}'=x_{11}-1$, $x_{21}'=x_{21}-1$ and $x_{32}'=x_{32}-1$. The inequalities give us the following:
    $b_{12}+b_{24}\geq b_{23}+b_{35}\Longleftrightarrow 2x_{21}-x_{32}
    \geq 2x_{11}-2x_{22}+x_{32}$, $b_{12}+b_{13}\geq b_{23}+b_{35}\Longleftrightarrow 2x_{21}-x_{32}\geq 2x_{11}-x_{33}$ and $b_{34}> b_{23}\Longleftrightarrow 2x_{21}-x_{32}> x_{31}$.
    \item $b_{11}'=b_{11}-1,\,\,\, b_{36}'=b_{36}+1,\,\,\,b_{13}'=b_{13}+1,\,\,\,b_{23}'=b_{23}-1,\,\,\,b_{25}'=b_{25}+1,\,\,\,b_{35}'=b_{35}-1$ if $b_{13}\geq b_{24},\,\,\,b_{12}+b_{24}<b_{34}+b_{35},\,\,\,b_{12}+b_{24}< b_{23}+b_{35}$.
     Applying the isomorphism, we get that all $x_{ij}$ stay the same except $x_{11}'=x_{11}-1$, $x_{22}'=x_{22}-1$ and $x_{32}'=x_{32}-1$. The inequalities give us the following:
    $b_{13}\geq b_{24}\Longleftrightarrow 2x_{11}-2x_{22}+x_{32}\geq 2x_{11}-x_{33} $, $b_{12}+b_{24}<b_{34}+b_{35}\Longleftrightarrow 2x_{11}-2x_{22}+x_{32}> 2x_{21}-x_{32}$ and $b_{12}+b_{24}> b_{34}+b_{35}\Longleftrightarrow x_{31}< 2x_{11}-2x_{22}+x_{32}$.
    \item $b_{11}'=b_{11}-1,\,\,\, b_{36}'=b_{36}+1,\,\,\,b_{14}'=b_{14}+1,\,\,\,b_{24}'=b_{24}-1,\,\,\,b_{25}'=b_{25}+1,\,\,\,b_{35}'=b_{35}-1$ if $b_{13}< b_{24},\,\,\,b_{12}+b_{13}< b_{34}+b_{35},\,\,\,b_{13}+b_{24}<2b_{35}$.
     Applying the isomorphism, we get that all $x_{ij}$ stay the same except $x_{11}'=x_{11}-1$, $x_{22}'=x_{22}-1$ and $x_{33}'=x_{33}-1$. The inequalities give us the following:
    $b_{13}< b_{24}\Longleftrightarrow 2x_{11}-2x_{22}+x_{32}< 2x_{11}-x_{33} $, $b_{12}+b_{13}< b_{34}+b_{35}\Longleftrightarrow 2x_{11}-x_{33}> x_{31}$ and $b_{24}\leq b_{35}\Longleftrightarrow 2x_{11}-2x_{22}+x_{32}< 2x_{11}-x_{33}$.
\end{enumerate}
Next we look at $\Tilde{e_1}(b)$. We have $b_{11}'=b_{11}+1,\,\,\,b_{12}'=b_{12}-1,\,\,\,b_{35}'=b_{35}+1,\,\,\,b_{36}'=b_{36}-1$. This gives $x_{11}'=x_{11}+1$. Every other $x_{ij}$ stays the same.
Next we consider the two cases of $\Tilde{e_2}(b)$. 
\begin{enumerate}
    \item $b_{12}'=b_{12}+1,\,\,\, b_{13}'=b_{13}-1,\,\,\,b_{24}'=b_{24}+1,\,\,\,b_{25}'=b_{25}-1$ if $b_{12}+b_{24}\geq b_{23}+b_{35}$. Applying the isomorphism, we have $x_{22}+1$ all others remaining the same if $x_{21}+x_{22}\geq x_{11}+x_{32}$.
    \item $b_{22}'=b_{22}+1,\,\,\, b_{23}'=b_{23}-1,\,\,\,b_{34}'=b_{34}+1,\,\,\,b_{35}'=b_{35}-1$ if $b_{12}+b_{24}< b_{23}+b_{35}$. Applying the isomorphism, we have $x_{22}-1$ all others remaining the same if $x_{21}+x_{22}<x_{11}+x_{32}$.
\end{enumerate}
Finally we consider the three cases of $\Tilde{e_3}(b)$.
\begin{enumerate}
    \item $b_{13}'=b_{13}+1,\,\,\,b_{14}'=b_{14}-1$ if $b_{23}\geq b_{34}$ and $b_{13}\geq b_{24}$. Applying the isomorphism, we get $x_{33}'=x_{33}+1$ and everything else stays the same. The inequalities give us:$b_{13}\geq b_{24}\Longleftrightarrow x_{32}+x_{33}\geq 2x_{22}$ and $b_{13}+b_{23}\geq b_{24}+b_{34}\Longleftrightarrow x_{31}+x_{32}\geq 2x_{21}$.
    \item $b_{23}'=b_{23}+1,\,\,\,b_{24}'=b_{24}-1$ if $b_{23}\geq b_{34}$ and $b_{13}< b_{24}$. Applying the isomorphism, we get $x_{32}'=x_{32}+1$ and everything else stays the same. The inequalities give us:$b_{13}< b_{24}\Longleftrightarrow x_{33}+x_{32}\geq 2x_{22}$ and $b_{23}\geq b_{34}\Longleftrightarrow 2x_{21}\leq x_{31}+x_{32}$.
   \item $b_{33}'=b_{33}+1,\,\,\,b_{34}'=b_{34}-1$ if $b_{13}+b_{23}< b_{24}+b_{34}$ and $b_{24}> b_{13}$. Applying the isomorphism, we get $x_{33}'=x_{33}+1$ and everything else stays the same. The inequalities give us:$b_{13}+b_{23}<b_{24}+b_{34}\Longleftrightarrow x_{31}+2x_{32}+x_{33}\leq 2x_{22}+2x_{21}$ and $b_{23}< b_{34}\Longleftrightarrow x_{31}+x_{32}< 2x_{21}$.
\end{enumerate}
Therefore, $\Omega(\Tilde{e_k}(b)=\Tilde{e_k}\Omega(b)$

Now we show 3.
\begin{align*}
 &\text{wt}_0(b)=b_{36}-b_{11}=-2x_{11}=\text{wt}_0(\Omega(b))\\
 &\text{wt}_1(b)=b_{35}-b_{12}=2x_{11}-x_{21}-x_{22}=\text{wt}_1(\Omega(b))\\
 &\text{wt}_2(b)=b_{22}+b_{14}+b_{12}-b_{33}-b_{25}-b_{23}=2x_{21}+2x_{22}-x_{11}-x_{31}-x_{32}-x_{33}=\text{wt}_2(\Omega(b))\\
 &\text{wt}_3(b)=b_{33}-b_{14}+b_{13}-b_{24}+b_{23}-b_{34}=2x_{31}+2x_{32}+2x_{33}-2x_{21}-2x_{22}=\text{wt}_3(\Omega(b))
\end{align*}
This proves $\text{wt}_k(\Omega(b))=\text{wt}_k(b)$. Finally we show 4.
\[\varepsilon_0(b)=\begin{cases}
-b_{36}-b_{12}-b_{13} & A<B,C,D,E\\
-b_{36}-b_{23}-b_{24} & B\geq A,\,\, B<C,D,E\\
-b_{36}-b_{34}-b_{35} & C\geq A,B\,\,C>D,E\\
-b_{36}-b_{12}-b_{24} & D\geq A,B,C,\,\,D>E\\
-b_{36}-b_{23}-b_{35} & E\geq A,B,C,D
\end{cases}\]
where $A=b_{12}+b_{13}$, $B=b_{23}+b_{24}$, $C=b_{34}+b_{35}$, $D=b_{12}+b_{24}$, and $E=b_{23}+b_{35}$.

Applying the isomorphism, we see this is equivalent to
\[\varepsilon_0(b)=\begin{cases}
2x_{11}-x_{33} & (1)>(2),(3),(4),(5)\\
x_{11}+x_{21}-x_{22} & (2)\leq (1),\,\,(2)>(3),(4),(5)\\
x_{31} & (3)\leq (1),(2),\,\,(3)>(4),(5)\\
2x_{11}-2x_{22}+x_{32} & (4)\leq (1),(2),(3),\,\,(4)>(5)\\
2x_{21}-x_{32} & (5)\leq (1),(2),(3),(4)
\end{cases}\]
Where $(1)=2x_{11}-x_{33}$, $(2)=x_{11}+x_{21}-x_{22}$, $(3)=x_{31}$, $(4)=2x_{11}-2x_{22}+x_{32}$, $(5)=2x_{21}-x_{32}$
$$\varepsilon_1(b)=b_{12}=x_{22}-x_{11}=\varepsilon_1(\Omega(b))$$
\[\varepsilon_2(b)=\begin{cases}
b_{25}-b_{14} &b_{23}>b_{12}\\
b_{25}-b_{14}+b_{23}-b_{12} & b_{23}\leq b_{12} \end{cases}\]
Applying the isomorphism, we see this is equivalent to
\[\varepsilon_2(b)=\begin{cases}
x_{33}-x_{22} & x_{21}+x_{22}>x_{11}+x_{32}\\  x_{33}+x_{32}+x_{11}-2x_{22}-x_{21} & x_{21}+x_{22}\leq x_{11}+x_{32}
\end{cases}\]
\[\varepsilon_3(b)=\begin{cases}
b_{14} & b_{24}\leq b_{13},\,\,\, b_{34}\leq b_{23}\\
b_{14}+b_{24}-b_{13} & b_{24}>b_{13},\,\,\, b_{34}\leq b_{23}\\
b_{14}+b_{24}-b_{13}+b_{34}-b_{23} & b_{34}>b_{23},\,\,\,b_{24}>b_{13}
\end{cases}\]
Applying the isomorphism, we see this equivalent to
\[\varepsilon_3(b)=\begin{cases}
-x_{33} & x_{32}+x_{33}\geq 2x_{22},\,\,\, x_{31}+x_{32}\geq 2x_{21}\\
2x_{22}-2x_{33}-x_{32} & x_{32}+x_{33}< 2x_{22},\,\,\, x_{31}+x_{32}\geq 2x_{21}\\
-2x_{33}+2x_{21}+2x_{22}-2x_{32}-x_{31} & x_{32}+x_{33}< 2x_{22},\,\,\, x_{31}+x_{32}< 2x_{21}\end{cases}\]

Therefore all of the conditions are satisfied, so $\Omega$ is an isomorphism.
\end{proof}
\subsection{n=4 Case}
The Langlands Dual of $C_4^{(1)}$ is $D_5^{(2)}$. We can parametrize the perfect crystal $B(l\Lambda_4)$ and its limit $B^{4,\infty}$ as follows:
$$B^{4,l}=\left\{(b_{ij})_{1\leq i\leq 4,i\leq j\leq i+4} \, \, \, \, \vline \begin{array}{l@{}}
(b_{ij})\in\mathbb{Z}^{\geq 0}, \sum_{j=i}^{i+4} b_{ij}=l, 1\leq i\leq 4\\
\sum_{j=i}^t b_{ij}\geq\sum_{j=i+1}^{t+1}b_{i+1,j}, 1\leq i\leq 3, 1\leq t\leq 5\\
b_{11}=b_{44}+b_{45}+b_{46}+b_{47}\\ b_{48}=b_{12}+b_{13}+b_{14}+b_{15},\\
b_{22}=b_{44}+b_{45}+b_{46},\,\,\, b_{37}=b_{13}+b_{14}+b_{15},\\
b_{22}+b_{23}=b_{33}+b_{34}+b_{35},\\b_{36}+b_{37}=b_{24}+b_{25}+b_{26},\\
b_{33}=b_{44}+b_{45},\,\,\,b_{26}=b_{14}+b_{15}
\end{array}\right\}$$
$$B^{4,\infty}=\left\{(b_{ij})_{1\leq i\leq 4,i\leq j\leq i+4} \, \, \, \, \vline \begin{array}{l@{}}
(b_{ij})\in\mathbb{Z}, \sum_{j=i}^{i+4} b_{ij}=0, 1\leq i\leq 4\\
b_{11}=b_{44}+b_{45}+b_{46}+b_{47},\\ b_{48}=b_{12}+b_{13}+b_{14}+b_{15},\\
b_{22}=b_{44}+b_{45}+b_{46},\,\,\, b_{37}=b_{13}+b_{14}+b_{15},\\
b_{22}+b_{23}=b_{33}+b_{34}+b_{35},\\b_{36}+b_{37}=b_{24}+b_{25}+b_{26},\\
b_{33}=b_{44}+b_{45},\,\,\,b_{26}=b_{14}+b_{15}
\end{array}\right\}$$

We give the $\Tilde{f}_k$ actions and $\phi_k,\varepsilon_k, \text{wt}_k$ functions.

We have the following actions and relations of this perfect crystal:
For $0\leq k\leq 4$, $b\in B^{4,\infty}$ $(B^{4,l})$, $\Tilde{e}_k(b)=(b_{ij}')$, where

For $k=0$ $b_{11}'=b_{11}-1$, $b_{48}'=b_{48}+1$ and 

\[\begin{cases}
    b_{12}'=b_{12}+1,\,\,\,b_{23}'=b_{23}+1,\,\,\,b_{34}'=b_{34}+1,\,\,\,b_{22}'=b_{22}-1,\,\,\,b_{33}'=b_{33}-1, b_{44}'=b_{44}-1 &(A1)\\
    b_{13}'=b_{13}+1,\,\,\,b_{24}'=b_{24}+1,\,\,\,b_{37}'=b_{37}+1,\,\,\,b_{23}'=b_{23}-1,\,\,\,b_{34}'=b_{34}-1, b_{47}'=b_{47}-1& (A2)\\
    b_{12}'=b_{12}+1,\,\,\,b_{24}'=b_{24}+1,\,\,\,b_{36}'=b_{36}+1,\,\,\,b_{22}'=b_{22}-1,\,\,\,b_{34}'=b_{34}-1,\,\,\,b_{46}'=b_{46}-1&(A3)\\
    b_{12}'=b_{12}+1,\,\,\,b_{23}'=b_{23}+1,\,\,\,b_{35}'=b_{35}+1,\,\,\,b_{22}'=b_{22}-1,\,\,\,b_{33}'=b_{33}-1,\,\,\,b_{45}'=b_{45}-1&(A4)\\
  b_{12}'=b_{12}+1,\,\,b_{25}'=b_{25}+1,\,\,b_{36}'=b_{36}+1,\,\,b_{22}'=b_{22}-1,\,\,b_{35}'=b_{35}-1,\,\,b_{46}'=b_{46}-1&(A5)\\
  b_{14}'=b_{14}+1,\,\,b_{26}'=b_{26}+1, b_{37}'=b_{37}+1,\,\,b_{24}'=b_{24}-1,\,\,b_{36}'=b_{36}-1,\,\,b_{47}'=b_{47}-1&(A6)\\
  b_{13}'=b_{13}+1,\,\,b_{25}'=b_{25}+1,\,\,b_{37}'=b_{37}+1,\,\,b_{23}'=b_{23}-1,\,\,b_{35}'=b_{35}-1,\,\,b_{47}'=b_{47}-1&(A7)\\
  b_{15}'=b_{15}+1,\,\,b_{26}'=b_{26}+1,\,\,b_{37}'=b_{37}+1,\,\,b_{25}'=b_{25}-1,\,\,b_{36}'=b_{36}-1,\,\,b_{47}'=b_{47}-1 &(A8)
    \end{cases}\]
\begin{align*}
  (A1)=&b_{12}+b_{35}+b_{36}\geq b_{45}+b_{46}+b_{47},\,\,\, b_{23}+b_{35}\geq b_{45}+b_{46},\,\,b_{34}\geq b_{45},\,\,b_{23}+b_{24}\geq b_{45}+b_{46},\\ &b_{12}+b_{13}+b_{25}\geq b_{45}+b_{46}+b_{47},\,\,b_{12}+b_{24}+b_{36}\geq b_{45}+b_{46}+b_{47},\\&b_{12}+b_{13}+b_{14}\geq b_{45}+b_{46}+b_{47}, (b_{11}>0, b_{22}>0,b_{33}>0,b_{44}>0)  \\
  (A2)=& b_{45}+b_{46}+b_{47}> b_{12}+b_{35}+b_{36},\,\,b_{23}+b_{47}\geq b_{12}+b_{36},\,\,b_{34}+b_{46}+b_{47}\geq b_{12}+b_{35}+b_{36},\\& b_{23}+b_{24}+b_{47}\geq b_{12}+b_{35}+b_{36},\,\,b_{13}+b_{25}\geq b_{35}+b_{36},\,\,b_{24}\geq b_{35},\,\,b_{13}+b_{14}\geq b_{35}+b_{36},\\&b_{23}+b_{24}+b_{47}> b_{12}+b_{35}+b_{36} \text{ or } b_{23}+b_{47}>b_{12}+b_{36}, (b_{11}>0, b_{23}>0,b_{34}>0,b_{47}>0)\\
  (A3)=& b_{45}+b_{46}>b_{23}+b_{35},\,\,b_{12}+b_{36}\geq b_{23}+b_{47},\,\,b_{34}+b_{46}\geq b_{23}+b_{35},\,\,b_{24}\geq b_{35},\\& b_{12}+b_{13}+b_{25}\geq b_{23}+b_{35}+b_{47},\,\,b_{12}+b_{13}+b_{14}\geq b_{23}+b_{35}+b_{47},\\& 
  b_{12}+b_{24}+b_{36}\geq b_{23}+b_{35}+b_{47},\,\,b_{12}+b_{36}\geq b_{23}+b_{47},\,\,b_{24}+b_{36}\geq b_{13}+b_{25}, \\
  &(b_{11},b_{22},b_{34},b_{46}>0)\\
  (A4)=& b_{45}>b_{34},\,\, b_{12}+b_{35}+b_{36}>b_{34}+b_{46}+b_{47},\,\,b_{23}+b_{35}>b_{34}+b_{46},\,\,b_{23}+b_{24}\geq b_{34}+b_{46},\\&b_{12}+b_{13}+b_{25}\geq b_{34}+b_{46}+b_{47},\,\,b_{12}+b_{24}+b_{36}\geq b_{34}+b_{46}+b_{47},\\& b_{12}+b_{13}+b_{14}\geq b_{34}+b_{46}+b_{47}, (b_{11},b_{22},b_{33},b_{45}>0)\\
  (A5)=& b_{45}+b_{46}>b_{23}+b_{24},\,\,b_{12}+b_{35}+b_{36}>b_{23}+b_{24}+b_{47},\,\,b_{35}>b_{24},\,\,b_{34}+b_{46}>b_{23}+b_{24},\\&b_{12}+b_{13}+b_{25}\geq b_{23}+b_{24}+b_{47},\,\,b_{12}+b_{36}\geq b_{23}+b_{47},\,\,b_{12}+b_{13}+b_{14}\geq b_{23}+b_{24}+b_{47},\\&b_{12}+b_{36}\geq b_{23}+b_{47},\text{ or } b_{12}+b_{35}+b_{36}\geq b_{34}+b_{46}+b_{47}, (b_{11},b_{22},b_{35},b_{46}>0)\\
\end{align*}
\begin{align*}
 (A6)=&b_{45}+b_{46}+b_{47}>b_{12}+b_{13}+b_{25},\,\,b_{35}+b_{36}>b_{13}+b_{25},\,\,b_{23}+b_{35}+b_{47}>b_{12}+b_{13}+b_{25},\\&b_{34}+b_{46}+b_{47}>b_{12}+b_{13}+b_{25},\,\,b_{23}+b_{24}+b_{47}>b_{12}+b_{13}+b_{25},\,\,b_{24}+b_{36}\geq b_{13}+b_{25},\\&
  b_{14}\geq b_{25}, (b_{11},b_{24},b_{36},b_{47}>0)\\
(A7)=&b_{45}+b_{46}+b_{47}>b_{12}+b_{24}+b_{36},\,\,b_{35}>b_{24},\,\,b_{23}+b_{35}+b_{47}>b_{12}+b_{24}+b_{36},\\& b_{34}+b_{46}+b_{47}>b_{12}+b_{24}+b_{36},\,\,b_{23}+b_{47}>b_{12}+b_{36},\,\,b_{13}+b_{25}>b_{24}+b_{36},\\& b_{13}+b_{14}\geq b_{24}+b_{36},\,\,b_{34}+b_{46}\leq b_{23}+b_{35}, (b_{11},b_{23},b_{35},b_{47}>0)\\
(A8)=& b_{45}+b_{46}+b_{47}>b_{12}+b_{13}+b_{14},\,\,b_{35}+b_{36}>b_{13}+b_{14},\,\,b_{23}+b_{35}+b_{47}>b_{12}+b_{13}+b_{14},\\&b_{34}+b_{46}+b_{47}>b_{12}+b_{13}+b_{14},\,\,b_{23}+b_{24}+b_{47}>b_{12}+b_{13}+b_{14},\,\,b_{25}>b_{14},\\& b_{36}+b_{24}>b_{13}+b_{14}, (b_{11},b_{25},b_{36},b_{47}>0)
\end{align*}
For $k=1$, $b_{11}'=b_{11}+1$, $b_{12}'=b_{12}-1$, $b_{47}'=b_{47}+1$, and $b_{48}'=b_{48}-1$ and if $b\in B^{4,l}$ we have the additional conditions $b_{12},b_{48}>0$.

For $k=2$, \[\begin{cases}
b_{12}'=b_{12}+1,\,\,b_{13}'=b_{13}-1,\,\,b_{36}'=b_{36}+1,\,\,b_{37}'=b_{37}-1 &b_{23}+b_{47}\leq b_{12}+b_{36}, (b_{13},b_{37}>0)\\
b_{22}'=b_{22}+1,\,\,b_{23}'=b_{23}-1,\,\,b_{46}'=b_{46}+1,\,\,b_{47}'=b_{47}-1&b_{23}+b_{47}>b_{12}+b_{36}, (b_{23},b_{47}>0)
\end{cases}\]

For $k=3$, \[\begin{cases}
b_{33}'=b_{33}+1,\,\,b_{34}'=b_{34}-1,\,\,b_{45}'=b_{45}+1,\,\,b_{46}'=b_{46}-1 &b_{25}+b_{35}<b_{36}+b_{46}\\&b_{35}<b_{46},\,\,(b_{34},b_{46}>0)\\
b_{23}'=b_{23}+1,\,\,b_{24}'=b_{24}-1,\,\,b_{35}'=b_{35}+1,\,\,b_{36}'=b_{36}-1 & b_{35}\geq b_{46},\,\,b_{13}<b_{24},\,\,(b_{24},b_{36}>0)\\
b_{13}'=b_{13}+1,\,\,b_{14}'=b_{14}-1,\,\,b_{25}'=b_{25}+1,\,\,b_{26}'=b_{26}-1& b_{13}\geq b_{24}\\&b_{25}+b_{35},\,\,\geq b_{36}+b_{46},\,\,(b_{14},b_{26}>0)
\end{cases}\]

For $k=4$, \[\begin{cases}
b_{14}'=b_{14}+1,\,\,b_{15}'=b_{15}-1 & b_{14}\geq b_{25},\,\,b_{14}+b_{24}\geq b_{25}+b_{35},\,\, b_{14}+b_{24}+b_{34}\geq b_{25}+b_{35}+b_{45},\\&(b_{15}>0)\\
b_{24}'=b_{24}+1,\,\,b_{25}'=b_{25}-1& b_{14}<b_{25},\,\, b_{24}\geq b_{35},\,\,b_{24}+b_{34}\geq b_{35}+b_{45},\,\,(b_{25}>0)\\
b_{34}'=b_{34}+1,\,\,b_{35}'=b_{35}-1 & b_{14}+b_{24}<b_{25}+b_{35},\,\,b_{24}<b_{35},\,\,b_{34}\geq b_{45},\,\,(b_{35}>0)\\
b_{44}'=b_{44}+1,\,\,b_{45}'=b_{45}-1& b_{14}+b_{24}+b_{34}<b_{25}+b_{35}+b_{45},\,\,b_{24}+b_{34}<b_{35}+b_{45},\,\,b_{34}<b_{45},\\&(b_{45}>0)
\end{cases}\]

For $0\leq k\leq 4$, $b\in B^{4,\infty}$, $\Tilde{f}_k(b)=(b_{ij}')$, where
For $k=0$ $b_{11}'=b_{11}+1$, $b_{48}'=b_{48}-1$ and 

\[\begin{cases}
    b_{12}'=b_{12}-1,\,\,\,b_{23}'=b_{23}-1,\,\,\,b_{34}'=b_{34}-1,\,\,\,b_{22}'=b_{22}+1,\,\,\,b_{33}'=b_{33}+1, b_{44}'=b_{44}+1 &(B1)\\
    b_{13}'=b_{13}-1,\,\,\,b_{24}'=b_{24}-1,\,\,\,b_{37}'=b_{37}-1,\,\,\,b_{23}'=b_{23}+1,\,\,\,b_{34}'=b_{34}+1, b_{47}'=b_{47}+1& (B2)\\
    b_{12}'=b_{12}-1,\,\,\,b_{24}'=b_{24}-1,\,\,\,b_{36}'=b_{36}-1,\,\,\,b_{22}'=b_{22}+1,\,\,\,b_{34}'=b_{34}+1,\,\,\,b_{46}'=b_{46}+1&(B3)\\
    b_{12}'=b_{12}-1,\,\,\,b_{23}'=b_{23}-1,\,\,\,b_{35}'=b_{35}-1,\,\,\,b_{22}'=b_{22}+1,\,\,\,b_{33}'=b_{33}+1,\,\,\,b_{45}'=b_{45}+1&(B4)\\
  b_{12}'=b_{12}-1,\,\,b_{25}'=b_{25}-1,\,\,b_{36}'=b_{36}-1,\,\,b_{22}'=b_{22}+1,\,\,b_{35}'=b_{35}+1,\,\,b_{46}'=b_{46}+1&(B5)\\
  b_{14}'=b_{14}-1,\,\,b_{26}'=b_{26}-1, b_{37}'=b_{37}-1,\,\,b_{24}'=b_{24}+1,\,\,b_{36}'=b_{36}+1,\,\,b_{47}'=b_{47}+1&(B6)\\
  b_{13}'=b_{13}-1,\,\,b_{25}'=b_{25}-1,\,\,b_{37}'=b_{37}-1,\,\,b_{23}'=b_{23}+1,\,\,b_{35}'=b_{35}+1,\,\,b_{47}'=b_{47}+1&(B7)\\
  b_{15}'=b_{15}-1,\,\,b_{26}'=b_{26}-1,\,\,b_{37}'=b_{37}-1,\,\,b_{25}'=b_{25}+1,\,\,b_{36}'=b_{36}+1,\,\,b_{47}'=b_{47}+1 &(B8)
    \end{cases}\]
\begin{align*}
  (B1)=&b_{12}+b_{35}+b_{36}> b_{45}+b_{46}+b_{47},\,\,\, b_{23}+b_{35}> b_{45}+b_{46},\,\,b_{34}> b_{45},\,\,b_{23}+b_{24}> b_{45}+b_{46},\\ &b_{12}+b_{13}+b_{25}> b_{45}+b_{46}+b_{47},\,\,b_{12}+b_{24}+b_{36}> b_{45}+b_{46}+b_{47},\\&b_{12}+b_{13}+b_{14}> b_{45}+b_{46}+b_{47},\,\,(b_{48},b_{12},b_{23},b_{34}>0)  \\
  (B2)=& b_{45}+b_{46}+b_{47}\geq b_{12}+b_{35}+b_{36},\,\,b_{23}+b_{47}> b_{12}+b_{36},\,\,b_{34}+b_{46}+b_{47}> b_{12}+b_{35}+b_{36},\\& b_{23}+b_{24}+b_{47}> b_{12}+b_{35}+b_{36},\,\,b_{13}+b_{25}> b_{35}+b_{36},\,\,b_{24}> b_{35},\,\,b_{13}+b_{14}> b_{35}+b_{36},\\&b_{23}+b_{24}+b_{47}\geq b_{12}+b_{35}+b_{36} \text{ or } b_{23}+b_{47}
  \geq b_{12}+b_{36},\,\,(b_{48},b_{13},b_{24},b_{37}>0)
  \end{align*}
  \begin{align*}
  (B3)=& b_{45}+b_{46}\geq b_{23}+b_{35},\,\,b_{12}+b_{36}> b_{23}+b_{47},\,\,b_{34}+b_{46}> b_{23}+b_{35},\,\,b_{24}> b_{35},\\& b_{12}+b_{13}+b_{25}> b_{23}+b_{35}+b_{47},\,\,b_{12}+b_{13}+b_{14}> b_{23}+b_{35}+b_{47},\\& 
  b_{12}+b_{24}+b_{36}> b_{23}+b_{35}+b_{47},\,\,b_{12}+b_{36}> b_{23}+b_{47},\,\,b_{24}+b_{36}> b_{13}+b_{25},\,\,(b_{48},\\&b_{12},b_{24},b_{36}>0)\\
  (B4)=& b_{45}\geq b_{34},\,\, b_{12}+b_{35}+b_{36}\geq b_{34}+b_{46}+b_{47},\,\,b_{23}+b_{35}\geq b_{34}+b_{46},\,\,b_{23}+b_{24}> b_{34}+b_{46},\\&b_{12}+b_{13}+b_{25}> b_{34}+b_{46}+b_{47},\,\,b_{12}+b_{24}+b_{36}> b_{34}+b_{46}+b_{47},\\& b_{12}+b_{13}+b_{14}> b_{34}+b_{46}+b_{47},\,\,(b_{48},b_{12},b_{23},b_{35}>0)\\
  (B5)=& b_{45}+b_{46}\geq b_{23}+b_{24},\,\,b_{12}+b_{35}+b_{36}\geq b_{23}+b_{24}+b_{47},\,\,b_{35}\geq b_{24},\,\,b_{34}+b_{46}\geq b_{23}+b_{24},\\&b_{12}+b_{13}+b_{25}> b_{23}+b_{24}+b_{47},\,\,b_{12}+b_{36}> b_{23}+b_{47},\,\,b_{12}+b_{13}+b_{14}> b_{23}+b_{24}+b_{47},\\&b_{12}+b_{36}> b_{23}+b_{47},\text{ or } b_{12}+b_{35}+b_{36}> b_{34}+b_{46}+b_{47},\,\,(b_{48},b_{12},b_{25},b_{36}>0)\\
  (B6)=&b_{45}+b_{46}+b_{47}\geq b_{12}+b_{13}+b_{25},\,\,b_{35}+b_{36}\geq b_{13}+b_{25},\,\,b_{23}+b_{35}+b_{47}\geq b_{12}+b_{13}+b_{25},\\&b_{34}+b_{46}+b_{47}\geq b_{12}+b_{13}+b_{25},\,\,b_{23}+b_{24}+b_{47}\geq b_{12}+b_{13}+b_{25},\,\,b_{24}+b_{36}> b_{13}+b_{25},\\&
  b_{14}> b_{25},\,\,(b_{48},b_{14},b_{26},b_{37}>0)\\
(B7)=&b_{45}+b_{46}+b_{47}\geq b_{12}+b_{24}+b_{36},\,\,b_{35}\geq b_{24},\,\,b_{23}+b_{35}+b_{47}\geq b_{12}+b_{24}+b_{36},\\& b_{34}+b_{46}+b_{47}\geq b_{12}+b_{24}+b_{36},\,\,b_{23}+b_{47}\geq b_{12}+b_{36},\,\,b_{13}+b_{25}\geq b_{24}+b_{36},\\& b_{13}+b_{14}> b_{24}+b_{36},\,\,b_{34}+b_{46}< b_{23}+b_{35},\,\,(b_{48},b_{13},b_{25},b_{37}>0)\\
(B8)=& b_{45}+b_{46}+b_{47}\geq b_{12}+b_{13}+b_{14},\,\,b_{35}+b_{36}\geq b_{13}+b_{14},\,\,b_{23}+b_{35}+b_{47}\geq b_{12}+b_{13}+b_{14},\\&b_{34}+b_{46}+b_{47}\geq b_{12}+b_{13}+b_{14},\,\,b_{23}+b_{24}+b_{47}\geq b_{12}+b_{13}+b_{14},\,\,b_{25}\geq b_{14},\\& b_{36}+b_{24}\geq b_{13}+b_{14},\,\,(b_{48},b_{15},b_{26},b_{37}>0)
\end{align*}
For $k=1$, $b_{11}'=b_{11}-1$, $b_{12}'=b_{12}+1$, $b_{47}'=b_{47}-1$, and $b_{48}'=b_{48}+1$ and if $b\in B^{4,l}$ we have the additional conditions $b_{11},b_{47}>0$.

For $k=2$, \[\begin{cases}
b_{12}'=b_{12}-1,\,\,b_{13}'=b_{13}+1,\,\,b_{36}'=b_{36}-1,\,\,b_{37}'=b_{37}+1 &b_{23}+b_{47}< b_{12}+b_{36}\\
&(b_{12},b_{36}>0)\\
b_{22}'=b_{22}-1,\,\,b_{23}'=b_{23}+1,\,\,b_{46}'=b_{46}-1,\,\,b_{47}'=b_{47}+1 & b_{23}+b_{47}\geq b_{12}+b_{36}\\
&(b_{22},b_{46}>0)
\end{cases}\]

For $k=3$, \[\begin{cases}
b_{33}'=b_{33}-1,\,\,b_{34}'=b_{34}+1,\,\,b_{45}'=b_{45}-1,\,\,b_{46}'=b_{46}+1 &b_{25}+b_{35}\leq b_{36}+b_{46}\\&b_{35}\leq b_{46},\,\,(b_{33},b_{45}>0)\\
b_{23}'=b_{23}-1,\,\,b_{24}'=b_{24}+1,\,\,b_{35}'=b_{35}-1,\,\,b_{36}'=b_{36}+1 & b_{35}> b_{46},\,\,b_{13}\leq b_{24},\,\,(b_{23},b_{35}>0)\\
b_{13}'=b_{13}-1,\,\,b_{14}'=b_{14}+1,\,\,b_{25}'=b_{25}-1,\,\,b_{26}'=b_{26}+1 &b_{25}+b_{35}> b_{36}+b_{46}\\&b_{13}> b_{24},\,\,(b_{13},b_{25}>0)
\end{cases}\]

For $k=4$, \[\begin{cases}
b_{14}'=b_{14}-1,\,\,b_{15}'=b_{15}+1 & b_{14}> b_{25},\,\,b_{14}+b_{24}> b_{25}+b_{35},\,\, b_{14}+b_{24}+b_{34}> b_{25}+b_{35}+b_{45}\\
&(b_{14}>0)\\
b_{24}'=b_{24}-1,\,\,b_{25}'=b_{25}+1& b_{14}\leq b_{25},\,\, b_{24}> b_{35},\,\,b_{24}+b_{34}> b_{35}+b_{45},\,\,(b_{24}>0)\\
b_{34}'=b_{34}-1,\,\,b_{35}'=b_{35}+1 & b_{14}+b_{24}\leq b_{25}+b_{35},\,\,b_{24}\leq b_{35},\,\,b_{34}> b_{45},\,\,(b_{34}>0)\\
b_{44}'=b_{44}-1,\,\,b_{45}'=b_{45}+1& b_{14}+b_{24}+b_{34}\leq b_{25}+b_{35}+b_{45},\,\,b_{24}+b_{34}\leq b_{35}+b_{45},\,\,b_{34}\leq b_{45},\\
&(b_{44}>0)
\end{cases}\]

Then we have the following formulas for $\varphi_i(b)$.
\[\varphi_0(b)=\begin{cases}
-b_{11}-\min\{b_{45}+b_{46}+b_{47},b_{34}+b_{46}+b_{47},b_{12}+b_{35}+b_{36}, b_{23}+b_{35}& b\in B^{4,\infty}\\ +b_{47},b_{23}+b_{24}+b_{47},b_{12}+b_{13}+b_{25},b_{12}+b_{24}+b_{36},b_{12}+b_{13}+b_{14}\}\\
l-b_{11}-\min\{b_{45}+b_{46}+b_{47},b_{34}+b_{46}+b_{47},b_{12}+b_{35}+b_{36}, b_{23}+b_{35}& b\in B^{4,l}\\+b_{47}, b_{23}+b_{24}+b_{47}, b_{12}+b_{13}+b_{25},b_{12}+b_{24}+b_{36},b_{12}+b_{13}+b_{14}\} 
\end{cases}\]

\begin{align*}
 &\varphi_1(b)=b_{47} \\ 
 &\varphi_2(b)=b_{46}+\max\{b_{12}-b_{23},0\}\\
 &\varphi_3(b)=b_{45}+\max\{b_{13}-b_{24},b_{13}+b_{23}-b_{24}-b_{34},0\}\\
 &\varphi_4(b)=b_{44}+\max\{b_{14}-b_{25},b_{14}+b_{24}-b_{25}-b_{35},b_{14}+b_{24}+b_{34}-b_{25}-b_{35}-b_{45},0\}
\end{align*}
And the following formulas for $\varepsilon_i(b)$.
\[\varepsilon_0(b)=\begin{cases}
-b_{48}-\min\{b_{45}+b_{46}+b_{47},b_{34}+b_{46}+b_{47},b_{12}+b_{35}+b_{36}, b_{23}+b_{35}& b\in B^{4,\infty}\\ +b_{47},b_{23}+b_{24}+b_{47},b_{12}+b_{13}+b_{25},b_{12}+b_{24}+b_{36},b_{12}+b_{13}+b_{14}\}\\
l-b_{48}-\min\{b_{45}+b_{46}+b_{47},b_{34}+b_{46}+b_{47},b_{12}+b_{35}+b_{36}, b_{23}+b_{35}& b\in B^{4,l}\\+b_{47}, b_{23}+b_{24}+b_{47}, b_{12}+b_{13}+b_{25},b_{12}+b_{24}+b_{36},b_{12}+b_{13}+b_{14}\} 
\end{cases}\]

\begin{align*}
 &\varepsilon_1(b)=b_{12} \\ 
 &\varepsilon_2(b)=b_{13}+\max\{b_{23}-b_{12},0\}\\
 &\varepsilon_3(b)=b_{14}+\max\{b_{24}-b_{13},b_{24}+b_{34}-b_{13}-b_{23},0\}\\
 &\varepsilon_4(b)=b_{15}+\max\{b_{25}-b_{14},b_{25}+b_{35}-b_{24}-b_{14},b_{25}+b_{35}+b_{45}-b_{14}-b_{24}-b_{34},0\}
\end{align*}
Finally, the formulas for $\text{wt}_i(b)$ are:
\begin{align*}
    &\text{wt}_0(b)=b_{48}-b_{11}\\
    &\text{wt}_1(b)=b_{47}-b_{12}\\
    &\text{wt}_2(b)=b_{46}-b_{13}+b_{12}-b_{23}\\
    &\text{wt}_3(b)=b_{45}-b_{14}+b_{13}-b_{24}+b_{23}-b_{34}\\
    &\text{wt}_4(b)=b_{44}-b_{15}+b_{14}-b_{25}+b_{24}-b_{35}+b_{34}-b_{45}
\end{align*}
In [], we showed this is a limit of a coherent family of perfect crystals. Now we will ultradiscretize the geometric crystal corresponding to $C_4^{(1)}$. If $\mathcal{X}=\mathcal{UD}(\mathcal{V})$, then $\mathcal{X}=\mathbb{Z}^{10}$ as sets and $\mathcal{UD}(\mathcal{V})$ is equipped with the following functions:
\begin{align*}
   \text{wt}_i(x)=\begin{cases}
   -2x_{11} & i=0\\
   2x_{11}-x_{21}-x_{22} & i=1\\
   2x_{21}+2x_{22}-x_{11}-x_{31}-x_{32}-x_{33} & i=2\\
   2x_{31}+2x_{32}+2x_{33}-x_{21}-x_{22}-x_{41}-x_{42}-x_{43}-x_{44} & i=3\\
   2x_{41}+2x_{42}+2x_{43}+2x_{44}-2x_{31}-2x_{32}-2x_{33}&i=4
   \end{cases}
\end{align*}
\begin{align*}
    \varepsilon_i(x)=\begin{cases}
    \max\{x_{41},2x_{11}-x_{44}, 2x_{31}-x_{42},2x_{11}+x_{42}-2x_{22},& i=0\\ 2x_{21}+x_{42}-2x_{32},2x_{21}-x_{43},2x_{11}+x_{43}-2x_{33},2x_{11}+2x_{32}-2x_{22}-x_{43}\}\\
    x_{22}-x_{11}&i=1\\
    \max\{x_{33}-x_{22},x_{33}+x_{32}+x_{11}-2x_{22}-x_{21}\}&i=2\\
    \max\{x_{44}-x_{33},x_{44}+x_{43}+x_{22}-2x_{33}-x_{32},&i=3\\ x_{44}+x_{43}+x_{42}+x_{22}+x_{21}-2x_{33}-2x_{32}-x_{31}\}\\
    \max\{-x_{44},2x_{33}-x_{43}-2x_{44},2x_{32}+2x_{33}-2x_{44}-2x_{43}-x_{42},& i=4\\ 2x_{31}+2x_{32}+2x_{33}-2x_{44}-2x_{43}-2x_{42}-x_{41}\}
    \end{cases}
\end{align*}
For $i=1$, we compute the ultra-discretization of $e_1^c$
\begin{align*}\mathcal(UD)(e_1^c(x))=(x_{44},x_{33},x_{43},x_{22},x_{32},x_{42},x_{11}+c,x_{21},x_{31},x_{41})\end{align*}

By restricting $c$ to -1 and 1, we get $\Tilde{f_1}$ and $\Tilde{e_1}$.
\begin{align*}&\Tilde{e_1}(x)=(x_{44},x_{33},x_{43},x_{22},x_{32},x_{42},x_{11}+1,x_{21},x_{31},x_{41})\\
&\Tilde{f_1}(x)=(x_{44},x_{33},x_{43},x_{22},x_{32},x_{42},x_{11}-1,x_{21},x_{31},x_{41})
\end{align*}
Now for $i=2$
\begin{align*}
    &\mathcal{UD}(e_2^c(x))=(x_{44},x_{33},x_{43},x_{22}+\max\{x_{22}+x_{21}+c, x_{11}+x_{32}\}-\max\{x_{22}+x_{21}, x_{11}+x_{32}\},\\&x_{32},x_{42},x_{11},x_{21}+c+\max\{x_{22}+x_{21}, x_{11}+x_{32}\}-\max\{x_{22}+x_{21}+c, x_{11}+x_{32}\},x_{31},x_{41})
\end{align*}
Again, restricting $c$ we obtain $\Tilde{f}_2$ and $\Tilde{e}_2$
\begin{align*}
    \Tilde{e}_2(x)=\begin{cases}
    (x_{44},x_{33},x_{43},x_{22}+1,x_{32},x_{42},x_{11},x_{21},x_{31},x_{41}) & x_{21}+x_{22}\geq x_{11}+x_{32}\\
    (x_{44},x_{33},x_{43},x_{22},x_{32},x_{42},x_{11},x_{21}+1,x_{31},x_{41}) & x_{21}+x_{22}< x_{11}+x_{32}
    \end{cases}
\end{align*}
\begin{align*}
    \Tilde{f}_2(x)=\begin{cases}
    (x_{44},x_{33},x_{43},x_{22}-1,x_{32},x_{42},x_{11},x_{21},x_{31},x_{41}) & x_{21}+x_{22}> x_{11}+x_{32}\\
    (x_{44},x_{33},x_{43},x_{22},x_{32},x_{42},x_{11},x_{21}-1,x_{31},x_{41}) & x_{21}+x_{22}\leq x_{11}+x_{32}
    \end{cases}
\end{align*}
For $i=3$, we have
\begin{align*}
    &\mathcal{UD}(e_3^c(x))=(x_{44},x_{33}+\max\{x_{44}-x_{33}+c, x_{44}+x_{43}+x_{22}-2x_{33}-x_{32},x_{44}+x_{43}+x_{42}\\&+x_{22}+x_{21}-2x_{33}-2x_{32}-x_{31}\}-\max\{x_{44}-x_{33}, x_{44}+x_{43}+x_{22}-2x_{33}-x_{32},x_{44}+x_{43}+x_{42}\\&+x_{22}+x_{21}-2x_{33}-2x_{32}-x_{31}\},x_{43},x_{22},x_{32}+\max\{x_{44}-x_{33}+c, x_{44}+x_{43}+x_{22}-2x_{33}\\&-x_{32}+c, x_{44}+x_{43}+x_{42}+x_{22}+x_{21}-2x_{33}-2x_{32}-x_{31}\}-\max\{x_{44}-x_{33}+c, x_{44}+x_{43}+x_{22}\\&-2x_{33}-x_{32},x_{44}+x_{43}+x_{42}+x_{22}+x_{21}-2x_{33}-2x_{32}-x_{31}\},x_{42},x_{11},x_{21},x_{31}+c+\max\{x_{44}\\&-x_{33}, x_{44}+x_{43}+x_{22}-2x_{33}-x_{32},x_{44}+x_{43}+x_{42}+x_{22}+x_{21}-2x_{33}-2x_{32}-x_{31}\}-\max\{x_{44}\\&-x_{33}+c, x_{44}+x_{43}+x_{22}-2x_{33}-x_{32}+c,x_{44}+x_{43}+x_{42}+x_{22}+x_{21}-2x_{33}-2x_{32}-x_{31}\},x_{41})
\end{align*}
Restricting $c$ we obtain $\Tilde{f}_3$ and $\Tilde{e}_3$
\begin{align*}
    \Tilde{e}_3(x)=\begin{cases}
    (x_{44},x_{33}+1,x_{43},x_{22},x_{32},x_{42},x_{11},x_{21},x_{31},x_{41}) & 0\geq x_{43}+x_{22}-x_{33}-x_{32},\\& 0\geq x_{43}+x_{42}+x_{22}+x_{21}\\&-x_{33}-2x_{32}-x_{31}\\
    (x_{44},x_{33},x_{43},x_{22},x_{32}+1,x_{42},x_{11},x_{21},x_{31},x_{41}) & 0< x_{43}+x_{22}-x_{33}-x_{32},\\& 0\geq x_{42}+x_{21}-x_{32}-x_{31}\\
    (x_{44},x_{33},x_{43},x_{22},x_{32},x_{42},x_{11},x_{21},x_{31}+1,x_{41}) & 0< x_{43}+x_{42}+x_{22}+x_{21}\\&-x_{33}-2x_{32}-x_{31}\\&0< x_{42}+x_{21}-x_{32}-x_{31}
    \end{cases}
\end{align*}
\begin{align*}
  \Tilde{f}_3(x)=\begin{cases}
    (x_{44},x_{33}-1,x_{43},x_{22},x_{32},x_{42},x_{11},x_{21},x_{31},x_{41}) & 0> x_{43}+x_{22}-x_{33}-x_{32},\\& 0> x_{43}+x_{42}+x_{22}+x_{21}\\&-x_{33}-2x_{32}-x_{31}\\
    (x_{44},x_{33},x_{43},x_{22},x_{32}-1,x_{42},x_{11},x_{21},x_{31},x_{41}) & 0\leq x_{43}+x_{22}-x_{33}-x_{32},\\& 0> x_{42}+x_{21}-x_{32}-x_{31}\\
    (x_{44},x_{33},x_{43},x_{22},x_{32},x_{42},x_{11},x_{21},x_{31}-1,x_{41}) & 0\leq x_{43}+x_{42}+x_{22}+x_{21}\\&-x_{33}-2x_{32}-x_{31}\\&0\leq x_{42}+x_{21}-x_{32}-x_{31}
    \end{cases}
\end{align*}

For $i=4$, let 
\begin{align*}
   \mathcal{UD}(c_4)=\Tilde{c}_4=&\max\{-x_{44},2x_{33}-2x_{44}-x_{43},2x_{32}+2x_{33}-2x_{44}-2x_{43}-x_{42},\\&2x_{31}+2x_{32}+2x_{33}-2x_{44}-2x_{43}-2x_{42}-x_{41}\} \\
    \mathcal{UD}(c_{41})=\Tilde{c}_{41}=&\max\{c-x_{44},2x_{33}-2x_{44}-x_{43},2x_{32}+2x_{33}-2x_{44}-2x_{43}-x_{42},\\&2x_{31}+2x_{32}+2x_{33}-2x_{44}-2x_{43}-2x_{42}-x_{41}\}\\ 
     \mathcal{UD}(c_{42})=\Tilde{c}_{42}=&\max\{c-x_{44},c+2x_{33}-2x_{44}-x_{43},2x_{32}+2x_{33}-2x_{44}-2x_{43}-x_{42},\\&2x_{31}+2x_{32}+2x_{33}-2x_{44}-2x_{43}-2x_{42}-x_{41}\}\\ 
      \mathcal{UD}(c_{43})=\Tilde{c}_{43}=&\max\{c-x_{44},c+2x_{33}-2x_{44}-x_{43},c+2x_{32}+2x_{33}-2x_{44}-2x_{43}-x_{42},\\&2x_{31}+2x_{32}+2x_{33}-2x_{44}-2x_{43}-2x_{42}-x_{41}\}
\end{align*}

Now we have 
$$\mathcal{UD}(e_4^c(x))=(\Tilde{c}_{41}-\Tilde{c}_4+x_{44},x_{33},\Tilde{c}_{42}-\Tilde{c}_{41}+x_{43},x_{22},x_{32},\Tilde{c}_{43}-\Tilde{c}_{42}+x_{42},x_{11},x_{21},x_{31},c+\Tilde{c}_4-\Tilde{c}_{43}+x_{41})$$

Let $L=-x_{44}$, $M=2x_{33}-2x_{44}-x_{43}$, $N=2x_{32}+2x_{33}-2x_{44}-2x_{43}-x_{42}$, and $P=2x_{31}+2x_{32}+2x_{33}-2x_{44}-2x_{43}-2x_{42}-x_{41}$. Restricting $c$ we obtain $\Tilde{e}_4$ and $\Tilde{f}_4$
\begin{align*}
    \Tilde{e}_4(x)=\begin{cases}
    (x_{44}+1,x_{33},x_{43},x_{22},x_{32},x_{42},x_{11},x_{21},x_{31},x_{41})& L\geq M, L\geq N, L\geq P\\
    (x_{44},x_{33},x_{43}+1,x_{22},x_{32},x_{42},x_{11},x_{21},x_{31},x_{41})& M>L, M\geq N, M\geq P\\
    (x_{44},x_{33},x_{43},x_{22},x_{32},x_{42}+1,x_{11},x_{21},x_{31},x_{41})& N>L, N>M, N\geq P\\
    (x_{44},x_{33},x_{43},x_{22},x_{32},x_{42},x_{11},x_{21},x_{31},x_{41}+1)& P>L, P>M, P>N\\
    \end{cases}
\end{align*}

\begin{align*}
    \Tilde{f}_4(x)=\begin{cases}
    (x_{44}-1,x_{33},x_{43},x_{22},x_{32},x_{42},x_{11},x_{21},x_{31},x_{41})& L>M, L>N, L>P\\
    (x_{44},x_{33},x_{43}-1,x_{22},x_{32},x_{42},x_{11},x_{21},x_{31},x_{41})& M\geq L, M>N, M> P\\
    (x_{44},x_{33},x_{43},x_{22},x_{32},x_{42}-1,x_{11},x_{21},x_{31},x_{41})& N\geq L, N\geq M, N>P\\
    (x_{44},x_{33},x_{43},x_{22},x_{32},x_{42},x_{11},x_{21},x_{31},x_{41}-1)& P\geq L, P\geq M, P\geq N\\
    \end{cases}
\end{align*}

Finally, we need to find $\Tilde{e}_0(x)$ and $\Tilde{f}_0(x)$.
Let $A=2x_{11}-2x_{22}+x_{42}$, $B=2x_{21}-2x_{32}+x_{42}$, $C=2x_{31}-x_{42}$, $X=2x_{21}-x_{43}$, $Y=2x_{11}-2x_{33}+x_{43}$, and $Z=2x_{11}+2x_{32}-2x_{22}-x_{43}$. By simplifying the inequalities $A>B$ and $X>Z$, we can see that these inequalities are mutually exclusive. 
\begin{align*}
\mathcal{UD}(e_0^c(x_{44})=&x_{44}-c+\max\{x_{41}+c, \max\{A,B,C\}+c, \max\{X,Y,Z\}+c, 2x_{11}-x_{44}\}\\&-\max\{x_{41}, \max\{A,B,C\}, \max\{X,Y,Z\}, 2x_{11}-x_{44}\} \\
\mathcal{UD}(e_0^c(x_{33})=&x_{33}-c-\max\{x_{41}, \max\{A,B,C\}, \max\{X,Y,Z\}, 2x_{11}-x_{44}\}\\&+\max\{x_{41}+c, \max\{A,B,C\}+c, \max\{X,Y,Z\}+c, 2x_{11}-x_{44}\}\\
\mathcal{UD}(e_0^c(x_{43})=&x_{43}+\max\{x_{41}+c, \max\{A,B,C\}+c, \max\{X,Y,Z\}, 2x_{11}-x_{44}\}\\&-\max\{x_{41}+c, \max\{A,B,C\}+c, \max\{X,Y,Z\}+c, 2x_{11}-x_{44}\}\\
\mathcal{UD}(e_0^c(x_{22})=&x_{22}-c-\max\{x_{41}, \max\{A,B,C\}, \max\{X,Y,Z\}, 2x_{11}-x_{44}\}\\&+\max\{x_{41}+c, \max\{A,B,C\}+c, \max\{X,Y,Z\}+c, 2x_{11}-x_{44}\}\\
\mathcal{UD}(e_0^c(x_{32}))=&x_{32}+\max\{x_{41}+c, \max\{A,B,C\}, \max\{X,Y,Z\}, 2x_{11}-x_{44}\}\\&-\max\{x_{41}+c, \max\{A,B,C\}+c, \max\{X,Y,Z\}+c, 2x_{11}-x_{44}\}\\
\mathcal{UD}(e_0^c(x_{42}))=&x_{42}+\max\{x_{41}+c, \max\{A,B,C\}, \max\{X,Y,Z\}, 2x_{11}-x_{44}\}\\&-\max\{x_{41}+c, \max\{A,B,C\}+c, \max\{X,Y,Z\}, 2x_{11}-x_{44}\}\\
\mathcal{UD}(e_0^c(x_{11}))=&x_{11}-c\\
\mathcal{UD}(e_0^c(x_{21}))=&x_{21}+\max\{x_{41}, \max\{A,B,C\}, \max\{X,Y,Z\}, 2x_{11}-x_{44}\}\\&-\max\{x_{41}+c, \max\{A,B,C\}+c, \max\{X,Y,Z\}+c, 2x_{11}-x_{44}\}\\
\mathcal{UD}(e_0^c(x_{31}))=&x_{31}+\max\{x_{41}, \max\{A,B,C\}, \max\{X,Y,Z\}, 2x_{11}-x_{44}\}\\&-\max\{x_{41}+c, \max\{A,B,C\}, \max\{X,Y,Z\}, 2x_{11}-x_{44}\}\\
\mathcal{UD}(e_0^c(x_{41}))=&x_{41}+\max\{x_{41}, \max\{A,B,C\}, \max\{X,Y,Z\}, 2x_{11}-x_{44}\}\\&-\max\{x_{41}+c, \max\{A,B,C\}, \max\{X,Y,Z\}, 2x_{11}-x_{44}\}
\end{align*}

From this, we can consider the different cases that occur based off which elements are maximal. Coupled with the fact that certain options are mutually exclusive, 8 cases for $\Tilde{e_0}$ and $\Tilde{f_0}$ are possible. These are as follows:
\begin{align*}
    \Tilde{e_0}(x)=\begin{cases}
    (x_{44},x_{33},x_{43},x_{22},x_{32},x_{42},x_{11}-1,x_{21}-1,x_{31}-1,x_{41}-1)& x_{41}\geq A, x_{41}\geq B, x_{41}\geq C, \\&x_{41}\geq X, x_{41}\geq Y, x_{41}\geq Z,\\& x_{41}\geq 2x_{11}-x_{44}\\
    (x_{44},x_{33},x_{43},x_{22},x_{32}-1,x_{42}-1,x_{11}-1,x_{21}-1,x_{31},x_{41})& B>x_{41}, B>A, B\geq C, \\& B\geq X,B\geq Y,B\geq Z \\&B\geq 2x_{11}-x_{44},X\geq Y,\\& X\geq Z\\
     (x_{44},x_{33},x_{43},x_{22}-1,x_{32}-1,x_{42}-1,x_{11}-1,x_{21},x_{31},x_{41})& A>x_{41}, A\geq B,C,X,Y,Z,\\& A\geq 2x_{11}-x_{44}\\& Y>X \text{ or } Z\geq X\\
     (x_{44},x_{33},x_{43},x_{22},x_{32},x_{42}-1,x_{11}-1,x_{21}-1,x_{31}-1,x_{41})& C>x_{41},A,B,\\& C\geq X,Y,Z,2x_{11}-x_{44}\\
     (x_{44},x_{33},x_{43}-1,x_{22},x_{32}-1,x_{42},x_{11}-1,x_{21}-1,x_{31},x_{41})& X>x_{41},A,B,C\\&X\geq Y,Z,2x_{11}-x_{44}\\& B,C\geq A\\
     (x_{44},x_{33}-1,x_{43}-1,x_{22}-1,x_{32},x_{42},x_{11}-1,x_{21},x_{31},x_{41})&Y>x_{41},A,B,C,X\\&Y\geq Z, 2x_{11}-x_{44}\\
     (x_{44},x_{33},x_{43}-1,x_{22}-1,x_{32}-1,x_{42},x_{11}-1,x_{21},x_{31},x_{41})&Z>x_{41},A,B,C,X,Y\\&Z\geq 2x_{11}-x_{44}\\&A>B\text{ or } C\geq B\\
     (x_{44}-1,x_{33}-1,x_{43},x_{22}-1,x_{32},x_{42},x_{11}-1,x_{21},x_{31},x_{41})&2x_{11}-x_{44}>x_{41},A,B,C,\\&2x_{11}-x_{44}>X,Y,Z
    \end{cases}
\end{align*}
\begin{align*}
    \Tilde{f_0}(x)=\begin{cases}
    (x_{44},x_{33},x_{43},x_{22},x_{32},x_{42},x_{11}+1,x_{21}+1,x_{31}+1,x_{41}+1)& x_{41}> A,B,C,X,Y,Z,\\& x_{41}> 2x_{11}-x_{44}\\
    (x_{44},x_{33},x_{43},x_{22},x_{32}+1,x_{42}+1,x_{11}+1,x_{21}+1,x_{31},x_{41})& B\geq x_{41},A, B>C,X,Y,Z \\&B> 2x_{11}-x_{44},X> Y,Z\\
     (x_{44},x_{33},x_{43},x_{22}+1,x_{32}+1,x_{42}+1,x_{11}+1,x_{21},x_{31},x_{41})& A\geq x_{41}, A> B,C,X,Y,Z,\\& A> 2x_{11}-x_{44}\\& Y\geq X \text{ or } Z> X\\
     (x_{44},x_{33},x_{43},x_{22},x_{32},x_{42}+1,x_{11}+1,x_{21}+1,x_{31}+1,x_{41})& C\geq x_{41},A,B,\\& C> X,Y,Z,2x_{11}-x_{44}\\
     (x_{44},x_{33},x_{43}+1,x_{22},x_{32}+1,x_{42},x_{11}+1,x_{21}+1,x_{31},x_{41})& X\geq x_{41},A,B,C\\&X> Y,Z,2x_{11}-x_{44}\\& B,C> A\\
     (x_{44},x_{33}+1,x_{43}+1,x_{22}+1,x_{32},x_{42},x_{11}+1,x_{21},x_{31},x_{41})&Y\geq x_{41},A,B,C,X\\&Y> Z, 2x_{11}-x_{44}\\
     (x_{44},x_{33},x_{43}+1,x_{22}+1,x_{32}+1,x_{42},x_{11}+1,x_{21},x_{31},x_{41})&Z\geq x_{41},A,B,C,X,Y\\&Z> 2x_{11}-x_{44}\\&A\geq B\text{ or } C>B\\
     (x_{44}+1,x_{33}+1,x_{43},x_{22}+1,x_{32},x_{42},x_{11}+1,x_{21},x_{31},x_{41})&2x_{11}-x_{44}\geq x_{41},A,B,C,\\&2x_{11}-x_{44}\geq X,Y,Z
    \end{cases}
\end{align*}

\begin{theorem}
 Let $\Omega:\mathcal{UD}(\mathcal{V})\to B^{4,\infty}$ be given by 
    \begin{align*}
      b_{11}=x_{11},\,\,\,  b_{12}=x_{22}-x_{11},\,\,\,
      b_{13}=x_{33}-x_{22},\,\,\,b_{14}=x_{44}-x_{33}\\
      b_{15}=-x_{44},\,\,\,b_{22}=x_{21},\,\,\,
     b_{23}=x_{32}-x_{21},\,\,\, b_{24}=x_{43}-x_{32}\\
      b_{25}=x_{33}-x_{43},\,\,\,b_{26}=-x_{33},\,\,\,
      b_{33}=x_{31},\,\,\,b_{34}=x_{42}-x_{31}\\
      b_{35}=x_{32}-x_{42},\,\,\,b_{36}=x_{22}-x_{32},\,\,\,
      b_{37}=-x_{22},\,\,\,b_{44}=x_{41}\\
      b_{45}=x_{31}-x_{41},\,\,\,b_{46}=x_{21}-x_{31},\,\,\,
      b_{47}=x_{11}-x_{21},\,\,\,b_{48}=-x_{11}
    \end{align*} and $\Omega^{-1}:B^{4,\infty}\to\mathcal{UD}(\mathcal{V})$ be given by $x_{11}=b_{11},\,\, x_{21}=b_{22},\,\, x_{22}=-b_{37},\,\,\, x_{31}=b_{33},\,\,\,x_{32}=b_{22}+b_{23},\,\,\,x_{33}=-b_{26}$, $x_{41}=b_{44},\,\,\,x_{42}=b_{33}+b_{34},\,\,\,x_{43}=b_{22}+b_{23}+b_{24},\,\,\,x_{44}=-b_{15}$ \\
    $\Omega$ is an isomorphism of crystals. 
\end{theorem}
\begin{proof}
Clearly we see that the map is bijective. We need to prove the following conditions to show that $\Omega$ is an isomorphism:
\begin{enumerate}
\item $\Omega(\Tilde{f_k}(b))=\Tilde{f_k}\Omega(b)$
\item $\Omega(\Tilde{e_k}(b))=\Tilde{e_k}\Omega(b)$
\item $\text{wt}_k(\Omega(b))=\text{wt}_k(b)$
\item $\varepsilon_k(\Omega(b)=\varepsilon_k(b)$
\end{enumerate}
Clearly if these hold, then $\varphi_k(\Omega(b))=\varepsilon_k(\Omega(b))+\text{wt}_k(\Omega(b))=\varepsilon_k(b)+\text{wt}_k(b)=\varphi_k(b)$. Now we check the above conditions for each $k=0,1,2,3,4$.
First we check 1. There are 8 cases for $\Tilde{f_0}(b)$.
\begin{enumerate}
    \item $b_{11}'=b_{11}+1,\,\,\,b_{48}'=b_{48}-1,\,\,\,b_{12}'=b_{12}-1,\,\,\,b_{23}'=b_{23}-1,\,\,\,b_{34}'=b_{34}=1,\,\,\,b_{22}'=b_{22}+1,\,\,\,b_{33}'=b_{33}+1,\,\,\,b_{44}'=b_{44}+1$ if $b_{12}+b_{35}+b_{36}>b_{45}+b_{46}+b_{47}$, $b_{23}+b_{35}>b_{45}+b_{46}$, $b_{34}>b_{45}$, $b_{23}+b_{34}>b_{45}+b_{46}$, $b_{12}+b_{13}+b_{25}>b_{45}+b_{46}+b_{47}$, $b_{12}+b_{24}+b_{36}>b_{45}+b_{46}+b_{47}$, and $b_{12}+b_{13}+b_{14}>b_{45}+b_{46}+b_{47}$. 
    
    Applying the isomorphism, we end up with the action $x_{11}'=x_{11}+1$,$x_{21}'=x_{21}+1,$, $x_{31}'=x_{31}+1$, and $x_{41}'=x_{41}+1$ along with the following conditions: $x_{41}>2x_{11}-2x_{22}+x_{42}$, $x_{41}>2x_{21}-2x_{32}+x_{42}$, $x_{41}>2x_{31}-x_{42}$, $x_{41}>2x_{21}-x_{43}$, $x_{41}>2x_{11}-2x_{33}+x_{43}$, $x_{41}>2x_{11}-2x_{22}+2x_{32}-x_{43}$, and $x_{41}>2x_{11}-x_{44}$. These exactly correspond to the first case of $\Tilde{f}_0(\Omega(b))$.
    \item $b_{11}'=b_{11}+1,\,\,\,b_{48}'=b_{48}-1,\,\,\,b_{13}'=b_{13}-1,\,\,\,b_{24}'=b_{24}-1,\,\,\,b_{37}'=b_{37}-1,\,\,\,b_{23}'=b_{23}+1,\,\,\,b_{34}'=b_{34}+1,\,\,\,b_{47}'=b_{47}+1$ if $b_{12}+b_{35}+b_{36}\leq b_{45}+b_{46}+b_{47}$, $b_{23}+b_{47}>b_{12}+b_{36}$, $b_{34}+b_{46}+b_{47}>b_{12}+b_{35}+b_{36}$, $b_{23}+b_{34}+b_{47}>b_{12}+b_{35}+b_{36}$, $b_{13}+b_{25}>b_{35}+b_{36}$, $b_{24}>b_{35}$, $b_{13}+b_{14}>b_{35}+b_{36}$, and $b_{23}+b_{24}+b_{47}\geq b_{12}+b_{13}+b_{25}$ or $b_{23}+b_{47}\geq b_{12}+b_{36}$ . 
    
    Applying the isomorphism, we end up with the action $x_{11}'=x_{11}+1$,$x_{22}'=x_{22}+1,$, $x_{32}'=x_{32}+1$, and $x_{42}'=x_{42}+1$ along with the following conditions: $x_{41}\leq 2x_{11}-2x_{22}+x_{42}$, $2x_{11}-2x_{22}+x_{42}>2x_{21}-2x_{32}+x_{42}$, $2x_{11}-2x_{22}+x_{42}>2x_{31}-x_{42}$, $2x_{11}-2x_{22}+x_{42}>2x_{21}-x_{43}$, $2x_{11}-2x_{22}+x_{42}>2x_{11}-2x_{33}+x_{43}$, $2x_{11}-2x_{22}+x_{42}>2x_{11}-2x_{22}+2x_{32}-x_{43}$, and $2x_{11}-2x_{22}+x_{42}>2x_{11}-x_{44}$, and $2x_{11}-2x_{33}+x_{43}\geq 2x_{21}-x_{43}$ or $2x_{11}+2x_{32}-2x_{22}-x_{43}\geq 2x_{21}-x_{43}$. These exactly correspond to the second case of $\Tilde{f}_0(\Omega(b))$.
    \item $b_{11}'=b_{11}+1,\,\,\,b_{48}'=b_{48}-1,\,\,\,b_{12}'=b_{12}-1,\,\,\,b_{24}'=b_{24}-1,\,\,\,b_{36}'=b_{36}-1,\,\,\,b_{22}'=b_{22}+1,\,\,\,b_{34}'=b_{34}+1,\,\,\,b_{46}'=b_{46}+1$ if $b_{45}+b_{46}\geq b_{23}+b_{35}$, $b_{12}+b_{36}>b_{23}+b_{47}$, $b_{34}+b_{46}>b_{23}+b_{35}$, $b_{24}>b_{35}$, $b_{12}+b_{13}+b_{25}>b_{23}+b_{35}+b_{47}$,$b_{12}+b_{24}+b_{36}>b_{23}+b_{35}+b_{47}$,   $b_{12}+b_{13}+b_{14}>b_{23}+b_{35}+b_{47}$, $b_{24}+b_{36}>b_{13}+b_{25}$, and $b_{12}+b_{36}>b_{23}+b_{47}$. 
    
    Applying the isomorphism, we end up with the action $x_{11}'=x_{11}+1$,$x_{21}'=x_{21}+1,$, $x_{32}'=x_{32}+1$, and $x_{42}'=x_{42}+1$ along with the following conditions: $x_{41}\leq 2x_{21}-2x_{32}+x_{42}$, $2x_{11}-2x_{22}+x_{42}\leq 2x_{21}-2x_{32}+x_{42}$, $2x_{21}-2x_{32}+x_{42}>2x_{31}-x_{42}$, $2x_{21}-2x_{32}+x_{42}>2x_{21}-x_{43}$, $2x_{21}-2x_{32}+x_{42}>2x_{11}-2x_{33}+x_{43}$, $2x_{21}-2x_{32}+x_{42}>2x_{11}-2x_{22}+2x_{32}-x_{43}$, and $2x_{21}-2x_{32}+x_{42}>2x_{11}-x_{44}$,  $2x_{11}-2x_{33}+x_{43}\geq 2x_{11}+2x_{32}-2x_{22}-x_{43}$, and $2x_{21}-x_{43}>2x_{11}+2x_{32}-2x_{22}-x_{43}$. These exactly correspond to the third case of $\Tilde{f}_0(\Omega(b))$.
    \item $b_{11}'=b_{11}+1,\,\,\,b_{48}'=b_{48}-1,\,\,\,b_{12}'=b_{12}-1,\,\,\,b_{23}'=b_{23}-1,\,\,\,b_{35}'=b_{35}-1,\,\,\,b_{22}'=b_{22}+1,\,\,\,b_{33}'=b_{33}+1,\,\,\,b_{45}'=b_{45}+1$ if $b_{45}\geq b_{34}$ $b_{12}+b_{35}+b_{36}\geq b_{34}+b_{46}+b_{47}$, $b_{23}+b_{35}\geq b_{34}+b_{46}$, $b_{23}+b_{34}>b_{34}+b_{46}$, $b_{12}+b_{13}+b_{25}>b_{34}+b_{46}+b_{47}$, $b_{12}+b_{24}+b_{36}>b_{34}+b_{46}+b_{47}$, and  $b_{12}+b_{13}+b_{14}>b_{34}+b_{46}+b_{47}$.  
    
    Applying the isomorphism, we end up with the action $x_{11}'=x_{11}+1$,$x_{21}'=x_{21}+1,$, $x_{31}'=x_{31}+1$, and $x_{42}'=x_{42}+1$ along with the following conditions: $x_{41}\leq 2x_{31}-x_{42}$, $2x_{11}-2x_{22}+x_{42}\leq 2x_{31}-x_{42}$, $2x_{21}-2x_{32}+x_{42}\leq 2x_{31}-x_{42}$, $2x_{31}-x_{42}>2x_{21}-x_{43}$, $2x_{31}-x_{42}>2x_{11}-2x_{33}+x_{43}$, $2x_{31}-x_{42}>2x_{11}-2x_{22}+2x_{32}-x_{43}$, and $2x_{31}-x_{42}>2x_{11}-x_{44}$. These exactly correspond to the fourth case of $\Tilde{f}_0(\Omega(b))$.
    \item $b_{11}'=b_{11}+1,\,\,\,b_{48}'=b_{48}-1,\,\,\,b_{12}'=b_{12}-1,\,\,\,b_{25}'=b_{25}-1,\,\,\,b_{36}'=b_{36}-1,\,\,\,b_{22}'=b_{22}+1,\,\,\,b_{35}'=b_{35}+1,\,\,\,b_{46}'=b_{46}+1$ if $b_{45}+b_{46}\geq b_{23}+ b_{24}$ $b_{12}+b_{35}+b_{36}\geq b_{23}+b_{24}+b_{47}$, $b_{35}\geq b_{24}$, $b_{34}+b_{46}\geq b_{23}+b_{24}$, $b_{12}+b_{13}+b_{25}>b_{23}+b_{24}+b_{47}$, $b_{12}+b_{36}>b_{23}+b_{47}$,   $b_{12}+b_{13}+b_{14}>b_{23}+b_{24}+b_{47}$, and $b_{12}+b_{36}>b_{23}+b_{47}$ or $b_{12}+b_{35}+b_{46}>b_{34}+b_{46}+b_{47}$. 
    
    Applying the isomorphism, we end up with the action $x_{11}'=x_{11}+1$,$x_{21}'=x_{21}+1,$, $x_{32}'=x_{32}+1$, and $x_{43}'=x_{43}+1$ along with the following conditions: $x_{41}\leq 2x_{21}-x_{43}$, $2x_{11}-2x_{22}+x_{42}\leq 2x_{21}-x_{43}$, $2x_{21}-2x_{32}+x_{42}\leq 2x_{21}-x_{43}$, $2x_{31}-x_{42}\leq 2x_{21}-x_{43}$, $2x_{21}-x_{43}>2x_{11}-2x_{33}+x_{43}$, $2x_{21}-x_{43}>2x_{11}-2x_{22}+2x_{32}-x_{43}$, $2x_{21}-x_{43}>2x_{11}-x_{44}$, and $2x_{21}-2x_{32}+x_{42}>2x_{11}-2x_{22}+x_{42}$ or $2x_{31}-x_{42}>2x_{11}-2x_{22}+x_{42}$. These exactly correspond to the fifth case of $\Tilde{f}_0(\Omega(b))$.
    \item $b_{11}'=b_{11}+1,\,\,\,b_{48}'=b_{48}-1,\,\,\,b_{14}'=b_{14}-1,\,\,\,b_{26}'=b_{26}-1,\,\,\,b_{37}'=b_{37}-1,\,\,\,b_{24}'=b_{24}+1,\,\,\,b_{36}'=b_{36}+1,\,\,\,b_{47}'=b_{47}+1$ if $b_{45}+b_{46}+b_{47}\geq b_{12}+b_{13}+b_{25}$ $b_{35}+b_{36}\geq b_{13}+b_{25}$, $b_{23}+b_{35}+b_{47}\geq b_{12}+b_{13}+b_{25}$, $b_{34}+b_{46}+b_{47}\geq b_{12}+b_{13}+b_{25}$, $b_{12}+b_{13}+b_{25}\leq b_{23}+b_{24}+b_{47}$, $b_{24}+b_{36}>b_{13}+b_{25}$, and $b_{14}>b_{25}$.
    
    Applying the isomorphism, we end up with the action $x_{11}'=x_{11}+1$,$x_{22}'=x_{22}+1,$, $x_{33}'=x_{33}+1$, and $x_{43}'=x_{43}+1$ along with the following conditions: $x_{41}\leq 2x_{11}-2x_{33}+x_{43}$, $2x_{11}-2x_{22}+x_{42}\leq 2x_{11}-2x_{33}+x_{43}$, $2x_{21}-2x_{32}+x_{42}\leq 2x_{11}-2x_{33}+x_{43}$, $2x_{31}-x_{42}\leq 2x_{11}-2x_{33}+x_{43}$, $2x_{21}-x_{43}\leq 2x_{11}-2x_{33}+x_{43}$, $2x_{11}-2x_{33}+x_{43}>2x_{11}-2x_{22}+2x_{32}-x_{43}$, and $2x_{11}-2x_{33}+x_{43}>2x_{11}-x_{44}$. These exactly correspond to the sixth case of $\Tilde{f}_0(\Omega(b))$.
    \item $b_{11}'=b_{11}+1,\,\,\,b_{48}'=b_{48}-1,\,\,\,b_{13}'=b_{13}-1,\,\,\,b_{25}'=b_{25}-1,\,\,\,b_{37}'=b_{37}-1,\,\,\,b_{23}'=b_{23}+1,\,\,\,b_{35}'=b_{35}+1,\,\,\,b_{47}'=b_{47}+1$ if $b_{45}+b_{46}+b_{47}\geq b_{12}+b_{24}+b_{36}$ $b_{35}\geq b_{24}$, $b_{23}+b_{35}+b_{47}\geq b_{12}+b_{24}+b_{36}$, $b_{34}+b_{46}+b_{47}\geq b_{12}+b_{24}+b_{36}$, $b_{23}+b_{47}\geq b_{12}+b_{36}$, $b_{24}+b_{36}\leq b_{13}+b_{25}$,   $b_{13}+b_{14}>b_{24}+b_{36}$, and $b_{23}+b_{35}>b_{34}+b_{46}$. 
    
    Applying the isomorphism, we end up with the action $x_{11}'=x_{11}+1$,$x_{22}'=x_{22}+1,$, $x_{32}'=x_{32}+1$, and $x_{43}'=x_{43}+1$ along with the following conditions: $x_{41}\leq 2x_{11}+2x_{32}-2x_{22}-x_{43}$, $2x_{11}-2x_{22}+x_{42}\leq 2x_{11}+2x_{32}-2x_{22}-x_{43}$, $2x_{21}-2x_{32}+x_{42}\leq 2x_{11}+2x_{32}-2x_{22}-x_{43}$, $2x_{31}-x_{42}\leq 2x_{11}+2x_{32}-2x_{22}-x_{43}$, $2x_{21}-x_{43}\leq 2x_{11}+2x_{32}-2x_{22}-x_{43}$, $2x_{11}-2x_{33}+x_{43}\leq 2x_{11}-2x_{22}+2x_{32}-x_{43}$,  $2x_{11}+2x_{32}-2x_{22}-x_{43}>2x_{11}-x_{44}$, and $2x_{31}-x_{42}>2x_{21}-2x_{32}+x_{42}$. These exactly correspond to the seventh case of $\Tilde{f}_0(\Omega(b))$.
     \item $b_{11}'=b_{11}+1,\,\,\,b_{48}'=b_{48}-1,\,\,\,b_{15}'=b_{15}-1,\,\,\,b_{26}'=b_{26}-1,\,\,\,b_{37}'=b_{37}-1,\,\,\,b_{25}'=b_{25}+1,\,\,\,b_{36}'=b_{36}+1,\,\,\,b_{47}'=b_{47}+1$ if $b_{45}+b_{46}+b_{47}\geq b_{12}+b_{13}+b_{14}$ $b_{35}+b_{36}\geq b_{13}+b_{14}$, $b_{23}+b_{35}+b_{47}\geq b_{12}+b_{13}+b_{14}$, $b_{34}+b_{46}+b_{47}\geq b_{12}+b_{13}+b_{14}$, $b_{23}+b_{24}+b_{47}\geq b_{12}+b_{13}+b_{14}$, $b_{24}+b_{36}\leq b_{13}+b_{14}$, and  $b_{25}\geq b_{14}$.
    
    Applying the isomorphism, we end up with the action $x_{11}'=x_{11}+1$,$x_{22}'=x_{22}+1,$, $x_{33}'=x_{33}+1$, and $x_{44}'=x_{44}+1$ along with the following conditions: $x_{41}\leq 2x_{11}-x_{44}$, $2x_{11}-2x_{22}+x_{42}\leq 2x_{11}-x_{44}$, $2x_{21}-2x_{32}+x_{42}\leq 2x_{11}-x_{44}$, $2x_{31}-x_{42}\leq 2x_{11}-x_{44}$, $2x_{21}-x_{43}\leq 2x_{11}-x_{44}$, $2x_{11}-2x_{33}+x_{43}\leq 2x_{11}-x_{44}$, and  $2x_{11}+2x_{32}-2x_{22}-x_{43}\leq 2x_{11}-x_{44}$. These exactly correspond to the eighth case of $\Tilde{f}_0(\Omega(b))$.
\end{enumerate}

Now we show that $\Omega(\Tilde{f_1}(b))=\Tilde{f_1}\Omega(b)$. The perfect crystal has $b_{11}'=b_{11}-1$, $b_{12}'=b_{12}+1$, $b_{47}'=b_{47}-1$, and $b_{48}'=b_{48}+1$. Applying the isomorphism, we can see that this corresponds with $x_{11}'=x_{11}-1$ and all other $x_{ij}$ remaining the same.

There are 2 cases for $f_2$:
\begin{enumerate}
    \item $b_{12}'=b_{12}-1$, $b_{36}'=b_{36}-1$, $b_{13}'=b_{13}+1$, and $b_{37}'=b_{37}+1$ if $b_{12}+b_{36}>b_{23}+b_{47}$. Applying the isomorphism, we get that this case corresponds to $x_{22}'=x_{22}-1$ along with the condition $x_{21}+x_{22}>x_{11}+x_{32}$. This corresponds to the first case of $\Tilde{f}_2(\Omega(b))$.
    \item $b_{22}'=b_{22}-1$, $b_{46}'=b_{46}-1$, $b_{23}'=b_{23}+1$, and $b_{47}'=b_{47}+1$ if $b_{12}+b_{36}\leq b_{23}+b_{47}$. Applying the isomorphism, we get that this case corresponds to $x_{21}'=x_{21}-1$ along with the condition $x_{21}+x_{22}\leq x_{11}+x_{32}$. This corresponds to the second case of $\Tilde{f}_2(\Omega(b))$.
\end{enumerate}

There are 3 cases for $f_3$:
\begin{enumerate}
    \item $b_{33}'=b_{33}-1$, $b_{45}'=b_{45}-1$, $b_{34}'=b_{34}+1$, and $b_{46}'=b_{46}+1$ if $b_{25}+b_{35}\leq b_{36}+b_{46}$ and $b_{35}\leq b_{46}$. Applying the isomorphism, we see that these correspond to $x_{31}'=x_{31}-1$ along with the conditions $x_{44}+x_{43}+x_{42}+x_{22}+x_{21}-2x_{33}-2x_{32}-x_{31}\geq x_{44}+x_{43}+x_{22}-2x_{33}-x_{32}$ and $x_{44}+x_{43}+x_{42}+x_{22}+x_{21}-2x_{33}-2x_{32}-x_{31}\geq x_{44}-x_{33}$. 
    \item $b_{23}'=b_{23}-1$, $b_{35}'=b_{35}-1$, $b_{24}'=b_{24}+1$, and $b_{36}'=b_{36}+1$ if $b_{35}>b_{46}$ and $b_{13}\leq b_{24}$ Applying the isomorphism, we see that these correspond to $x_{32}'=x_{32}-1$ along with the conditions $x_{44}+x_{43}+x_{42}+x_{22}+x_{21}-2x_{33}-2x_{32}-x_{31}< x_{44}+x_{43}+x_{22}-2x_{33}-x_{32}$ and $x_{44}+x_{43}+x_{22}-2x_{33}-x_{32}\geq x_{44}-x_{33}$. 
    \item $b_{13}'=b_{13}-1$, $b_{25}'=b_{25}-1$, $b_{14}'=b_{14}+1$, and $b_{26}'=b_{26}+1$ if $b_{25}+b_{35}>b_{36}+b_{46}$ and $b_{13}> b_{24}$ Applying the isomorphism, we see that these correspond to $x_{33}'=x_{33}-1$ along with the conditions $x_{44}+x_{43}+x_{42}+x_{22}+x_{21}-2x_{33}-2x_{32}-x_{31}> x_{44}+x_{43}+x_{22}-2x_{33}-x_{32}$ and $x_{44}+x_{43}+x_{42}+x_{22}+x_{21}-2x_{33}-2x_{32}-x_{31}> x_{44}-x_{33}$. 
\end{enumerate}
These all correspond to the cases of $\Tilde{f}_3(\Omega(b))$

Finally we consider $f_4$, which has 4 cases:
\begin{enumerate}
    \item $b_{14}'=b_{14}-1$, $b_{15}'=b_{15}+1$ if $b_{14}>b_{25}$, $b_{14}+b_{24}>b_{25}+b_{35}$, and $b_{14}+b_{24}+b_{34}>b_{25}+b_{35}+b_{45}$. Applying the isomorphism, we see that this corresponds to $x_{44}'=x_{44}-1$ along with the conditions $-x_{44}>2x_{33}-2x_{44}-x_{43}$, $-x_{44}>2x_{32}+2x_{33}-2x_{44}-2x_{43}-x_{42}$, and $-x_{44}>2x_{31}+2x_{32}+2x_{33}-2x_{44}-2x_{43}-2x_{42}-x_{41}$
    \item $b_{24}'=b_{24}-1$, $b_{25}'=b_{25}+1$ if $b_{14}\leq b_{25}$, $b_{24}>b_{35}$, and $b_{24}+b_{34}>b_{35}+b_{45}$. Applying the isomorphism, we see that this corresponds to $x_{43}'=x_{43}-1$ along with the conditions $-x_{44}\leq 2x_{33}-2x_{44}-x_{43}$, $2x_{33}-2x_{44}-x_{43}>2x_{32}+2x_{33}-2x_{44}-2x_{43}-x_{42}$, and $2x_{33}-2x_{44}-x_{43}>2x_{31}+2x_{32}+2x_{33}-2x_{44}-2x_{43}-2x_{42}-x_{41}$
    \item $b_{34}'=b_{34}-1$, $b_{35}'=b_{35}+1$ if $b_{14}+b_{24}\leq b_{25}+b_{35}$, $b_{24}\leq b_{35}$, and $b_{34}>b_{45}$. Applying the isomorphism, we see that this corresponds to $x_{42}'=x_{42}-1$ along with the conditions $-x_{44}\leq 2x_{32}+2x_{33}-2x_{44}-2x_{43}-x_{42}$, $2x_{33}-2x_{44}-x_{43}\leq 2x_{32}+2x_{33}-2x_{44}-2x_{43}-x_{42}$, and $2x_{32}+2x_{33}-2x_{44}-2x_{43}-x_{42}>2x_{31}+2x_{32}+2x_{33}-2x_{44}-2x_{43}-2x_{42}-x_{41}$
    \item $b_{44}'=b_{44}-1$, $b_{45}'=b_{45}+1$ if $b_{14}+b_{24}+b_{34}\leq b_{25}+b_{35}+b_{45}$, $b_{24}+b_{34}\leq b_{35}+b_{45}$, and $b_{34}\leq b_{45}$. Applying the isomorphism, we see that this corresponds to $x_{41}'=x_{41}-1$ along with the conditions $-x_{44}\leq 2x_{31}+2x_{32}+2x_{33}-2x_{44}-2x_{43}-2x_{42}-x_{41}$, $2x_{33}-2x_{44}-x_{43}\leq 2x_{31}+2x_{32}+2x_{33}-2x_{44}-2x_{43}-2x_{42}-x_{41}$, and $2x_{32}+2x_{33}-2x_{44}-2x_{43}-x_{42}\leq 2x_{31}+2x_{32}+2x_{33}-2x_{44}-2x_{43}-2x_{42}-x_{41}$
\end{enumerate}
This corresponds to all of the cases of $\Tilde{f}_4(\Omega(b))$. Thus 1. is proved.

Next we check 2. There are 8 cases for $\Tilde{\varepsilon_0}(b)$.
\begin{enumerate}
    \item $b_{11}'=b_{11}-1,\,\,\,b_{48}'=b_{48}+1,\,\,\,b_{12}'=b_{12}+1,\,\,\,b_{23}'=b_{23}+1,\,\,\,b_{34}'=b_{34}+1,\,\,\,b_{22}'=b_{22}-1,\,\,\,b_{33}'=b_{33}-1,\,\,\,b_{44}'=b_{44}-1$ if $b_{12}+b_{35}+b_{36}\geq b_{45}+b_{46}+b_{47}$, $b_{23}+b_{35}\geq b_{45}+b_{46}$, $b_{34}\geq b_{45}$, $b_{23}+b_{34}\geq b_{45}+b_{46}$, $b_{12}+b_{13}+b_{25}\geq b_{45}+b_{46}+b_{47}$, $b_{12}+b_{24}+b_{36}\geq b_{45}+b_{46}+b_{47}$, and $b_{12}+b_{13}+b_{14}\geq b_{45}+b_{46}+b_{47}$. 
    
    Applying the isomorphism, we end up with the action $x_{11}'=x_{11}-1$,$x_{21}'=x_{21}-1,$, $x_{31}'=x_{31}-1$, and $x_{41}'=x_{41}-1$ along with the following conditions: $x_{41}\geq 2x_{11}-2x_{22}+x_{42}$, $x_{41}\geq 2x_{21}-2x_{32}+x_{42}$, $x_{41}\geq 2x_{31}-x_{42}$, $x_{41}\geq 2x_{21}-x_{43}$, $x_{41}\geq 2x_{11}-2x_{33}+x_{43}$, $x_{41}\geq 2x_{11}-2x_{22}+2x_{32}-x_{43}$, and $x_{41}\geq 2x_{11}-x_{44}$. These exactly correspond to the first case of $\Tilde{e}_0(\Omega(b))$.
    \item $b_{11}'=b_{11}-1,\,\,\,b_{48}'=b_{48}+1,\,\,\,b_{13}'=b_{13}+1,\,\,\,b_{24}'=b_{24}+1,\,\,\,b_{37}'=b_{37}+1,\,\,\,b_{23}'=b_{23}-1,\,\,\,b_{34}'=b_{34}-1,\,\,\,b_{47}'=b_{47}-1$ if $b_{12}+b_{35}+b_{36}< b_{45}+b_{46}+b_{47}$, $b_{23}+b_{47}\geq b_{12}+b_{36}$, $b_{34}+b_{46}+b_{47}\geq b_{12}+b_{35}+b_{36}$, $b_{23}+b_{34}+b_{47}\geq b_{12}+b_{35}+b_{36}$, $b_{13}+b_{25}\geq b_{35}+b_{36}$, $b_{24}\geq b_{35}$, $b_{13}+b_{14}\geq b_{35}+b_{36}$, and $b_{23}+b_{24}+b_{47}> b_{12}+b_{13}+b_{25}$ or $b_{23}+b_{47}> b_{12}+b_{36}$ . 
    
    Applying the isomorphism, we end up with the action $x_{11}'=x_{11}-1$,$x_{22}'=x_{22}-1,$, $x_{32}'=x_{32}-1$, and $x_{42}'=x_{42}-1$ along with the following conditions: $x_{41}< 2x_{11}-2x_{22}+x_{42}$, $2x_{11}-2x_{22}+x_{42}\geq 2x_{21}-2x_{32}+x_{42}$, $2x_{11}-2x_{22}+x_{42}\geq 2x_{31}-x_{42}$, $2x_{11}-2x_{22}+x_{42}\geq 2x_{21}-x_{43}$, $2x_{11}-2x_{22}+x_{42}\geq 2x_{11}-2x_{33}+x_{43}$, $2x_{11}-2x_{22}+x_{42}\geq 2x_{11}-2x_{22}+2x_{32}-x_{43}$, and $2x_{11}-2x_{22}+x_{42}\geq 2x_{11}-x_{44}$, and $2x_{11}-2x_{33}+x_{43}> 2x_{21}-x_{43}$ or $2x_{11}+2x_{32}-2x_{22}-x_{43}> 2x_{21}-x_{43}$. These exactly correspond to the second case of $\Tilde{e}_0(\Omega(b))$.
    \item $b_{11}'=b_{11}-1,\,\,\,b_{48}'=b_{48}+1,\,\,\,b_{12}'=b_{12}+1,\,\,\,b_{24}'=b_{24}+1,\,\,\,b_{36}'=b_{36}+1,\,\,\,b_{22}'=b_{22}-1,\,\,\,b_{34}'=b_{34}-1,\,\,\,b_{46}'=b_{46}-1$ if $b_{45}+b_{46}> b_{23}+b_{35}$, $b_{12}+b_{36}\geq b_{23}+b_{47}$, $b_{34}+b_{46}\geq b_{23}+b_{35}$, $b_{24}\geq b_{35}$, $b_{12}+b_{13}+b_{25}\geq b_{23}+b_{35}+b_{47}$,$b_{12}+b_{24}+b_{36}\geq b_{23}+b_{35}+b_{47}$,   $b_{12}+b_{13}+b_{14}\geq b_{23}+b_{35}+b_{47}$, $b_{24}+b_{36}\geq b_{13}+b_{25}$, and $b_{12}+b_{36}\geq b_{23}+b_{47}$. 
    
    Applying the isomorphism, we end up with the action $x_{11}'=x_{11}-1$,$x_{21}'=x_{21}-1,$, $x_{32}'=x_{32}-1$, and $x_{42}'=x_{42}-1$ along with the following conditions: $x_{41}< 2x_{21}-2x_{32}+x_{42}$, $2x_{11}-2x_{22}+x_{42}< 2x_{21}-2x_{32}+x_{42}$, $2x_{21}-2x_{32}+x_{42}\geq 2x_{31}-x_{42}$, $2x_{21}-2x_{32}+x_{42}\geq 2x_{21}-x_{43}$, $2x_{21}-2x_{32}+x_{42}\geq 2x_{11}-2x_{33}+x_{43}$, $2x_{21}-2x_{32}+x_{42}\geq 2x_{11}-2x_{22}+2x_{32}-x_{43}$, and $2x_{21}-2x_{32}+x_{42}\geq 2x_{11}-x_{44}$,  $2x_{11}-2x_{33}+x_{43}> 2x_{11}+2x_{32}-2x_{22}-x_{43}$, and $2x_{21}-x_{43}\geq 2x_{11}+2x_{32}-2x_{22}-x_{43}$. These exactly correspond to the third case of $\Tilde{e}_0(\Omega(b))$.
    \item $b_{11}'=b_{11}-1,\,\,\,b_{48}'=b_{48}+1,\,\,\,b_{12}'=b_{12}+1,\,\,\,b_{23}'=b_{23}+1,\,\,\,b_{35}'=b_{35}+1,\,\,\,b_{22}'=b_{22}-1,\,\,\,b_{33}'=b_{33}-1,\,\,\,b_{45}'=b_{45}-1$ if $b_{45}> b_{34}$ $b_{12}+b_{35}+b_{36}> b_{34}+b_{46}+b_{47}$, $b_{23}+b_{35}> b_{34}+b_{46}$, $b_{23}+b_{34}\geq b_{34}+b_{46}$, $b_{12}+b_{13}+b_{25}\geq b_{34}+b_{46}+b_{47}$, $b_{12}+b_{24}+b_{36}\geq b_{34}+b_{46}+b_{47}$, and  $b_{12}+b_{13}+b_{14}\geq b_{34}+b_{46}+b_{47}$.  
    
    Applying the isomorphism, we end up with the action $x_{11}'=x_{11}-1$,$x_{21}'=x_{21}-1,$, $x_{31}'=x_{31}-1$, and $x_{42}'=x_{42}-1$ along with the following conditions: $x_{41}< 2x_{31}-x_{42}$, $2x_{11}-2x_{22}+x_{42}< 2x_{31}-x_{42}$, $2x_{21}-2x_{32}+x_{42}< 2x_{31}-x_{42}$, $2x_{31}-x_{42}\geq 2x_{21}-x_{43}$, $2x_{31}-x_{42}\geq 2x_{11}-2x_{33}+x_{43}$, $2x_{31}-x_{42}\geq 2x_{11}-2x_{22}+2x_{32}-x_{43}$, and $2x_{31}-x_{42}\geq 2x_{11}-x_{44}$. These exactly correspond to the fourth case of $\Tilde{e}_0(\Omega(b))$.
    \item $b_{11}'=b_{11}-1,\,\,\,b_{48}'=b_{48}+1,\,\,\,b_{12}'=b_{12}+1,\,\,\,b_{25}'=b_{25}+1,\,\,\,b_{36}'=b_{36}+1,\,\,\,b_{22}'=b_{22}-1,\,\,\,b_{35}'=b_{35}-1,\,\,\,b_{46}'=b_{46}-1$ if $b_{45}+b_{46}>b_{23}+ b_{24}$ $b_{12}+b_{35}+b_{36}> b_{23}+b_{24}+b_{47}$, $b_{35}> b_{24}$, $b_{34}+b_{46}> b_{23}+b_{24}$, $b_{12}+b_{13}+b_{25}\geq b_{23}+b_{24}+b_{47}$, $b_{12}+b_{36}\geq b_{23}+b_{47}$,   $b_{12}+b_{13}+b_{14}\geq b_{23}+b_{24}+b_{47}$, and $b_{12}+b_{36}\geq b_{23}+b_{47}$ or $b_{12}+b_{35}+b_{46}\geq b_{34}+b_{46}+b_{47}$. 
    
    Applying the isomorphism, we end up with the action $x_{11}'=x_{11}-1$,$x_{21}'=x_{21}-1,$, $x_{32}'=x_{32}-1$, and $x_{43}'=x_{43}-1$ along with the following conditions: $x_{41}< 2x_{21}-x_{43}$, $2x_{11}-2x_{22}+x_{42}< 2x_{21}-x_{43}$, $2x_{21}-2x_{32}+x_{42}< 2x_{21}-x_{43}$, $2x_{31}-x_{42}< 2x_{21}-x_{43}$, $2x_{21}-x_{43}\geq 2x_{11}-2x_{33}+x_{43}$, $2x_{21}-x_{43}
    \geq 2x_{11}-2x_{22}+2x_{32}-x_{43}$, $2x_{21}-x_{43}\geq 2x_{11}-x_{44}$, and $2x_{21}-2x_{32}+x_{42}\geq 2x_{11}-2x_{22}+x_{42}$ or $2x_{31}-x_{42}\geq 2x_{11}-2x_{22}+x_{42}$. These exactly correspond to the fifth case of $\Tilde{e}_0(\Omega(b))$.
    \item $b_{11}'=b_{11}-1,\,\,\,b_{48}'=b_{48}+1,\,\,\,b_{14}'=b_{14}+1,\,\,\,b_{26}'=b_{26}+1,\,\,\,b_{37}'=b_{37}+1,\,\,\,b_{24}'=b_{24}-1,\,\,\,b_{36}'=b_{36}-1,\,\,\,b_{47}'=b_{47}-1$ if $b_{45}+b_{46}+b_{47}> b_{12}+b_{13}+b_{25}$ $b_{35}+b_{36}> b_{13}+b_{25}$, $b_{23}+b_{35}+b_{47}> b_{12}+b_{13}+b_{25}$, $b_{34}+b_{46}+b_{47}> b_{12}+b_{13}+b_{25}$, $b_{12}+b_{13}+b_{25}< b_{23}+b_{24}+b_{47}$, $b_{24}+b_{36}\geq b_{13}+b_{25}$,  and $b_{14}\geq b_{25}$. 
    
    Applying the isomorphism, we end up with the action $x_{11}'=x_{11}-1$,$x_{22}'=x_{22}-1,$, $x_{33}'=x_{33}-1$, and $x_{43}'=x_{43}-1$ along with the following conditions: $x_{41}< 2x_{11}-2x_{33}+x_{43}$, $2x_{11}-2x_{22}+x_{42}< 2x_{11}-2x_{33}+x_{43}$, $2x_{21}-2x_{32}+x_{42}< 2x_{11}-2x_{33}+x_{43}$, $2x_{31}-x_{42}< 2x_{11}-2x_{33}+x_{43}$, $2x_{21}-x_{43}< 2x_{11}-2x_{33}+x_{43}$, $2x_{11}-2x_{33}+x_{43}\geq 2x_{11}-2x_{22}+2x_{32}-x_{43}$, and $2x_{11}-2x_{33}+x_{43}\geq 2x_{11}-x_{44}$. These exactly correspond to the sixth case of $\Tilde{e}_0(\Omega(b))$.
    \item $b_{11}'=b_{11}-1,\,\,\,b_{48}'=b_{48}+1,\,\,\,b_{13}'=b_{13}+1,\,\,\,b_{25}'=b_{25}+1,\,\,\,b_{37}'=b_{37}+1,\,\,\,b_{23}'=b_{23}-1,\,\,\,b_{35}'=b_{35}-1,\,\,\,b_{47}'=b_{47}-1$ if $b_{45}+b_{46}+b_{47}> b_{12}+b_{24}+b_{36}$ $b_{35}> b_{24}$, $b_{23}+b_{35}+b_{47}> b_{12}+b_{24}+b_{36}$, $b_{34}+b_{46}+b_{47}> b_{12}+b_{24}+b_{36}$, $b_{23}+b_{47}> b_{12}+b_{36}$, $b_{24}+b_{36}< b_{13}+b_{25}$,   $b_{13}+b_{14}\geq b_{24}+b_{36}$, and $b_{23}+b_{35}\geq b_{34}+b_{46}$. 
    
    Applying the isomorphism, we end up with the action $x_{11}'=x_{11}-1$,$x_{22}'=x_{22}-1,$, $x_{32}'=x_{32}-1$, and $x_{43}'=x_{43}-1$ along with the following conditions: $x_{41}< 2x_{11}+2x_{32}-2x_{22}-x_{43}$, $2x_{11}-2x_{22}+x_{42}< 2x_{11}+2x_{32}-2x_{22}-x_{43}$, $2x_{21}-2x_{32}+x_{42}< 2x_{11}+2x_{32}-2x_{22}-x_{43}$, $2x_{31}-x_{42}< 2x_{11}+2x_{32}-2x_{22}-x_{43}$, $2x_{21}-x_{43}< 2x_{11}+2x_{32}-2x_{22}-x_{43}$, $2x_{11}-2x_{33}+x_{43}< 2x_{11}-2x_{22}+2x_{32}-x_{43}$,  $2x_{11}+2x_{32}-2x_{22}-x_{43}\geq 2x_{11}-x_{44}$, and $2x_{31}-x_{42}\geq 2x_{21}-2x_{32}+x_{42}$. These exactly correspond to the seventh case of $\Tilde{e}_0(\Omega(b))$.
     \item $b_{11}'=b_{11}-1,\,\,\,b_{48}'=b_{48}+1,\,\,\,b_{15}'=b_{15}+1,\,\,\,b_{26}'=b_{26}+1,\,\,\,b_{37}'=b_{37}+1,\,\,\,b_{25}'=b_{25}-1,\,\,\,b_{36}'=b_{36}-1,\,\,\,b_{47}'=b_{47}-1$ if $b_{45}+b_{46}+b_{47}> b_{12}+b_{13}+b_{14}$ $b_{35}+b_{36}> b_{13}+b_{14}$, $b_{23}+b_{35}+b_{47}> b_{12}+b_{13}+b_{14}$, $b_{34}+b_{46}+b_{47}> b_{12}+b_{13}+b_{14}$, $b_{23}+b_{24}+b_{47}> b_{12}+b_{13}+b_{14}$, $b_{24}+b_{36}> b_{13}+b_{14}$, and  $b_{25}> b_{14}$.
    
    Applying the isomorphism, we end up with the action $x_{11}'=x_{11}-1$,$x_{22}'=x_{22}-1,$, $x_{33}'=x_{33}-1$, and $x_{44}'=x_{44}-1$ along with the following conditions: $x_{41}< 2x_{11}-x_{44}$, $2x_{11}-2x_{22}+x_{42}< 2x_{11}-x_{44}$, $2x_{21}-2x_{32}+x_{42}< 2x_{11}-x_{44}$, $2x_{31}-x_{42}< 2x_{11}-x_{44}$, $2x_{21}-x_{43}< 2x_{11}-x_{44}$, $2x_{11}-2x_{33}+x_{43}< 2x_{11}-x_{44}$, and  $2x_{11}+2x_{32}-2x_{22}-x_{43}< 2x_{11}-x_{44}$. These exactly correspond to the eighth case of $\Tilde{e}_0(\Omega(b))$.
\end{enumerate}

Now we show that $\Omega(\Tilde{e_1}(b))=\Tilde{e_1}\Omega(b)$. The perfect crystal has $b_{11}'=b_{11}+1$, $b_{12}'=b_{12}-1$, $b_{47}'=b_{47}+1$, and $b_{48}'=b_{48}-1$. Applying the isomorphism, we can see that this corresponds with $x_{11}'=x_{11}+1$ and all other $x_{ij}$ remaining the same.

There are 2 cases for $\Tilde{e}_2$:
\begin{enumerate}
    \item $b_{12}'=b_{12}+1$, $b_{36}'=b_{36}+1$, $b_{13}'=b_{13}-1$, and $b_{37}'=b_{37}-1$ if $b_{12}+b_{36}\geq b_{23}+b_{47}$. Applying the isomorphism, we get that this case corresponds to $x_{22}'=x_{22}+1$ along with the condition $x_{21}+x_{22}\geq x_{11}+x_{32}$. This corresponds to the first case of $\Tilde{e}_2(\Omega(b))$.
    \item $b_{22}'=b_{22}+1$, $b_{46}'=b_{46}+1$, $b_{23}'=b_{23}-1$, and $b_{47}'=b_{47}-1$ if $b_{12}+b_{36}< b_{23}+b_{47}$. Applying the isomorphism, we get that this case corresponds to $x_{21}'=x_{21}+1$ along with the condition $x_{21}+x_{22}< x_{11}+x_{32}$. This corresponds to the second case of $\Tilde{e}_2(\Omega(b))$.
\end{enumerate}

There are 3 cases for $\Tilde{e}_3$:
\begin{enumerate}
    \item $b_{33}'=b_{33}+1$, $b_{45}'=b_{45}+1$, $b_{34}'=b_{34}-1$, and $b_{46}'=b_{46}-1$ if $b_{25}+b_{35}< b_{36}+b_{46}$ and $b_{35}<b_{46}$. Applying the isomorphism, we see that these correspond to $x_{31}'=x_{31}+1$ along with the conditions $x_{44}+x_{43}+x_{42}+x_{22}+x_{21}-2x_{33}-2x_{32}-x_{31}> x_{44}+x_{43}+x_{22}-2x_{33}-x_{32}$ and $x_{44}+x_{43}+x_{42}+x_{22}+x_{21}-2x_{33}-2x_{32}-x_{31}> x_{44}-x_{33}$. 
    \item $b_{23}'=b_{23}+1$, $b_{35}'=b_{35}+1$, $b_{24}'=b_{24}-1$, and $b_{36}'=b_{36}-1$ if $b_{35}\geq b_{46}$ and $b_{13}< b_{24}$. Applying the isomorphism, we see that these correspond to $x_{32}'=x_{32}+1$ along with the conditions $x_{44}+x_{43}+x_{42}+x_{22}+x_{21}-2x_{33}-2x_{32}-x_{31}\leq x_{44}+x_{43}+x_{22}-2x_{33}-x_{32}$ and $x_{44}+x_{43}+x_{22}-2x_{33}-x_{32}> x_{44}-x_{33}$. 
    \item $b_{13}'=b_{13}+1$, $b_{25}'=b_{25}+1$, $b_{14}'=b_{14}-1$, and $b_{26}'=b_{26}-1$ if $b_{25}+b_{35}\geq b_{36}+b_{46}$ and $b_{13}\geq b_{24}$ Applying the isomorphism, we see that these correspond to $x_{33}'=x_{33}+1$ along with the conditions $x_{44}+x_{43}+x_{42}+x_{22}+x_{21}-2x_{33}-2x_{32}-x_{31}\geq  x_{44}+x_{43}+x_{22}-2x_{33}-x_{32}$ and $x_{44}+x_{43}+x_{42}+x_{22}+x_{21}-2x_{33}-2x_{32}-x_{31}\geq  x_{44}-x_{33}$. 
\end{enumerate}
These all correspond to the cases of $\Tilde{e}_3(\Omega(b))$

Finally we consider $\Tilde{e}_4$, which has 4 cases:
\begin{enumerate}
    \item $b_{14}'=b_{14}+1$, $b_{15}'=b_{15}-1$ if $b_{14}\geq b_{25}$, $b_{14}+b_{24}\geq b_{25}+b_{35}$, and $b_{14}+b_{24}+b_{34}\geq b_{25}+b_{35}+b_{45}$. Applying the isomorphism, we see that this corresponds to $x_{44}'=x_{44}+1$ along with the conditions $-x_{44}\geq 2x_{33}-2x_{44}-x_{43}$, $-x_{44}\geq 2x_{32}+2x_{33}-2x_{44}-2x_{43}-x_{42}$, and $-x_{44}\geq 2x_{31}+2x_{32}+2x_{33}-2x_{44}-2x_{43}-2x_{42}-x_{41}$
    \item $b_{24}'=b_{24}+1$, $b_{25}'=b_{25}-1$ if $b_{14}< b_{25}$, $b_{24}\geq b_{35}$, and $b_{24}+b_{34}\geq b_{35}+b_{45}$. Applying the isomorphism, we see that this corresponds to $x_{43}'=x_{43}+1$ along with the conditions $-x_{44}< 2x_{33}-2x_{44}-x_{43}$, $2x_{33}-2x_{44}-x_{43}\geq 2x_{32}+2x_{33}-2x_{44}-2x_{43}-x_{42}$, and $2x_{33}-2x_{44}-x_{43}\geq 2x_{31}+2x_{32}+2x_{33}-2x_{44}-2x_{43}-2x_{42}-x_{41}$
    \item $b_{34}'=b_{34}+1$, $b_{35}'=b_{35}-1$ if $b_{14}+b_{24}< b_{25}+b_{35}$, $b_{24}< b_{35}$, and $b_{34}\geq b_{45}$. Applying the isomorphism, we see that this corresponds to $x_{42}'=x_{42}+1$ along with the conditions $-x_{44}< 2x_{32}+2x_{33}-2x_{44}-2x_{43}-x_{42}$, $2x_{33}-2x_{44}-x_{43}< 2x_{32}+2x_{33}-2x_{44}-2x_{43}-x_{42}$, and $2x_{32}+2x_{33}-2x_{44}-2x_{43}-x_{42}\geq 2x_{31}+2x_{32}+2x_{33}-2x_{44}-2x_{43}-2x_{42}-x_{41}$
    \item $b_{44}'=b_{44}+1$, $b_{45}'=b_{45}-1$ if $b_{14}+b_{24}+b_{34}< b_{25}+b_{35}+b_{45}$, $b_{24}+b_{34}< b_{35}+b_{45}$, and $b_{34}< b_{45}$. Applying the isomorphism, we see that this corresponds to $x_{41}'=x_{41}+1$ along with the conditions $-x_{44}< 2x_{31}+2x_{32}+2x_{33}-2x_{44}-2x_{43}-2x_{42}-x_{41}$, $2x_{33}-2x_{44}-x_{43}< 2x_{31}+2x_{32}+2x_{33}-2x_{44}-2x_{43}-2x_{42}-x_{41}$, and $2x_{32}+2x_{33}-2x_{44}-2x_{43}-x_{42}< 2x_{31}+2x_{32}+2x_{33}-2x_{44}-2x_{43}-2x_{42}-x_{41}$
\end{enumerate}
This corresponds to all of the cases of $\Tilde{e}_4(\Omega(b))$. Thus 2. is proved.
Now we prove $\text{wt}_k(\Omega(b))=\text{wt}_k(b)$ for $k=0,1,2,3,4$.
\begin{align*}
   &\text{wt}_0(b)=b_{48}-b_{11}=-2x_{11}=\text{wt}_0(\Omega(b)) \\
   &\text{wt}_1(b)=b_{47}-b_{12}=2x_{11}-x_{21}-x_{22}=\text{wt}_1(\Omega(b))\\
   &\text{wt}_2(b)=b_{12}-b_{23}+b_{46}-b_{13}=2x_{21}+2x_{22}-x_{11}-x_{31}-x_{32}-x_{33}=\text{wt}_2(\Omega(b))\\
   &\text{wt}_3(b)=b_{45}-b_{14}+b_{13}-b_{24}+b_{23}-b_{34}\\&=2x_{31}+2x_{32}+2x_{33}-x_{21}-x_{22}-x_{41}-x_{42}-x_{43}-x_{44}=\text{wt}_3(\Omega(b))\\
   &\text{wt}_4(b)=b_{44}-b_{15}+b_{34}-b_{45}+b_{24}-b_{35}+b_{14}-b_{25}\\&=2x_{41}+2x_{42}+2x_{43}+2x_{44}-2x_{31}-2x_{32}-2x_{33}=\text{wt}_4(\Omega(b))
\end{align*}
This proves 3.
Finally, we need to show $\varepsilon_k(\Omega(b))=\varepsilon_k(b)$.

For $k=0$, we have
\begin{align*}
    \varepsilon_k(b)=\begin{cases}
    -b_{48}-b_{47}-b_{46}-b_{45}& 1<2,3,4,5,6,7,8\\
    -b_{48}-b_{12}-b_{35}-b_{36}& 2\geq 1, 2<3,4,5,6,7,8\\
    -b_{48}-b_{47}-b_{35}-b_{23}&3\geq 1,2, 3<4,5,6,7,8\\
    -b_{48}-b_{47}-b_{46}-b_{34}& 4\geq 1,2,3, 4<5,6,7,8\\
    -b_{48}-b_{47}-b_{24}-b_{23}& 5\geq 1,2,3,4, 5<6,7,8\\
    -b_{48}-b_{25}-b_{13}-b_{12}& 6\geq 1,2,3,4,5, 6<7,8\\
    -b_{48}-b_{36}-b_{24}-b_{12}& 7\geq 1,2,3,4,5,6, 7<8\\
    -b_{48}-b_{14}-b_{13}-b_{12}& 8\geq 1,2,3,4,5,6,7
    \end{cases}
\end{align*}
 where $1=b_{45}+b_{46}+b_{47}$, $2=b_{12}+b_{35}+b_{36}$, $3=b_{23}+b_{35}+b_{47}$, $4=b_{34}+b_{46}+b_{47}$, $5=b_{23}+b_{24}+b_{47}$, $6=b_{12}+b_{13}+b_{25}$, $7=b_{12}+b_{24}+b_{36}$, and $8=b_{12}+b_{13}=b_{14}$.
Applying the isomorphism, we get the following, which is the same as $\varepsilon_0(\Omega(b))$
\begin{align*}
    \varepsilon_0(\Omega(b))=\begin{cases}
x_{41} & x_{41}>A,B,C,X,Y,Z,2x_{11}-x_{44}\\
2x_{11}-2x_{22}+x_{42}&A\geq x_{41}, A>B,C,X,Y,Z, 2x_{11}-x_{44}\\
2x_{21}-2x_{32}+x_{42}&B\geq x_{41},A, B>C,X,Y,Z, 2x_{11}-x_{44}\\
2x_{31}-x_{42}&C\geq x_{41},A,B, C> X,Y,Z, 2x_{11}-x_{44}\\
2x_{21}-x_{43}&X\geq x_{41},A,B,C, X>Y,Z, 2x_{11}-x_{44} \\
2x_{11}-2x_{33}+x_{43}&Y\geq x_{41},A,B,C,X, Y>Z, 2x_{11}-x_{44}\\
2x_{11}+2x_{32}-2x_{22}-x_{43}&Z\geq x_{41},A,B,C,X,Y, Z>2x_{11}-x_{44}\\
2x_{11}-x_{44}& 2x_{11}-x_{44}\geq x_{41},A,B,C,X,Y,Z
\end{cases}
\end{align*}
where $A=2x_{11}-2x_{22}+x_{42}$, $B=2x_{21}-2x_{32}+x_{42}$, $C=2x_{31}-x_{42}$, $X=2x_{21}-x_{43}$, $Y=2x_{11}-2x_{33}+x_{43}$, and $Z=2x_{11}+2x_{32}-2x_{22}-x_{43}$.

For $k=1$, $\varepsilon_1(b)=b_{12}$. Applying the isomorphism, we get $x_{22}-x_{11}$ which is equal to $\varepsilon_1(\Omega(b))$.

For $k=2$, \begin{align*}
    \varepsilon_2(b)=\begin{cases}
    b_{13} & b_{23}-b_{12}<0\\
    b_{13}+b_{23}-b_{12}& b_{23}-b_{12}\geq 0
    \end{cases}
\end{align*}

Applying the isomorphism, we get exactly $\varepsilon_2(\Omega(b))$ as follows: 
\begin{align*}
    \varepsilon_2(\Omega(b))=\begin{cases}
    x_{33}-x_{22} & x_{33}-x_{22}>x_{33}+x_{32}+x_{11}-2x_{22}-x_{21}\\
    x_{33}+x_{32}+x_{11}-2x_{22}-x_{21}& x_{33}-x_{22}\leq x_{33}+x_{32}+x_{11}-2x_{22}-x_{21}
    \end{cases}
\end{align*}

For $k=3$, \begin{align*}
    \varepsilon_3(b)=\begin{cases}
    b_{14} & b_{24}-b_{13}<0,\,\,\, b_{24}+b_{34}-b_{13}-b_{23}<0\\
    b_{14}+b_{24}-b_{13}& b_{23}-b_{12}\geq 0,\,\,\,b_{34}-b_{23}<0\\
    b_{14}+b_{24}-b_{13}+b_{34}-b_{23} & b_{34}-b_{23}\geq 0,\,\,\, b_{24}+b_{34}-b_{13}-b_{23}\geq 0
    \end{cases}
\end{align*}

Applying the isomorphism, we get exactly $\varepsilon_3(\Omega(b))$ as follows: 
\begin{align*}
    \varepsilon_3(\Omega(b))=\begin{cases}
    x_{44}-x_{33} & D>E,F \\
    x_{44}+x_{43}+x_{22}-2x_{33}-x_{32}& D\leq E, E>F\\
    x_{44}+x_{43}+x_{42}+x_{21}+x_{22}-2x_{33}-2x_{32}-x_{31}& F\geq D,E
    \end{cases}
\end{align*}
where $D=x_{44}-x_{33}$, $E=x_{44}+x_{43}+x_{22}-2x_{33}-x_{32}$, and $F=x_{44}+x_{43}+x_{42}+x_{21}+x_{22}-2x_{33}-2x_{32}-x_{31}$.

Finally we check the relation for $k=4$.

\begin{align*}
    \varepsilon_4(b)=\begin{cases}
    b_{15} & (1)>(2),(3),(4)\\
    b_{15}+b_{22}-b_{14}& (2)\geq (1), (2)>(3),(4)\\
    b_{15}+b_{25}-b_{14}+b_{35}-b_{24} & (3)\geq (2),(1), (3)>(4)\\
     b_{15}+b_{25}-b_{14}+b_{35}-b_{24}+b_{45}-b_{34} & (4)\geq (3),(2),(1)\\
    \end{cases}
\end{align*}
where $(1)=b_{15}$, $(2)=b_{15}+b_{25}-b_{14}$, $(3)=b_{15}+b_{25}-b_{14}+b_{35}-b_{24}$, and $(4)= b_{15}+b_{25}-b_{14}+b_{35}-b_{24}+b_{45}-b_{34}$.

Applying the isomorphism, we get exactly $\varepsilon_4(\Omega(b))$ as follows: 
\begin{align*}
    \varepsilon_4(\Omega(b))=\begin{cases}
    -x_{44} & A>B,C,D \\
    2x_{33}-2x_{44}-x_{43}& B\geq A, B>C,D\\
   2x_{33}+2x_{32}-2x_{44}-2x_{43}-x_{42}& C\geq A,B, C>D\\
   2x_{33}+2x_{32}+2x_{31}-2x_{44}-2x_{43}-2x_{42}-x_{41}& D\geq A,B,C
    \end{cases}
\end{align*}
where $A=-x_{44}$, $B=2x_{33}-2x_{44}-x_{43}$, and $C= 2x_{33}+2x_{32}-2x_{44}-2x_{43}-x_{42}$, $D=2x_{33}+2x_{32}+2x_{31}-2x_{44}-2x_{43}-2x_{42}-x_{41}$.

This proves 4. Therefore, $\Omega$ is an isomorphism.
\end{proof}

%\bibliographystyle{amsplain}
\bibliographystyle{amsalpha}

\end{document}